\documentclass[a4paper, 12pt]{amsart}

\usepackage[font=small,labelfont=bf]{caption}
 \usepackage{amsmath, amsthm, amssymb, amsfonts}
\usepackage[utf8]{inputenc}
\usepackage[normalem]{ulem}
\usepackage[colorlinks = true,
            linkcolor = black,
            urlcolor  = blue,
            citecolor = red,
            anchorcolor = red]{hyperref}
\usepackage{url}
\usepackage{comment}
\usepackage{mathrsfs}  
\usepackage{longtable} %tables with pagebreak
\usepackage{enumitem}
\usepackage{relsize}
\usepackage{bm}
\usepackage{xcolor}
\usepackage{graphicx}
\usepackage{mathabx}
\usepackage{mathdots}
\usepackage{mathtools}
\usepackage{color}
\usepackage{tikz-cd,tikz}
\usepackage{braket}
\usepackage{stmaryrd}
\usepackage{amssymb}
\usepackage{tipa}
\usepackage[margin=1.2 in]{geometry}
\usepackage{appendix}
\usepackage{physics}

\makeatletter
\newcommand{\leqnomode}{\tagsleft@true}
\newcommand{\reqnomode}{\tagsleft@false}
\makeatother

\usetikzlibrary{arrows, decorations.markings}
\usetikzlibrary{positioning}
\usetikzlibrary{matrix}
\usetikzlibrary{shapes}
\usetikzlibrary{arrows,decorations.markings, tikzmark}
\usepackage{stmaryrd}

\pgfarrowsdeclare{bad to}{bad to}
{
  \pgfarrowsleftextend{-2\pgflinewidth}
  \pgfarrowsrightextend{\pgflinewidth}
}
{
  \pgfsetlinewidth{0.8\pgflinewidth}
  \pgfsetdash{}{0pt}
  \pgfsetroundcap
  \pgfsetroundjoin
  \pgfpathmoveto{\pgfpoint{-3\pgflinewidth}{4\pgflinewidth}}
  \pgfpathcurveto
  {\pgfpoint{-2.75\pgflinewidth}{2.5\pgflinewidth}}
  {\pgfpoint{0pt}{0.25\pgflinewidth}}
  {\pgfpoint{0.75\pgflinewidth}{0pt}}
  \pgfpathcurveto
  {\pgfpoint{0pt}{-0.25\pgflinewidth}}
  {\pgfpoint{-2.75\pgflinewidth}{-2.5\pgflinewidth}}
  {\pgfpoint{-3\pgflinewidth}{-4\pgflinewidth}}
  \pgfusepathqstroke
}
\theoremstyle{definition}
\newtheorem{definition}{Definition}[section]

\newtheorem{theorem}[definition]{Theorem}
\newtheorem{example}[definition]{Example}
\newtheorem{proposition}[definition]{Proposition}

\newtheorem{corollary}[definition]{Corollary}

\newtheorem{lemma}[definition]{Lemma}

\newtheorem{remark}[definition]{Remark}

\newtheorem{thm}{Theorem}

\newcommand{\full}{\textup{full}}

\newcommand{\BC}{{\mathbb{C}}}
\newcommand{\BD}{{\mathbb{D}}}

\newcommand{\BF}{{\mathbb{F}}}
\newcommand{\BG}{{\mathbb{G}}}

\newcommand{\BP}{{\mathbb{P}}}
\newcommand{\BQ}{{\mathbb{Q}}}
\newcommand{\BR}{{\mathbb{R}}}

\newcommand{\BV}{{\mathbb{V}}}
\newcommand{\BW}{{\mathbb{W}}}

\newcommand{\BZ}{{\mathbb{Z}}}

\newcommand\useleqno{\renewcommand\@eqnnum{\hb@xt@.01\p@{}%
                      \rlap{\normalfont\normalcolor
                        \hskip -\displaywidth(\theequation)}}}

\newcommand{\CC}{{\mathcal C}}

\newcommand{\CE}{{\mathcal E}}
\newcommand{\CF}{{\mathcal F}}

\newcommand{\CH}{{\mathcal H}}

\newcommand{\CK}{{\mathcal K}}
\newcommand{\CL}{{\mathcal L}}
\newcommand{\CM}{{\mathcal M}}

\newcommand{\CO}{{\mathcal O}}
\newcommand{\CP}{{\mathcal P}}

\newcommand{\CR}{{\mathcal R}}

\newcommand{\CV}{{\mathcal V}}
\newcommand{\CW}{{\mathcal W}}

\newcommand{\1}{\mathbf 1}

\DeclareMathOperator{\Aut}{Aut}

\DeclareMathOperator{\id}{id}

\DeclareMathOperator*{\res}{Res\,}

\newcommand{\Ext}{\mathcal{E}\text{xt}}

\newcommand{\Pic}{\mathop{\rm Pic}\nolimits}

\newcommand{\bL}{{\mathsf{L}}}

\newcommand{\bT}{{\mathsf{T}}}

\newcommand{\bR}{\mathsf{R}}

\newcommand{\pt}{{\mathsf{pt}}}

\newcommand{\congpf}{\xymatrix@1@=15pt{\ar[r]^-\sim&}}

\newcommand{\ch}{\mathrm{ch}}
\newcommand{\chh}{\mathrm{ch}^H}
\newcommand{\td}{\mathrm{td}}

\newcommand{\Hom}{\mathrm{Hom}}

\newcommand{\st}{\text{st}}
\renewcommand{\ss}{\text{ss}}

\newcommand{\sym}{\text{sym}}

\newcommand{\RHom}{\mathrm{RHom}}

\newcommand{\rk}{\mathrm{rk}}

\newcommand{\rig}{\mathrm{rig}}
\newcommand{\inv}{{\mathrm{wt}_0}}

\newcommand{\paar}{\mathrm{par}}
\newcommand{\Vpar}{\mathbf{V}^\mathrm{par}}
\newcommand{\Vpartr}{\mathbf{V}^\mathrm{par}_{\tr}}
\newcommand{\V}{\mathbf{V}}
\newcommand{\qpar}{\mathrm{qpar}}
\usepackage{tasks}

\usetikzlibrary{decorations.pathmorphing}
\AtBeginDocument{%
	\def\MR#1{}
}

\renewcommand{\Ext}{\textup{Ext}}

\usepackage{tikz}
\usepackage{tikz-cd}
\usetikzlibrary{decorations.pathmorphing}
\AtBeginDocument{%
	\def\MR#1{}
}

\begin{document}

\baselineskip=16pt
\parskip=5pt

\title[On the intersection theory of moduli spaces of parabolic bundles]{On the intersection theory of moduli spaces of parabolic bundles}
\author{Miguel Moreira}

\email{}

\date{\today}
\maketitle

\begin{abstract}
This paper concerns the intersection numbers of tautological classes on moduli spaces of parabolic bundles on a smooth projective curve. We show that such intersection numbers are completely determined by wall-crossing formulas, Hecke isomorphisms, and flag bundle structures and resulting Weyl symmetry. As applications of these ideas, we prove the Newstead--Earl--Kirwan vanishing -- and a natural strengthening in terms of Chern filtrations -- and the Virasoro constraints for parabolic bundles. Both of these results were already known for moduli of stable bundles without parabolic structure, but even in those cases our proofs are new and independent of the existing ones. We use a Joyce style vertex algebra formulation of wall-crossing, and define intersection numbers even in the presence of strictly semistable parabolic bundles; all of our results hold in that setting as well. 
\end{abstract}

\tableofcontents

\newpage
\section{Introduction}

Let $C$ be a smooth projective surface and let $M_{r,d}$ be the moduli space of stable vector bundles on $C$ of rank $r$ and degree $d$, where $r\geq 1, d\in \BZ$ are coprime integers. The moduli spaces $M_{r,d}$, and in particular their topology and their intersection theory, have been a topic of interest for more than 50 years. Let us mention, in a non-exhaustive way, some of the results that motivate the present paper.

Newstead \cite{newsteadrk2} (in rank 2) and Atiyah--Bott \cite{AB} (in higher rank) identified a set of generators for the cohomology of $M_{r,d}$, which we will refer to as descendents, or tautological classes. Two natural questions arise: what are the relations among descendents, and what are their intersection numbers. Both questions have been completely answered in the 90s. The study of relations culminated in the proof \cite{EKrel} that the set of Mumford relations is complete.

The question of determining intersection numbers -- which we also refer to as ``integrals'' -- of descendents was first solved in rank 2 by Thaddeus \cite{thaddeuscft} and Donaldson \cite{donaldsongluing}. Thaddeus proved his formula for integrals of descendents as a consequence of the Verlinde formula for the number of sections of determinant line bundles. Witten \cite{witintegrals} proposed a formula for such integrals in arbitrary rank, which was later proved by Jeffrey--Kirwan \cite{JK}. 

In this paper we propose a different approach to determining integrals on moduli spaces of bundles and, more generally, on moduli spaces of parabolic bundles \cite{seshadri, MS}. Given a marked point $p\in C$, a parabolic bundle on $(C,p)$ is essentially a vector bundle $V$ on $C$ together with a flag on the vector space $V_p$ (and also a choice of weights). Parabolic bundles can be thought of as vector bundles on an orbifold curve where we introduce a sigularity at $p$ via a root stack construction \cite{biswas}. It turns out that enriching bundles by adding a parabolic structure is very profitable: it introduces new features, namely wall-crossing, Hecke modifications, and Weyl symmetry, that allow us to completely determine all the integrals. This statement is the content of Theorem \ref{thm: reconstruction}, which we call the reconstruction theorem. It is analogous and greatly inspired by the recent proof of the Verlinde formula by Szenes--Trapeznikova \cite{ST}.

Although we do not pursue explicit formulas for the integrals, in principle these can be extracted from the reconstruction theorem. We give two applications of these ideas, namely a proof of the Newstead--Earl--Kirwan vanishing for parabolic bundles -- and a natural generalization in the context of Chern filtrations --, and the Virasoro constraints for parabolic bundles proposed in \cite{thesis}.

\subsection{Main results}
Given $r\in \BZ_{\geq 1}$, $d\in \BZ$, $f_\bullet=(f_1,\ldots, f_{l-1})$ and a choice of weights
\[0=c_{0}<c_1<c_2<\ldots<c_{l}<c_{l+1}=1\,,\]
we have a moduli space $M_{r,d,f_\bullet}(c)$ parametrizing $c$-semistable parabolic bundles with rank $r$, degree $d$ and dimension vector of the flag $f_\bullet$. We often write $\alpha$ for the data $(r,d,f_\bullet, c)$ and $M_\alpha=M_{r,d,f_\bullet}(c)$. For appropriate choices of $\alpha$ (which we call regular), all the semistable objects in this moduli space are stable, in which case $M_\alpha$ is a smooth projective variety of dimension
\[\dim M_\alpha=r^2(g-1)+1+\sum_{i=1}^{l} f_{i-1}(f_i-f_{i-1})\,.\]
There is a natural way to construct tautological classes (also called descendents) on the cohomology of the moduli spaces $M_\alpha$. Indeed, such moduli spaces carry a universal rank $r$ vector bundle $\BV$ on $M_\alpha\times C$, and a universal flag
\[0=\BF_0\subseteq \BF_1\subseteq \ldots \subseteq \BF_l=\BV_{|M_\alpha\times \{p\}}\]
of vector bundles on $M_\alpha$. The tautological classes that we consider are the Kunneth components of the Chern classes of $\BV$ and the Chern classes of $\BF_i$. It is well known that these generate $H^\ast(M_\alpha)$ as an algebra. We refer the reader to Section \ref{sec: prel} for more detail on parabolic bundles, stability, moduli spaces and tautological classes.

In this paper we are interested in intersection numbers of tautological classes, i.e.
\begin{equation}
\label{eq: integrals}
\int_{M_{r,d,f_\bullet}(c)}D
\end{equation}
where $D$ is a polynomial in the classes described before. Our main result is that such numbers can be fully reconstructed from some key properties of parabolic bundles:

\begin{thm}[Reconstruction theorem]\label{thm: reconstruction}
Integrals of descendents \eqref{eq: integrals} over moduli spaces of semistable parabolic bundles are completely determined by the following properties:
\begin{enumerate}\addtocounter{enumi}{-1}
\item The rank $r=1$ case, i.e. integrals over the Jacobian $M_{1,d}$;
\item The wall-crossing formula \eqref{eq: fullwcformula};
\item The flag bundle formula \eqref{eq: flagbundleclass} and Weyl anti-symmetry \eqref{eq: weylantisymmetry};
\item The Hecke compatibility \eqref{eq: heckeclass}.
\end{enumerate} 
\end{thm}

The theorem will be proved in Section \ref{subsec: proofreconstruction}. We will briefly explain what (1), (2) and (3) mean, and the key ideas in the proof, in Sections \ref{subsec: wallcrossingintro} and \ref{subsec: ingredientsintro}. This theorem is inspired by the proof of the Verlinde formula given by Szenes--Trapeznikova \cite{ST} (see also \cite{trapeznikova}), which is itself inspired by work of Teleman--Woodward \cite{TW}.

\begin{remark}
We hope that the reconstruction theorem will lead to a new proof of the Witten--Jeffrey--Kirwan formulas for integrals on $M_{r,d}$ as iterated residues, and to a generalization of those to parabolic bundles. Once there is a proposal for a formula for any moduli space of parabolic bundles, showing that it is correct just amounts to showing that it satisfies the same properties (0)-(3).
\end{remark}

As an application of these ideas, we prove two vanishing results for intersection numbers on moduli spaces of parabolic bundles. The first is addressed in Section~\ref{sec: newstead}, and is related to a conjecture of Newstead, which was proved for moduli spaces of semistable bundles by Earl--Kirwan \cite{EKvanishing} (based on the integral formulas from \cite{JK}). For parabolic bundles, a non-optimal result of the same type was proved by Gamse--Weitsman \cite{GW} with a more direct geometric/combinatorial argument. 

We state our result in terms of the Chern filtration (see Definition \ref{def: chernfiltration}). In the case of bundles, this formulation already appears in \cite[Theorem 0.2, Corollary 1.7, Proposition 1.8]{LMP}, where a proof is given following \cite{EKvanishing}.

\begin{thm}[=Theorem \ref{thm: newstead}]\label{thm: B}
If $D$ is a product of tautological classes with Chern degree $\leq \dim M_\alpha+r(g-1)$, then 
\[\int_{M_\alpha}D =0\,.\]
Moreover, if $D$ has Chern degree $\dim M_\alpha+r(g-1)+1$, then
\[(-1)^{(r-1)d}\int_{M_\alpha} D\,\]
does not depend on $d$ and on $c$.
\end{thm}

We refer the reader to Section 0.3 and Remark 1.11 in \cite{LMP} for a discussion concerning the relevance of the Chern filtration and on how this result is analogous to expected properties of moduli spaces of 1-dimensional sheaves on del Pezzo surfaces. 

A consequence of this result (cf. Corollary \ref{cor: newstead}) and Poincaré duality is the cohomological vanishing 
\begin{equation}
\label{eq: vanishingchernflag}
c_{k_1}(\BF_{j_1})\ldots c_{k_m}(\BF_{j_m})=0 \quad \in H^\ast(M_{r,d,f_\bullet}(c))
\end{equation}
for any $k_i\geq 0$ such that\footnote{For full parabolic bundles we get vanishing when the sum of the $k_i$ is at least $r(r-1)g-\binom{r}{2}+1$. The vanishing of \cite[Theorem 1.6]{GW} is the same, except that they only show it when the sum of the $k_i$ is bounded below by $r(r-1)g-r+2$, so our result is stronger for $r>2$. Indeed, our bound is optimal, see Remark \ref{rmk: optimal}.} 
\[k_1+\ldots+k_m\geq \dim M_\alpha-r(g-1)\] and $j_i\in \{1, 2, \ldots, l\}$.

The next application is a proof of the Virasoro constraints for moduli spaces of parabolic bundles. The Virasoro constraints are a set of linear relations among integrals of tautological on different moduli spaces.

\begin{thm}[=Theorem \ref{thm: virasoro}]\label{thm: C}
Every moduli space of parabolic bundles $M_\alpha$ satisfies the Virasoro constraints.
\end{thm}

Virasoro constraints for moduli spaces of sheaves and related objects have been investigated in \cite{moop, moreira,bree, blm, bojkoquivers, lmquivers}. The constraints for moduli spaces of bundles on curves, which are also a particular case of Theorem \ref{thm: virasoro}, were proved in \cite{blm}. The parabolic version was formulated in the author's thesis \cite{thesis}. A path of proof was sketched in loc. cit.; there, the proof suggested was by reducing to the case of bundles. However, the present proof does not depend on the case of bundles, and therefore is independent of the argument in \cite{blm} (but it does rely on the ideas there, most notably the wall-crossing compatibility).

\begin{remark}
Parabolic bundles with weights in $\frac{1}{N}\BZ$ correspond to orbibundles on the root stack $\mathcal C=\sqrt[N]{(C,p)}$, which introduces an orbifold singularity locally looking like $[\BC^1/\BZ_N]$ at the point $p$. In particular, the Virasoro constraints for parabolic bundles are the first instance of Virasoro constraints for sheaves on orbifolds. Note that if $I=\frac{1}{N}\BZ\cap [0,1]$ then the vector space $H_I(C)$ defined in Section \ref{subsec: descendents} can be naturally identified with the orbifold cohomology of $\mathcal C$. This hints that, in general, the orbifold descendents should be modeled in the orbifold cohomology. This is also what Lin proposes in \cite[Definition 5.23]{lin}, in the context of Pandharipande--Thomas invariants for orbifolds of dimension 3.
\end{remark}

\subsection{Wall-crossing, Joyce style vertex algebras, and stable$\neq$semistable}
\label{subsec: wallcrossingintro}
One of the key aspects of parabolic bundles, with particular importance in the present paper, is wall-crossing. By wall-crossing we mean the dependence of the moduli space $M_{r,d,f_\bullet}(c)$, and associated invariants \eqref{eq: integrals}, on $c$. The weight vector $c$ can be regarded as a stability condition; indeed, $M_{r,d,f_\bullet}(c)$ is constructed as a GIT quotient where the linearization is determined by the choice of $c$, cf. Section~\ref{subsubsec: stability}. The space of possible weight vectors $c$ (later denoted by $S_{r,d,f_\bullet}$) admits a wall-chamber decomposition; i.e., there are finitely many codimension 1 walls that partition the space of stability conditions into chambers, and the moduli space $M_{r,d,f_\bullet}(c)$ is the same when $c$ varies within a chamber, but might change when a wall is crossed. We say that $c$ (or, equivalently, $\alpha=(r,d,f_\bullet, c)$) is regular if it does not lie in any wall, in which case there are no strictly $c$-semistable sheaves, and the moduli space $M_{r,d, f_\bullet}(c)$ is smooth and proper.  

In Section \ref{sec: fullwc}, we prove Theorem \ref{thm: fullwc} (see also Theorems \ref{thm: simplewcI} and \ref{thm: simplewcII}), which provides a wall-crossing formula for integrals of tautological classes in moduli spaces of full parabolic bundles, i.e. when $f_\bullet=(1, 2, \ldots, r-1)$. Full parabolic bundles have the advantage that generic weight vectors $c$ are regular.  In concrete terms, letting $\alpha=(r,d,f_\bullet, c)$ and $\alpha'=(r,d,f_\bullet, c')$, the wall-crossing formula expresses the difference
\[\int_{M_\alpha}D-\int_{M_\alpha'}D\]
between integrals of tautological classes $D\in \BD^\paar$ (see Section \ref{subsec: descendents}) over moduli spaces defined with respect to two different regular stability conditions $c,c'$ as a sum of integrals over product of moduli spaces of parabolic bundles of smaller rank
\[\int_{M_{\alpha_1}\times \cdots \times M_{\alpha_l}} \Box\,;\]
the sum runs over some of the partitions\footnote{A partition of $\alpha$ is the same as the data of partitions
\[r=r_1+\ldots+r_l,\, d=d_1+\ldots+d_l,\,\{1, \ldots, r\}=J_1\sqcup \ldots \sqcup J_l\,,\]
with $|J_i|=r_i$.}
\[\alpha=\alpha_1+\ldots+\alpha_l\, ,\]
and the tautological classes $\Box$ inside the integrals are determined from $D$ by an explicit formula. 

Following Joyce, we write the precise wall-crossing formulas (in particular, describing $\Box$) using the vertex algebras $\Vpar$ and $\Vpar_{\tr}$, and the corresponding Lie algebras $\widecheck{\V}^\paar$ and $\widecheck{\V}^\paar_{\tr}$. The vertex algebra $\Vpar$ is exactly the vertex algebra constructed from the $\BC$-linear exact category of parabolic bundles following \cite{Jo17, Jo21}. On the other hand, $\Vpar_{\tr}$ is defined more concretely in terms of the descendent algebra and it is isomorphic to a lattice vertex algebra. Note that $\Vpar$ is naturally a vertex subalgebra of $\Vpar_{\tr}$ (cf. Proposition \ref{prop: embeddingVA}) and $\widecheck{\V}^\paar$ is a Lie subalgebra of $\widecheck{\V}^\paar_{\tr}$. For regular $\alpha$, the moduli space $M_\alpha$ defines a class
\[[M_\alpha]\in \widecheck{\V}^\paar\subseteq \widecheck{\V}^\paar_{\tr}\]
containing information about integrals
\[\int_{M_\alpha}D\,\]
of tautological classes $D\in \BD^\paar$ (see Section \ref{subsec: descendents}). The wall-crossing formula expresses the difference between two such classes in terms of Lie brackets of lower rank classes. 

An important feature of Joyce's wall-crossing theory, which is also present here, is that one also defines classes $[M_{\alpha}]$ for non-regular $\alpha$, which amounts to giving a meaning to the integrals \eqref{eq: integrals} when $c$ is non-regular. Such classes may appear in the wall-crossing formula, even between regular stability conditions. Our main results  -- Theorems \ref{thm: reconstruction}, \ref{thm: B} and \ref{thm: C} -- apply to non-regular weights when interpreted in this way. Our definition of $[M_{r,d,f_\bullet}(c)]$ utilizes full parabolic bundles, where we can always perturb a non-regular $c$ to a regular $c$. We prove in Theorem \ref{thm: agreement} that our definition, in the particular case of bundles $M_{r,d}$ (when $\gcd(r,d)\neq 1$ every $c$ is non-regular), agrees with Joyce's.

\begin{remark}\label{rmk: wconly1stability}
In his work, Joyce considers families of stability conditions on a fixed exact category. The nature of parabolic bundles, however, is different. The weight vector is part of the data of a parabolic bundle, and there is only one stability condition. Even if $c, c'$ are in the same chamber, the moduli spaces $M_{r,d,f_\bullet}(c)$, $M_{r,d,f_\bullet}(c')$ technically parametrize different parabolic bundles\footnote{They parametrize the same quasi-parabolic bundles. But, crucially, quasi-parabolic bundles do not form an exact category -- there is no natural notion of morphisms or subobjects for quasi-parabolic bundles, for instance -- and in particular the homology of the stack of quasi-parabolic bundles does not have a natural vertex algebra structure.} because they have different weights; indeed, the objects that they parametrize even have different classes in the $K$-theory of parabolic bundles. Hence, Joyce's setup does not immediately apply to parabolic bundles, although most of the ideas adapt easily with minor modifications, for instance requiring the wall-crossing formulas to hold only after forgetting weights (cf. Definition \ref{def: forgetweights}). 
\end{remark}

\subsection{Flag bundles, Hecke isomorphisms and (affine) Weyl symmetry}
\label{subsec: ingredientsintro}

We now briefly explained the remaining structures mentioned in the formulation of Theorem \ref{thm: reconstruction}.

\subsubsection{Flag bundles}
Often, there are flag bundle maps between different moduli spaces of parabolic bundles. The simplest example is when $\gcd(r,d)=1$; if $c$ is a \textit{small} stability condition, meaning that $c_{i+1}-c_i\ll 1$, then a parabolic bundle $(V, F, c)$ with rank $r$ and degree $d$ is stable if and only if $V$ is stable. In particular, there is a morphism $M_{r,d, f_\bullet}(c)\to M_{r,d}$ that forgets the data of the flag. See Proposition \ref{prop: flagbundle} for a more general statement. The intersection theory of the total space of a flag bundle can be related to the intersection theory of the base. This observation reduces Theorems \ref{thm: reconstruction}, \ref{thm: B} and \ref{thm: C} to the case of full parabolic bundles.

\subsubsection{Weyl symmetry}
The Weyl (anti-)symmetry appears when we consider a full flag bundle $M_{r,d}^\full(c)\to M_{r,d}$. As we explain in Section \ref{sec: hecke}, the cohomology $M_{r,d}^\full(c)$ admits an action of the symmetric group $\Sigma_r$ which permutes $t_1, \ldots, t_r\in H^2(M_{r,d}^\full(c))$; these classes $t_i$ are the first Chern classes of the line bundles $\BF_i/\BF_{i-1}$ obtained as successive quotients of the universal flag. Another consequence of the flag bundle structure is the anti-invariance of the integration functional:
\begin{equation}
\int_{M_{r,d}^\full(c)}\sigma\cdot D=(-1)^{\textup{sgn}(\sigma)}\int_{M_{r,d}^\full(c)}D\,
\end{equation}
for any $\sigma\in \Sigma_r$ and $D$ a polynomial in tautological classes. See Proposition~\ref{prop: weylflagbundle} and~\ref{cor: weylantisymmetry}. This anti-invariance property is what we mean in (2) of Theorem \ref{thm: reconstruction}. 

\subsubsection{Hecke isomorphisms}
The last ingredient in the reconstruction theorem are the Hecke isomorphisms, which we discuss in Section \ref{sec: hecke}. A feature of full parabolic bundles is that moduli spaces of parabolic bundles with different degrees are isomorphic, provided we change the weights accordingly. The basic idea behind this fact is the construction of Hecke modifications. If $(V, F_\bullet)$ is a parabolic bundle, then
\[\widetilde V=\ker\big(V\to V_p\to V_p/F_i\big)\]
is a new vector bundle of degree $\deg(V)-r+\dim(F_i)$ which inherits a parabolic structure from $V$. Such Hecke modifications preserve stability as we show in Lemma~\ref{lem: heckestability}. The upshot of Section \ref{sec: hecke} is Corollary \ref{cor: heckeiso}: given a rank $r$, two degrees $d, d'$ and a stability condition $c$ then we can always find another stability condition $c'$ such that
\[M_{r,d}^\full(c)\simeq M_{r, d'}^\full(c')\,.\]

\subsubsection{Affine Weyl symmetry}

One insight from \cite{TW, ST} is that by combining Hecke isomorphisms and the Weyl symmetry we obtain a much larger group of symmetries, namely the affine Weyl group. Consider the following diagram, where $c, c'$ are small stability conditions and $c''$ is so that we have an Hecke isomorphism as explained above.

\begin{center}
\begin{tikzcd}[column sep=2cm, row sep=1cm]
M_{r,1}^\full(c)\arrow[loop, dashed, distance=3em, "\Sigma_r"'] \arrow[r, leftrightsquigarrow, "\textup{wall-crossing}"] \arrow[d, "\textup{flag bundle}"'] & M_{r,1}^\full(c'')\arrow[loop, dashed, distance=3em, "\widetilde \Sigma_r"']\arrow[r, leftrightarrow, "\textup{Hecke}"]& M_{r,-1}^\full(c') \arrow[loop, dashed, distance=3em, "\Sigma_r"']\arrow[d , "\textup{flag bundle}"]\\
M_{r,1} & & M_{r,-1}
\end{tikzcd}
\end{center}

 The flag bundle structures imposes the $\Sigma_r$ anti-invariance of the integral functionals on $M_{r,1}^\full(c)$ and $M_{r,-1}^\full(c')$. Using the Hecke isomorphism, we obtain anti-invariance of the integral functional on $M_{r,1}^\full(c'')$ with respect to a different action of the symmetric group, which we now denote by $\widetilde \Sigma_r$; this action is obtained from conjugating the $\Sigma_r$ action by the Hecke isomorphism.
 
When the two groups $\Sigma_r, \widetilde \Sigma_r$ are put together by regarding both as groups of automorphisms of the formal descendent algebra $\BD^{\qpar}$, they generate a much larger group, the affine Weyl group of $\mathfrak{sl}_r$ (see Remark \ref{rmk: affineweyl}).

\begin{example}
Consider the two actions of the symmetric group $\Sigma_2$ and $\widetilde \Sigma_2$ on the lattice $\BZ$, generated by the involutions $x\mapsto -x$ and $x\mapsto 1-x$, respectively. Then, $\Sigma_2$ and $\widetilde \Sigma_2$ generate the entire group of automorphisms of $\BZ$, consisting of reflections and translations. This is isomorphic to the affine Weyl group of $\mathfrak{sl}_2$.
\end{example}

\subsubsection{Proofs of main results}\label{subsec: proofsintro}

The proof of the reconstruction theorem uses induction on $r$. Using wall-crossing, the Hecke operators and the flag bundle structures, it is enough to determine integrals on $M_{r,1}^\full(c)$. The proof of the reconstruction theorem then boils down to an algebraic fact that we prove in Proposition \ref{prop: affineweylvanishing}.  If it were the case that $c, c''$ were in the same chamber of $S_{r,d}$, the functional $\int_{M_{r,1}^\full(c)}D$ would be anti-invariant with respect to the actions of both $\Sigma_r, \widetilde \Sigma_r$ and hence Proposition \ref{prop: affineweylvanishing} would say that all the integrals vanish. Although this is never the case, wall-crossing gives the difference between integrals
\[\int_{M_{r,1}^\full(c)}D-\int_{M_{r,1}^\full(c'')}D=\textup{lower rank invariants}\]
on the two moduli spaces, so $\Sigma_r$ anti-invariance for $M_{r,1}^\full(c)$ and $\widetilde \Sigma_r$ anti-invariance for  $M_{r,1}^\full(c'')$ are still enough to recover the integrals over $M_{r,1}^\full(c)$. 

The proofs of Theorems \ref{thm: B} and \ref{thm: C} follow a strikingly similar strategy. In both cases, induction on $r$ and the wall-crossing formula are used to first show a wall-crossing independence statement, and then the affine Weyl anti-symmetry is used to show vanishing. In both cases we need a result showing compatibility with wall-crossing; for compatibility with the Chern filtration, that is the content of Lemma \ref{lem: cdegreewc}, and for the Virasoro constraints such statement is already known from \cite{blm}. Compatibility with Hecke operators and the flag bundle morphisms are also necessary. In that regard, the flag bundle compatibility for the Virasoro constraints is the most subtle part, and it is shown in Theorem \ref{thm: voamap}, using an identity from \cite[Appendix A]{klmp}.

\subsection{Rank 3 example}
\label{subsec: exampleintro}

Let us illustrate some of the structures explained above when $r=3$. We encourage the reader to come back to this example while reading Sections~\ref{sec: fullwc} to \ref{sec: weyl}. 

We first consider the case of full parabolic bundles with rank 3. In that case, the moduli space $M_{3, d}(c)$ depends on a choice of weights
\[0<c_1<c_2<c_3<1\,.\]
Actually, $M_{3, d}(c)$ only depends on the differences 
\[\lambda_1=c_2-c_1\textup{ and }\lambda_2=c_3-c_2\,.\]
In the picture below, we have the spaces of stability conditions
\[S_{3,d}=\{(\lambda_1, \lambda_2)\colon \lambda_1, \lambda_2>0,\, \lambda_1+\lambda_2<1\}\,\]
for $d=1, 0, -1$ and their respective wall-chamber structure:

\begin{center}
\includegraphics[scale=0.8]{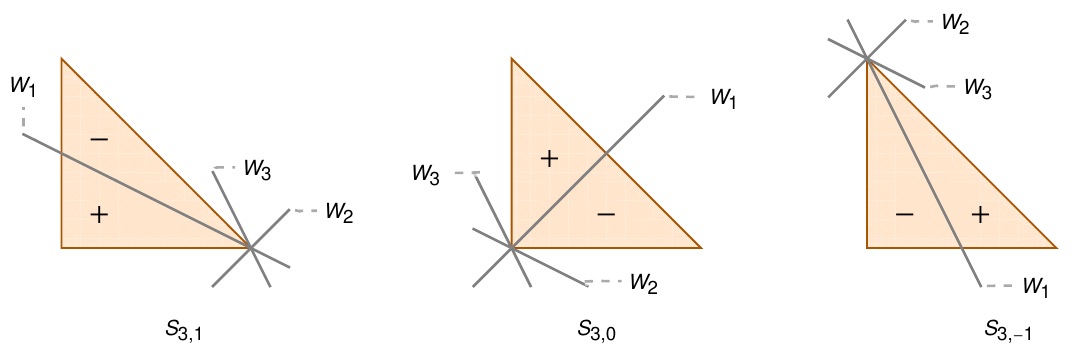}
\captionof{figure}{The spaces of stability conditions $S_{3, d}$ for $d=1, 0, -1$.}
\end{center}
\setlength{\parindent}{1em}

\bigskip

In each of the three cases $d=-1, 0, 1$, there is a unique wall $W_1$ which intersects $S_{3,d}$ (but we also draw the walls $W_2, W_3$ that intersect the closure $\overline S_{3,d}$) and divides the space of stability conditions into two chambers. In general, the walls are determined by the existence of a partition $\alpha=\alpha_1+\alpha_2$ of topological types such that $\alpha_1$ and $\alpha_2$ have the same slope (see the definition in Section \ref{subsubsec: stability}). For example, in the case of $d=1$, the topological types
\[\alpha_1=\big(r_1=2, d_1=1, f_1\big)\textup{ and }\alpha_2=\big(r_2=1,d_2=0, f_2\big)\]
where
\[f_1(t)=\begin{cases}
0\textup{ if }0\leq t\leq c_2\\
1\textup{ if }c_2<t\leq c_3\\
2\textup{ if }c_3<t\leq 1
\end{cases}\textup{ and }\quad f_2(t)=\begin{cases}
0\textup{ if }0\leq t\leq c_1\\
1\textup{ if }c_1<t\leq 1
\end{cases}\]
have the same slope if and only if $\lambda_1+2\lambda_2=1$, which is the equation for the wall~$W_1$.

Let $c^+$ be a stability conditions in the ``$+$'' chamber of $S_{3,1}$ (below the wall $W_1$) and $c^-$ in the ``$-$'' chamber (above the wall $W_1$). The moduli spaces $M_1=M_{\alpha_1}$ and $M_2=M_{\alpha_2}$ control the difference between $M_+=M_{3,1}^\full(c^+)$ and $M_-=M_{3,1}^\full(c^-)$, in the sense that $\int_{M_+}D-\int_{M_-}D$ can be expressed as an integral over $M_1\times M_2$. In the vertex algebraic formulation of wall-crossing, we have an identity (cf. Theorem \ref{thm: simplewcII})
\[[M_+]-[M_-]=\big[[M_1], [M_2]\big]\,\]
in $\widecheck \V^\paar$. If $c^0$ is on the wall $W_1$, then the moduli space $M_0=M^\full_{3,1}(c^0)$ includes strictly semistable objects, so $M_0$ is no longer smooth and does not admit a universal sheaf. We formally define the class $[M_0]$ -- i.e., give a definition of ``intersection numbers'' over $M_0$ -- by imposing the wall-crossing formula, which in this case reads as
\[[M_0]\coloneqq [M_-]+\frac 12 \big[[M_1], [M_2]\big]=[M_+]-\frac 12 \big[[M_1], [M_2]\big]\,.\]
Note that $[M_0]$ can be defined by using either the wall-crossing formula comparing $c^0$ and $c^-$ or $c^0$ and $c^+$, and this is consistent by the wall-crossing formula comparing $c^-$ and $c^+$.

The spaces of stability conditions $S_{3,d}$ for different $d$ can all be naturally identified via Hecke modifications, see Section \ref{subsec: heckespacesstability}. This correspondence identifies the 3 walls labeled as $W_1, W_2, W_3$ and the chambers labeled by ``$\pm$''. The moduli spaces $M_{3,d}(c)$ are isomorphic to $M_+$ for $c$ in the ``$+$'' chamber in any of the 3 spaces, and isomorphic to $M_-$ in the ``$-$'' chamber.  

We now address partial parabolic bundles. It is convenient to think of weight vectors for partial parabolic bundles as degenerate cases of weight vectors for full parabolic bundles where we allow some of the weights to be equal, or equivalently some of the differences $\lambda_i$ to be 0; see Section \ref{subsec: spacestaibilitiespartial}. Take for example $d=1$ and $f_\bullet=(2)$. Then $M_{3,1,(2)}(c)$ depends on the choice of two weights $0<c_1<c_2<1$, and we let $\overline c=(c_1, c_1, c_2)$ be a point on the side $\lambda_1=0$ of the compactified triangle $\overline S_{3,1}$. The wall-chamber structure on the space of stability conditions $S_{3,1,(2)}$ for partial parabolic bundles is obtained by restricting the wall-chamber structure on $\overline S_{3,1}$ to this edge. In particular, there are two different moduli spaces $M_{3,1,(2)}(c)$, depending on whether $\lambda_2$ is $>1/2$ or $<1/2$. On the other hand, since there is no wall intersecting the side $\lambda_2=0$, the moduli spaces $M_{3,1,(1)}(c)$ are all isomorphic, regardless of the choice of $c$. The vertex $\lambda_1=\lambda_2=0$ corresponds to the most degenerate case $M_{3,1}$, where we do not parametrize any flag.

When we perturb a regular $c\in \overline S_{r,d}$ to a ``less degenerate'' stability condition, we obtain flag bundle morphisms. For example, in the picture with $r=3, d=1$ we get the following fibrations:

\begin{center}
\begin{tikzcd}[column sep=3cm]
&M_+ \arrow[ld, "\BP^1\textup{ bundle}"']\arrow[rd,"\BP^1\textup{ bundle}"] \arrow[dd, "{\parbox{2cm}{$\mathsf{Fl}(\BC^3 ;1{,}2)$\\ \scriptsize bundle}}"]&\\
M_{3,1,(2)}(\lambda_2<1/2) \arrow[rd, "\BP^2 \textup{ bundle}"'] & & M_{3,1,(1)}(\lambda_1)  \arrow[ld, "\BP^2 \textup{ bundle}"]\\
& M_{3,1} &
\end{tikzcd}
\end{center}

There is also a $\BP^1$-bundle $M_-\to M_{3,1,(2)}(\lambda_2>1/2)$. These bundles can be used to calculate the intersection numbers on the base from the intersection numbers of the total space. This fact is nicely encoded by the vertex algebra homomorphisms $\Omega$ defined in Section \ref{subsec: mapsVA}. For example, 
\[[M_{3,1}]=\frac{1}{3!}\Omega([M_+])\,.\]
Such homomorphisms are also used to define classes $[M_{r,d, f_\bullet}(c)]$ for partial parabolic bundles when $c$ is not regular, see Definition \ref{def: joyceclassesboundary}. For example, 
\[[M_{3,0}]\coloneqq \frac{1}{3!}\Omega([M_{3,0}^\full(0)]\]
where the class $M_{3,0}^\full(0)$ is formally defined by a wall-crossing formula. Concretely, let $c$ be a stability condition in the ``$+$'' chamber of $S_{3,0}$, i.e. above $W_1$, and let $\alpha_i=(1,0,f_i)$ where $f_i(t)=0$ for $t\in [0, c_i]$ and $f_i(t)=1$ for $t\in (c_i, 1]$. Then
\begin{align*}[M_{3,0}^\full(0)]\coloneqq M_{\alpha_1+\alpha_2+\alpha_3}&-\frac{1}{2}\big[M_{\alpha_1}, M_{\alpha_2+\alpha_3}\big]-\frac{1}{2}\big[M_{\alpha_2}, M_{\alpha_1+\alpha_3}\big]+\frac{1}{2}\big[M_{\alpha_3}, M_{\alpha_1+\alpha_2}\big]\\
&+\frac{1}{3}\big[\big[M_{\alpha_1}, M_{\alpha_2}\big], M_{\alpha_3}\big]-\frac{1}{6}\big[\big[M_{\alpha_1}, M_{\alpha_3}\big], M_{\alpha_2}\big]\,.
\end{align*}
On the right hand side we abuse notation and write $M_\alpha$ for the class $[M_\alpha]$, to make the formula more readable. The complexity of the formula, when compared to the previous wall-crossing formula, is due to the fact that the origin $\lambda_1=\lambda_2=0$ is in the intersection of the 3 walls
\begin{align*}
W_1&\colon \mu(\alpha_2)=\mu(\alpha_1+\alpha_3)\\
W_2&\colon \mu(\alpha_3)=\mu(\alpha_1+\alpha_2)\\
W_3&\colon \mu(\alpha_1)=\mu(\alpha_2+\alpha_3)\,.
\end{align*}

We have the Hecke isomorphism $M_{\alpha_1+\alpha_2+\alpha_3}=M_{3,0}^\full(c)\simeq M_+$. The part of Theorem \ref{thm: newstead} concerning the independence of $d$ can be illustrated from these formulas. Integrals on $M_{3,1}$ are related to integrals on $M_+$; the same is true for integrals against the formally defined class $[M_{3,0}]$, modulo the wall-crossing corrections that vanish when we impose the condition on the Chern degree.

The full flag bundle $M_+\to M_{3,1}$ is responsible for the existence of an action of the Weyl group $\Sigma_3$ on $H^\ast(M_+)$, and for the anti-invariance of the functional $D\mapsto \int_{M_+} D$ with respect to this action. On the other hand, by making use of the Hecke isomorphism, we also have a full flag bundle $M_-\to M_{3,-1}$ which is responsible for the $\widetilde \Sigma_3$ action on $H^\ast(M_-)$ and analogous anti-invariance.

\subsection{Acknowledgements}
This paper draws great inspiration from \cite{ST} and the author benefited from conversations with O. Trapeznikova and A. Szenes, in particular during a visit to Université de Genève in 2022, which is more or less when the idea for this project started. The author also thanks D. Joyce, I. Karpov, W. Lim, D. Maulik, W. Pi and R. Pandharipande for conversations related to different aspects of this paper. The author was supported during part of the project by ERC-2017-AdG-786580-MACI. The project received funding from the European Research Council (ERC) under the European Union Horizon 2020 research and innovation programme (grant agreement 786580).

\subsection{Table of notation}

\begin{tasks}[label = {}, column-sep=-10em, item-indent=1.5em, after-item-skip=0.3em](2)

         \task $C$, $p$
        
         \task Fixed smooth projective curve $C$ and point $p\in C$
         
         \task $W_\bullet=(V, F_\bullet)$
        
         \task Quasi-parabolic bundle ($V$ vector bundle, $F_\bullet$ flag)
          
         \task $W=(V, F_\bullet, c)=(V, F)$ 
                 
         \task Parabolic bundle ($c$ weight vector)
         
         \task $\alpha_\bullet=(r,d,f_\bullet) $
         
         \task Topological type of a quasi-parabolic bundle

         \task $\alpha=(r,d,f_\bullet,c)=(r,d,f) $
         
         \task Topological type of a parabolic bundle    
         
         \task $I\subseteq [0,1]$
         
         \task Set of weights
         
         \task $C(I)$ (resp. $K(I)$)
         
         \task Monoid of topological types of parabolic bundles with weights in $I$ (resp. its Grothendieck group)
         
         \task $\CP=\CP_I$ 
         
         \task Category of parabolic bundles with weights in $I$
         
         \task $\chi_{\CP}$ (resp. $\chi_{\CP}^\sym$)
         
         \task Euler pairing on $K(I)$ (resp. its symmetrization)
         
         \task $\CM^\paar=\CM^I$ (res. $\CM^\qpar$)
         
         \task Moduli stack of parabolic bundles with weights in $I$ (resp. quasi-parabolic bundles)
                 
         \task $M_\alpha=M_{r,d,f_\bullet}(c)$ 
         
         \task Coarse moduli space of stable parabolic bundles of type $\alpha=(r,d,f_\bullet, c)$
         
         \task $M^\full_{r,d}(c)$ 
         
         \task Coarse moduli space of full stable parabolic bundles, i.e. $f_\bullet=(1, \ldots, r-1)$
         
          \task $M_{r,d}$ 
         
         \task Coarse moduli space of stable bundles, i.e. $f_\bullet=\emptyset$
         
         \task $\1$, $\pt$
        
         \task Unit class in $H^0(C)$ and Poincaré dual to a point in $H^2(C)$

         \task $\BD^\paar=\BD^I$
         
         \task Descendent algebra for parabolic bundles with \\ weights in $I$
         
         \task $\BD^\paar_\alpha$ (resp. $\BD^\qpar_{\alpha_\bullet}$)
         
         \task Descendent algebra for parabolic (resp. quasi-parabolic) bundles of topological type $\alpha$ (resp. $\alpha_\bullet$)
         
         \task $\ch_k(\gamma)$ and $\ch_k(e(t))$
         
         \task Generators of $\BD^\paar$  (or  $\BD^\paar_\alpha$)
         
         \task $\ch_k(\gamma)$ and $\ch_k(e_i)$
         
         \task Generators of  $\BD^\qpar_\alpha$
         
         \task $\CW=(\CV, \CF)$ 
         
         \task Universal parabolic bundle on stacks $\CM^\paar$
         
         \task $\BW=(\BV, \BF)$ 
         
         \task Universal parabolic bundle on moduli spaces $M_\alpha$
         
         \task $\Vpar=\V^I$
         
         \task Joyce's vertex algebra on $H_\ast(\CM^I)$
         
         \task $\Vpar_{\tr}=\V^I_{\tr}\supseteq \V^I$
         
         \task Descendent style vertex algebra
         
         \task $\widecheck \V^\paar=\widecheck \V^I$, $\widecheck \V^\paar_{\tr}=\widecheck \V^I_{\tr}$
         
         \task Associated Lie algebras
         
         \task $[M_\alpha]\in \widecheck \V^\paar\subseteq \widecheck \V^\paar_{\tr}$
         
         \task Class of a moduli space $M_\alpha$ in the Lie algebra
         
         \task $S_{r,d}$ (resp. $\overline S_{r,d}$)
         
         \task Space of stability conditions for full parabolic bundles (resp. its compactification)
         
         \task $\Omega_{I, I'}$
         
         \task Vertex algebra morphism $\V^I\to \V^{I'}$ (or induced Lie algebra homomorphism $\widecheck \V^I\to \widecheck \V^{I'}$)
         
          \task $h_\tau$
         
         \task Hecke modification, automorphism of the category $\CP_I$, of the stack $\CM^I$, or of the abelian group $K(I)$
         
         \task $(h_\tau)_\ast$ (resp. $(h_\tau)_\dagger$)
         
         \task Induced Hecke operator on $\V^\paar$ or $\widecheck \V^\paar$ (resp. $\V^\paar_{\tr}$ or $\widecheck \V^\paar_{\tr}$)
         
         \task $\deg$, $\deg^C$
       
       	 \task Cohomological and Chern degrees on $\BD^\paar$ (or variations), respectively

         \task $C_\bullet \BD^\paar, C_\bullet H^\ast(M_\alpha)$
       
       	 \task Chern filtration on the descendent algebra or cohomology of moduli spaces
         
         \task $\bL_n=\bR_n+\bT_n$
         
         \task Virasoro operators on $\BD^\paar$ 
         
         \task $L_n$ 
         
         \task Dual Virasoro operators on $\V^\paar_{\tr}$
         
         \task $\widecheck P_0\subseteq \widecheck \V^\paar_{\tr}$ 
         
         \task Lie subalgebra of primary states
\end{tasks}

\section{Preliminaries on parabolic bundles}\label{sec: prel}

Let $C$ be a smooth projective curve over $\BC$ and $p\in C$ a point that will be fixed throughout. We start with the definitions of (quasi-)parabolic bundles on $(C, p)$.

\begin{definition}\label{def: quasipar}
A quasi-parabolic bundle $W_\bullet=(V, F_\bullet)$ is a vector bundle $V$ on $C$ together with a flag
\[0=F_0\subsetneq F_1\subsetneq \ldots \subsetneq F_{l-1}\subseteq F_l=V_p\,.\]
We call $l$ the length of $W_\bullet$. The topological type of a quasi-parabolic bundle is the triple $\alpha_\bullet=(r, d, f_\bullet)$ where $r=\rk(V)$, $d=\deg(V)$ and $f_\bullet=(f_1,\ldots, f_{l-1})$ is the vector encoding the dimensions of the vector spaces in the flag, $f_j=\dim F_j$; in particular, the length is part of the topological type.  For convenience, we will also denote $f_0=\dim F_0=0$ and $f_l=\dim F_l=r$. A quasi-parabolic bundle with $f_\bullet=(1,2, \ldots, r-1)$ is called a full quasi-parabolic bundle. We say that $\alpha_\bullet$ is a full topological type if $f_\bullet=(1,2, \ldots, r-1)$.
\end{definition}

In what follows, $I$ is a subset of the interval $[0,1]$ such that $0,1\in I$. This interval will be fixed throughout most of paper, except in Section \ref{sec: partial}, so we will often suppress $I$ from the notation.

\begin{definition}\label{def: parabolic}
A parabolic bundle with weights in $I$ is a quasi-parabolic bundle $W_\bullet$ together with a choice of weights $c_0, \ldots, c_{l+1}\in I$ such that
\[0=c_{0}<c_1<c_2<\ldots<c_{l}<c_{l+1}=1\,.\]
The weights of a parabolic bundle define a partition of $I$. It is useful to encode the information of a parabolic bundles as a pair $W=(V, F)$ of a vector bundle $V$ together with an increasing lower semicontinuous function
\[F\colon I\to \{\textup{subspaces of }V_p\}\] defined by
\[F(t)=F_j \quad \textup{for }c_{j}<t\leq c_{j+1}\]
and $F(0)=0$. Such function is increasing in the sense that $t<t'$ implies that $F(t)\subseteq F(t')$.
\end{definition}

It is convenient to also define $F_+(t)$ by $F(t)=F_j$ for $c_j\leq t<c_{j+1}$ and $F_+(1)=F_l$, and let
\[\partial F(t)\coloneqq F_+(t)/F(t)=\begin{cases}F_{j}/F_{j-1} & \textup{for }t=c_j\,, j=1, \ldots, l\\
0 & \textup{ otherwise}
\end{cases}\,.\]

The topological type of a parabolic bundle $W=(V, F)$ is the triple $\alpha=(r,d,f)$ where $f\colon I\to \BZ_{\geq 0}$ is the increasing function $f(t)=\dim F(t)$. Given such function $f$, we denote by $\textup{im} f=f_\bullet=(f_1, \ldots, f_{l-1})$ the vector of images of $f$, excluding $0$ and $r$. Since the information contained in $f$ is exactly the same as the information of $f_\bullet$ and choice of weights, we sometimes write $\alpha=(r,d,f_\bullet, c)$ instead. We let $C(I)$ be the monoid of topological types of parabolic bundles with weights in $I$, i.e.
\begin{align*}C(I)=\{(r, d, f)\vert\, r&\in \BZ_{>0},\, d\in \BZ,\, f\colon I\to \BZ_{\geq 0} \textup{ increasing with }f(0)=0, f(1)=r\}\\
&\cup \{(0,0,0)\}\,.
\end{align*}
We also consider its Grothendieck group $K(I)$, i.e.
\begin{align*}K(I)=\{(r, d, f)\vert\, r&\in \BZ,\, d\in \BZ,\, f\colon I\to \BZ \textup{ such that } f(0)=0, f(1)=r\}\,.
\end{align*}

\subsection{The category of parabolic bundles}\label{subsec: paraboliccat}
Given $I$, we have a $\BC$-linear exact category $\CP=\CP_I$ of parabolic bundles with weights in $I$. A morphism $W=(V, F)\to W'=(V', F')$ in $\CP$ is a morphism of vector bundles $\varphi\colon V\to V'$ such that the induced $\varphi_p\colon V_p\to V'_p$ preserves the respective flags, in the sense that 
\[\varphi_p(F(t))\subseteq F'(t)\,,\quad t\in I\,.\]

The direct product of $W\oplus W'$ is defined by $(V\oplus V', F\oplus F')$ where 
$(F\oplus F')(t)=F(t)\oplus F'(t)$; exact sequences are defined similarly.

This category can be embedded into the abelian category of parabolic coherent sheaves, which has homological dimension 1. Indeed, if $W=(V, F)$ and $W'=(V', F')$ are two parabolic bundles, its Ext groups can be calculated \cite[Theorem 5.1]{mozgovoyquiverrepab} from the long exact sequence
\begin{align}
0\to \Hom_{\CP}(W, W')\to \Hom_C(V, V')\oplus \bigoplus_{t\in I}\Hom(F(t), F'(t)) \label{eq: homext}\\
\to \bigoplus_{t\in I}\Hom(F(t), F'_+(t))\to \Ext^1_{\CP}(W, W')\to \Ext^1_C(V, V')\to 0\,.\nonumber
\end{align}

In particular, the Euler pairing for parabolic bundles is given by
\begin{align*}\chi_{\CP}(W, W')&\coloneqq \dim \Hom_{\CP}(W, W')-\dim \Ext^1_\CP(W, W')\\
&=\chi_C(V,V')-\sum_{t\in I}\dim(F(t))\dim\big(\partial F'(t)\big)
\end{align*}
where $\chi_C$ is the usual Euler pairing between the underlying vector bundles on the curve. Note that $\chi_\CP(W, W')$ only depends on the topological types of $W, W'$ and defines a non-degenerate pairing in $K(I)$.
 
\subsection{Stack of parabolic bundles and stability}\label{subsec: moduliparabolic}
There is a stack $\CM^\paar=\CM^I$ parametrizing the objects in the exact category $\CP=\CP_I$; since $\CP$ has homological dimension 1, this stack is smooth. It admits a decomposition into connected components according to the topological type:
\[\CM^\paar=\bigsqcup_{\alpha\in C(I)}\CM^\paar_\alpha\,.\]

The stack $\CM^\paar$ comes equipped with a universal parabolic bundle $(\CV, \CF)$ consisting of a vector bundle $\CV$ on $\CM^\paar\times C$ and a collection of vector bundles
\[\CF\colon I\to \{\textup{vector bundles over }\CM^\paar\}\]
with the property that $\CF(0)=0$, $\CF(1)=\CV_{|\CM^\paar\times \{p\}}$. 

Similarly, there is a stack of quasi-parabolic bundles 
\[\CM^\qpar=\bigsqcup_{r, d, f_\bullet}\CM^\qpar_{r, d, f_\bullet}\,.\]

The connected component $\CM^\qpar_{r, d, f_\bullet}$ comes with a universal quasi-parabolic bundle $(\CV, \CF_\bullet)$ where
\[0=\CF_0\subseteq \CF_1\subseteq\ldots\subseteq \CF_{l-1}\subseteq \CF_l=\CV_{\CM^\paar\times \{p\}}\]
is a flag of vector bundles on $\CM^\qpar_{r, d, f_\bullet}$. We can understand $\CM^\qpar_{r, d, f_\bullet}$ explicitly as follows: if $\CM_{r,d}$ is the stack of vector bundles on $C$ and $\CV$ is the universal bundle on $\CM_{r,d}\times C$, then $\CM^\qpar_{r, d, f_\bullet}$ is the flag bundle
\[\CM^\qpar_{r, d, f_\bullet}=\textup{Fl}_{\CM_{r,d}}(\CV_{|\CM\times \{p\}}; f_1, f_2, \ldots, f_{l-1})\,.\]

Note that, for any $f$ with $\textup{im}(f)=f_\bullet$, the connected components of the stacks of parabolic and quasi-parabolic bundles are identified:
\[\CM^{\paar}_{r,d, f}=\CM^\qpar_{r, d, f_\bullet}\,.\]
In particular, there is a canonical map
\[\CM^{\paar}\to \CM^{\qpar}\]
that forgets the weights. 

\subsubsection{Stability}
\label{subsubsec: stability}
There is a slope function 
\[\mu\colon C(I)\to \BR\]
defined by
\[\mu(r,d,f)=\frac{d-\int_{0}^1 f(t)\,\textup{d}t}{r}\,\]
where in the integral we understand $f(t)$ to be defined for every $t\in [0,1]$ by extending $f$. In other words,
\[\int_{0}^1 f(t)=\sum_{j=0}^\ell (c_{j+1}-c_{j})f_j\,.\]
\begin{remark}
If $f_\bullet=(1, 2, \ldots, r-1)$ then 
\[\mu(r, d, f)=\frac{d+\sum_{j=1}^{r} c_j}{r}-1\,.\]
\end{remark}
A parabolic bundle is said to be semistable if $\mu(W')\leq \mu(W)$ for any $0\neq W'\subsetneq W$, and stable if the inequality is strict for any such $W'$.\footnote{If $V'\subseteq V$ is a subbundle, there is an induced parabolic subbundle $W'=(V', F')$ where $F'(t)=W'_p\cap F(t)$, called the saturated parabolic subbundle. Since $W'$ has larger slope than any other parabolic subbundle with the same underlying vector bundle $V'$, it is enough to test stability against saturated parabolic subbundles.}  We denote $\CM^{\st}_\alpha\subseteq \CM^{\ss}_\alpha\subseteq \CM_\alpha$ the stacks of stable and semistable objects. When  $\CM^{\st}_\alpha= \CM^{\ss}_\alpha$, we say that $\alpha$ is \textit{regular} (or, if $r, d, f_\bullet$ are fixed, we say that the weight vector $c$ is regular), and we write $M_\alpha$ for the coarse moduli space of stable parabolic bundles, which is the rigidification of this stack.\footnote{Note that the stabilizer group at any point $0\neq [W]\in \CM^\paar$ always contains a copy of $\BG_m$ given by scalar multiplication, and if $W$ is stable that is the entire stabilizer. The rigidification is the procedure of removing this $\BG_m$, see \cite[Appendix A]{AOV}.} The moduli space is a smooth projective variety of dimension
\[r^2(g-1)+1+\sum_{i=1}^{l-1}f_i\,.\]

We also denote
\[M_{r, d, f_\bullet}(c)=M_{r,d, f}\]
where $c$ are the weights of $f$. The moduli spaces $M_{r, d, f_\bullet}(c)$ are constructed as GIT quotients
\[X\sslash_{L_c} G\]
where $X, G$ depend only on $r,d,f_{\bullet}$ and $L_c$ depends on the weights. We refer to \cite{MS} for the details on the construction.  

\subsection{Descendents}\label{subsec: descendents}
One may extract enumerative invariants out of the moduli spaces $M_\alpha$ by integrating tautological classes constructed out of the Chern characters of $\CV$ and $\CF(t)$. We keep track of this information by introducing the formal descendent algebra. First, we define the $\BZ/2$-graded $\BQ$-vector space\footnote{Unless explicitly stated, cohomology $H^\ast(-)=H^\ast(-; \BQ)$ is always implicitly taken with $\BQ$ coefficients.}
\[H_I(C)=\Big(H^\ast(C)\oplus \bigoplus_{t\in I} \BQ\cdot e(t)\Big)/\langle e(0), \pt- e(1)\rangle\]
where $e(t)$ are formal symbols. The $\BZ/2$ grading is obtained from the $\BZ/2$ grading on $H^\ast(C)$ and by setting the generators $e(t)$ to be even. Note in particular that if $I=\{0,1\}$ then $H_{\{0,1\}}(C)=H^\ast(C)$. The $I$-parabolic descendent algebra $\BD^I$ is the free supercommutative algebra generated by the symbols
\[\ch_k(g)\,,\quad \textup{for }g\in H_I(C)\,,\]
which are understood to be linear in $g$. As usual, if $I$ is understood to be fixed, we write $\BD^\paar=\BD^I$. The algebra $\BD^\paar$ admits a natural cohomological grading given by
\[\deg(\ch_k(\gamma))=2k-2+\deg(\gamma)\textup{ and }\deg(\ch_k(e(t))=2k\,.\]
Let $p\colon \CM^\paar\times C\to \CM^\paar$ and $q\colon \CM^\paar\times C\to C$ be the projections of the product. We have a realization morphism
\[\xi\colon \BD^\paar\to H^\ast(\CM^{\paar})\,\]
which is a homomorphism of graded algebras and sends
\[\ch_k(\gamma)\mapsto p_\ast \big(\ch_k(\CV)q^\ast \gamma \big)\,,\quad\textup{for }\gamma\in H^\ast(C)\subseteq H_I(C)\]
and 
\[\ch_k(e(t))\mapsto \ch_k(\CF(t))\,,\quad\textup{for }t\in I\,.\]

Given $\alpha=(r,d,f) \in C(I)$, the realization morphism $\xi\colon \BD^\paar\to H^\ast(\CM^{\paar}_\alpha)$ factors through
\[\BD^\paar_\alpha=\frac{\BD^\paar}{\big\langle \ch_0(e(t))= f(t), \ch_1(\1)=d,  \ch_0(\gamma^{<2})=0\big\rangle}\]
where $\gamma^{<2}$ denotes an arbitrary class in $H^{<2}(C)$.

Similarly, given a quasi-parabolic topological type $\alpha_\bullet=(r, d, f_\bullet)$, we also have a quasi-parabolic descendent algebra
\[\BD_{\alpha_\bullet}^\qpar\]
generated by symbols $\ch_k(\gamma)$ for $\gamma\in H^\ast(C)$ and $\ch_k(e_j)$, for $j=0, 1, \ldots, l, l+1$ modulo relations $e_0=0$, $e_{l+1}=\pt$ and 
\[\ch_0(e_j)=f_j\,,\ch_1(\1)=d\,,\ch_0(\gamma^{<2})=0\,.\]
There is a realization map $\xi\colon \BD_{\alpha_\bullet}^\qpar\to H^\bullet(\CM^\qpar_{\alpha_\bullet})$ defined in a similar way to its parabolic version. There is also a natural map $\BD_{\alpha}^\paar\to \BD_{\alpha_\bullet}^\qpar$ which sends $\ch_k(e(t))$ to $\ch_k(e_j)$ where $j$ is such that $c_j<t\leq c_{j+1}$. This makes the square
\begin{center}
\begin{tikzcd}\BD_{\alpha}^\paar \arrow[r, "\xi"]\arrow[d] & H^\bullet(\CM^\paar_{\alpha})
\\
\BD_{\alpha_\bullet}^\qpar \arrow[r, "\xi"]&H^\bullet(\CM^\qpar_{\alpha_\bullet}) \arrow[u, equal]
\end{tikzcd}
\end{center}
commute.

\subsubsection{Weight 0 descendents}\label{subsubsec: wt0}
Let $\alpha\in C(I)$ be such that  $\CM^{\st}_\alpha=\CM^{\ss}_\alpha$. On the moduli space $M_\alpha$, there is a non-unique universal parabolic bundle $\BW=(\BV, \BF)$ where $\BV$ is a rank $r$ vector bundle on $M_\alpha\times C$ and $\BF(t)$ is a vector bundle of rank $f(t)$ on $M_\alpha$ for each $t\in I$. The universal property of $\BW$ is that its pullback along the rigidification map $\CM_\alpha^{\ss}\to M_\alpha$ coincides with the restriction of $\CW$ to $\CM_\alpha^{\ss}$. The non-uniqueness of $\BW$ stems from the fact that if $\BW$ is a universal parabolic bundle and $L$ is a line bundle on $M_\alpha$, then $\BW\otimes L\coloneqq (\BW\otimes q^\ast L, \BF\otimes L)$ is also a universal parabolic bundle. 

Given a choice of $\BW$, one obtains a realization homomorphism $\xi_\BW\colon \BD^\paar_\alpha\to H^\ast(M_\alpha)$. Often, one forces uniqueness of $\BW$ by requiring some normalization condition. Instead, we follow the approach in \cite{blm} and consider weight zero descendents. 

Let $\bR_{-1}$ be the derivation on $\BD^\paar$ or $\BD^\paar_{\alpha}$ which is defined on generators as
\[\bR_{-1}(\ch_k(g))=\ch_{k-1}(g)\,.\]

\begin{lemma}[{\cite[Lemma 2.8]{blm}}]\label{lem: R-1}
Given a parabolic universal bundle $\BW$, we have
\[\xi_{L\otimes \BW}=\sum_{j\geq 0}\frac{c_1(L)^j}{j!}\xi_{\BW}\circ \bR_{-1}^j\,.\]
\end{lemma}
It follows from the Lemma that if $\bR_{-1}(D)=0$ then $\xi_\BW(D)$ does not depend on the choice of $\BW$. Hence, we obtain canonical -- in the sense that they do not depend on the choice of $\BW$ -- realization homomorphisms
\[\xi\colon \BD^\paar_\inv\to H^\ast(M_\alpha)\quad \textup{or}\quad \xi\colon \BD^\paar_{\alpha,\inv}\to H^\ast(M_\alpha)\]
from the subalgebras of weight zero descendents
\[\BD^\paar_\inv=\ker(\bR_{-1})\subseteq \BD^\paar\,,\,\BD^\paar_{\alpha, \inv}=\ker(\bR_{-1})\subseteq \BD^\paar_\alpha\,.\]

The operator $\bR_{-1}$ can also be defined on $\BD^\qpar_{\alpha_\bullet}$, and we also define the algebra
\[\BD^\qpar_{\alpha_\bullet, \inv}=\ker(\bR_{-1})\subseteq \BD^\qpar_{\alpha_\bullet}\,.\]

\section{Vertex algebra and descendent integrals}\label{sec: vertexalgebras}

Joyce associates in a canonical way a vertex algebra to $\BC$-linear additive categories; in particular, this applies to the category of parabolic bundles $\CP_I$. We briefly recall here his construction, and then give a more concrete version of this vertex algebra in terms of the descendent algebra. We refer to \cite{Ka98} for an introduction to vertex algebras; see also \cite[Section 3]{blm} for a short presentation which contains everything that will be needed in this paper. 

\subsection{Joyce's vertex algebra}\label{subsec: joyceVA}

Joyce defines a vertex algebra with underlying vector space being the homology $\Vpar=H_\ast(\CM^\paar)$ of the stack parametrizing parabolic bundles. The vertex algebra structure is obtained from the following data on $\CM^\paar$:

\begin{enumerate}
\item The direct sum map $\Sigma\colon \CM^{\paar}\times \CM^\paar\to \CM^\paar$ sending a pair of parabolic bundles $(W_1, W_2)$ to $W_1\oplus W_2$. 
\item The $B\BG_m$ action $\rho\colon B\BG_m\times \CM^{\paar}\to \CM^{\paar}$ induced by scalar multiplication. 
\item A complex $\Ext_{12}$ (respectively $\Ext_{21}$) on $\CM^\paar\times \CM^\paar$ which, at a point $(W_1, W_2)$ is given by $\RHom_{\CP}(W_1, W_2)$ (respectively $\RHom_{\CP}(W_2, W_1)$). The restriction of $\Ext_{12}$ (respectively $\Ext_{21}$) to $\CM_{\alpha_1}\times \CM_{\alpha_2}$ has constant rank $\chi_\CP(\alpha_1, \alpha_2)$ (respectively  $\chi_\CP(\alpha_2, \alpha_1)$). 
\end{enumerate}

By \eqref{eq: homext}, the complex $\Ext_{12}$ is given in the $K$-theory of $\CM^\paar\times \CM^\paar$ by
\begin{equation}\Ext_{12}=Rp_\ast\mathcal{RH}\textup{om}(\CV_1, \CV_2)-\sum_{t\in I}\mathcal{H}\textup{om}(\CF_1(t), \partial\CF_2(t))\label{eq: extcomplex}
\end{equation}
where $\CV_1, \CV_2$ are the pullbacks of $\CV$ along the two possible projections $\CM^\paar\times \CM^\paar\times C\to \CM^\paar\times C$ and $p\colon \CM^\paar\times \CM^\paar\times C\to \CM^\paar\times \CM^\paar$.

From this data, one defines the vacuum vector $\ket{0}$ as the class of a point in $H_0(\CM_0)=H_0(\{0\})=\BQ$. The translation operator $T\colon \Vpar\to \Vpar$ is defined by
\[T(u)=\rho_\ast(t\otimes u)\]
where $t$ is the generator of $H_2(B\BG_m)$. Finally, the state-field correspondence is defined by the formula
\[Y(u,z)v=(-1)^{\chi_\CP(\alpha, \beta)}z^{\chi_\CP^\sym(\alpha, \beta)}\Sigma_\ast\big[(e^{zT}\boxtimes \textup{id})(c_{z^{-1}}(\Theta)\cap (u\boxtimes v)\big]\] 
for $u\in H_\ast(\CM^\paar_{\alpha})\subseteq \Vpar$ and $v\in H_\ast(\CM^\paar_{\beta})\subseteq \Vpar$, where $\Theta$ and $\chi_\CP^\sym$ are symmetrizations of the $\Ext$ complex and $\chi_\CP$, i.e.
\[\Theta=\Ext_{12}^\vee+\Ext_{21}\,\]
and
\[\chi_\CP^\sym(\alpha, \beta)=\chi_\CP(\alpha, \beta)+\chi_\CP(\beta, \alpha)\,.\]
In particular, $\Theta$ has rank $\chi_\CP^\sym(\alpha, \beta)$.

\subsection{Descendent vertex coalgebra}\label{subsec: coVA}

It is possible to define a similar but more concrete vertex algebra in terms of the descendent algebra, in which $\Vpar$ embeds. Indeed, the construction explained below produces a vertex algebra isomorphic to the one which Joyce would construct on the homology of the higher stack parametrizing objects in the derived category of parabolic bundles. 

The underlying vector space of the new vertex algebra is\footnote{The subscript $\tr$ stands for triangulated, and is motivated by the fact that this vertex algebra can be regarded as coming from the derived category of parabolic bundles, which is a triangulated category. More precisely, let $\CM^\paar_{\tr}$ be the (higher) moduli stack parametrizing parabolic complexes. This stack admits universal complexes $\CV$ and $\CF(t)$ on $\CM^\paar_{\tr}\times C$ and $\CM^\paar_{\tr}$, respectively; the induced realization homomorphism $\bigoplus_{\alpha\in K(I)} \BD_\alpha^{\paar}\to H^\ast(\CM^\paar_{\tr})$ is an isomorphism.}
\[\Vpartr=\bigoplus_{\alpha\in K(I)}(\BD^\paar_{\alpha})^\vee\,.\]
The vacuum vector is the obvious augmentation on $\BD^\paar_0$ which sends $1\in \BD^\paar_0$ to $1$ and any non-empty product of descendents to 0. The translation operator is defined to be the dual 
\[T=\bR_{-1}^\vee\colon (\BD^\paar_\alpha)^\vee\to (\BD^\paar_\alpha)^\vee\,.\]

To write down the state-field correspondence, we let
\[\Sigma^\ast \colon \BD^\paar_{\alpha+\beta}\to \BD^\paar_{\alpha}\otimes \BD^\paar_\beta\]
be the homomorphism defined on generators  by
\[\Sigma^\ast\ch_k(g)=\ch_k(g)\otimes 1+1\otimes \ch_k(g)\,.\]
Taking into consideration the formula \eqref{eq: extcomplex}, the Chern classes $c_{z^{-1}}(\Theta)\in H^\ast(\CM_\alpha\times \CM_\beta)\llbracket z^{-1}\rrbracket$ can be written in terms of descendents in a canonical way by applying Grothendieck--Riemann--Roch. We denote by $C_{z^{-1}}\in (\BD_{\alpha}\otimes {\BD_\beta})\llbracket z^{-1}\rrbracket$ this canonical lift of $c_{z^{-1}}(\Theta)\in H^\ast(\CM_\alpha\times \CM_\beta)\llbracket z^{-1}\rrbracket$ along the realization homomorphism
\[\xi\otimes \xi\colon \BD^\paar_{\alpha}\otimes {\BD^\paar_\beta}\to H^\ast(\CM_\alpha\times \CM_\beta)\,.\]
See for example \cite[(48)]{blm} for this lift written explicitly in a slightly different context.

The state-field correspondence is defined to be the dual of \[Y^\vee\colon \BD^\paar_{\alpha+\beta}\to (\BD^\paar_{\alpha}\otimes \BD^\paar_\beta)\llbracket z, z^{-1}\rrbracket\] given by the formula
\[Y^\vee=(-1)^{\chi_\CP(\alpha, \beta)}z^{\chi^{\sym}_\CP(\alpha, \beta)}C_{z^{-1}}( e^{z\bR_{-1}}\boxtimes\id)\Sigma^\ast\,.\]

\begin{proposition}\label{prop: latticeva}
The construction above makes $\Vpar_{\tr}$ into a vertex algebra isomorphic to a lattice vertex algebra in the sense of \cite[Theorem 3.5]{blm}.\footnote{Unlike in loc. cit., we use coefficients in $\BQ$ rather than $\BC$.} More precisely, following the notation in loc. cit., it is the vertex algebra associated to the data
\[\Lambda_{\mathrm{sst}}=K(I)\,,\Lambda_{\overline 0}=K(I)\otimes_\BZ \BQ\,, \Lambda_{\overline 1}=K^1(C)\otimes_\BZ \BQ\]
and $q$ is the Euler pairing, i.e.
\[q_{\overline 0}=\chi_{\CP}\textup{ and }q_{\overline 1}=(\chi_{C})_{|K^1(C)^{\otimes 2}}\,.\]
\end{proposition}
\begin{proof}
The proof for non-parabolic bundles (i.e. the $I=\{0,1\}$ case) can be found in \cite[Theorem 4.7]{blm}, and the parabolic case is an easy adaptation, so we will only say a few words. Recall that the descendent algebra is defined in terms of
\[H_I(C)=\frac{H^\ast(C)\oplus \BQ^I}{\langle e(0), e(1)-\pt\rangle }\]
and note that
\[\Lambda=\{(u, f)\in K^\ast(C)\oplus \BQ^I\colon f(0)=0\,, f(1)=\rk(u)\}\,.\]
There is a non-degenerate pairing $H^\ast(C)\otimes K^\ast(C)\to \BQ$ defined by the Chern character isomorphism together with Poincaré duality:
\[K^\ast(C)\simeq H^\ast(C)\simeq H^\ast(C)^\vee\,.\]
We also have the standard pairing $\BQ^I\otimes \BQ^I\to \BQ$. Together, these induce a non-degenerate pairing
\[H_I(C)\otimes \Lambda\to \BQ\,.\]
Following the proof of \cite[Lemma 4.8]{blm}, this non-degenerate pairing induces an identification of the vector spaces
\[(\BD_\alpha^\paar)^\vee\simeq \BD_\Lambda\]
and hence establishes an isomorphism of vector spaces between $\V^\paar_{\tr}$ and the vertex algebra associated to the data described. The fact that this isomorphism of vector spaces is also an isomorphism of vertex algebras is a calculation similar to the one in the proof of \cite[Theorem 4.7]{blm} (see also \cite[Theorem 5.13, Proposition 5.20]{lmquivers} for a fairly general statement, which does not quite apply in our setting due to the presence of odd cohomology classes), which crucially uses the formula \eqref{eq: extcomplex} for the $\Ext$ complex.
\end{proof}

The structure of the vertex algebra $\Vpar$ is best understood by embedding it into $\Vpar_{\tr}$. Recall that we have a realization homomorphism $\xi\colon \BD_\alpha^\paar\to H^\ast(\CM_\alpha)$. Taking duals, one obtains a canonical map
\[\xi^\vee \colon \Vpar\to \Vpar_{\tr}\,.\]

\begin{proposition}\label{prop: embeddingVA}
The canonical map $\Vpar\to \Vpar_{\tr}$ is an embedding of vertex algebras.
\end{proposition}
\begin{proof}
The map being a morphism of vertex algebras can be easily checked. For instance, the fact that it respects the translation operators can be shown as in \cite[Lemma 4.9]{blm} (see also Lemma \ref{lem: R-1}). For the state-field, observe that $\Sigma^\ast\colon \colon \BD^\paar_{\alpha+\beta}\to \BD_{\alpha}^\paar\otimes \BD_\beta^\paar$ and the pullback $\Sigma^\ast \colon H^\ast(\CM^\paar_{\alpha+\beta})\to H^\ast(\CM^\paar_{\alpha}\times \CM^\paar_{\beta})$ are intertwined by the realization map. 

The map $\Vpar\to \Vpar_{\tr}$ being injective is equivalent to the realization map $\xi\colon \BD^\paar_\alpha\to H^\ast(\CM_\alpha^\paar)$ being surjective, i.e. $H^\ast(\CM_\alpha^\paar)$ being generated by descendents. Recall that $\CM_\alpha^\paar$ is a flag bundle over the moduli stack $\CM_{r,d}$ of vector bundles on $C$; it follows that its cohomology is generated by elements pulled back from $H^\ast(\CM_{r,d})$ and the cohomology classes $\ch_k(\CF(t))$. Since the cohomology $H^\ast(\CM_{r,d})$ is well known to be generated by descendents $\xi(\ch_k(\gamma))$ (cf. \cite{AB}), we are done.\qedhere

\end{proof}

\subsection{Associated Lie algebra and intersection theory}\label{subsec: liealgebra}

Given a vertex algebra $\V$, the quotient $\widecheck \V=\V/T(\V)$ by the image of the translation operator admits a Lie algebra structure \cite{Borcherds}; the bracket is defined by
\[[\overline u, \overline v]=\overline{\Res_{z=0}Y(u,z)v}\,,\]
where $\overline u\in \widecheck \V$ denotes the image of $u\in \V$ in the quotient. 

If $\alpha$ is regular, the moduli spaces $M_\alpha$ define elements in the Lie algebras $\widecheck \V^\paar$ or $\widecheck \V_{\tr}^\paar$, as we now explain. Since $M_\alpha$ is the rigidification of $\CM_\alpha^{\ss}\subseteq \CM^\paar$, it admits an open embedding into the rigidification $M_\alpha\hookrightarrow (\CM^\paar)^\rig$. It is shown in \cite[Proposition 3.24]{Jo17} that $H_\ast((\CM^\paar)^\rig)$ is naturally isomorphic to $\widecheck \V^\paar$, so pushing forward the fundamental class of $M_\alpha$ along the embedding above defines a class
\[[M_\alpha]\in H_\ast((\CM^\paar)^\rig)\simeq \widecheck \V^\paar\,.\]

The class $[M_\alpha]\in \widecheck \V_{\tr}^\paar$ can be described in the alternative way that we now explain. Since by definition the translation operator on $\widecheck \V_{\tr}^\paar$ is the dual of $\bR_{-1}$, we have
\[ \widecheck\V_{\tr}^\paar=\bigoplus_{\alpha\in K(I)}(\BD^\paar_{\alpha, \inv})^\vee\,.\]

The moduli space $M_\alpha$ determines a functional 
\[[M_\alpha]\in(\BD^\paar_{\alpha, \inv})^\vee\subseteq \widecheck\V_{\tr}^\paar \]
given by 
\[D\mapsto \int_{M_\alpha}D\coloneqq \int_{M_\alpha}\xi(D)\]
where $\xi$ is the canonical realization homomorphism explained in Section \ref{subsubsec: wt0}. It is immediate that the two definitions of $[M_\alpha]$ are compatible with respect to the embedding $\widecheck\V^\paar\subseteq \widecheck\V_{\tr}^\paar$.

Note that the class $[M_\alpha]$, by definition, contains complete information about the intersection theory of $M_\alpha$. In Joyce's theory, wall-crossing formulas are expressed as formulas in the Lie algebras $\widecheck\V^\paar, \widecheck\V_{\tr}^\paar$. Our goal in the next sections is to prove such formulas for parabolic bundles.

\section{Wall-crossing formula for full parabolic bundles}\label{sec: fullwc}

 In this section, we will focus on the case of full parabolic bundles, i.e. $f_\bullet=(1, 2, \ldots, r-1)$, and denote the moduli spaces by $M_{r, d}^{\textup{full}}(c)$.

Given fixed $r$ and $d$, the moduli space $M^{\textup{full}}_{r,d}(c)$ depends on the choice of weights $c_1, \ldots, c_{r}$. Note that changing $c_1, \ldots, c_{r}$ to $c_1+s, c_1+s, \ldots, c_{r}+s$ does not affect stability; in other words, stability only depends on the difference of weights
\[\lambda_1=c_2-c_1\,,\lambda_2=c_3-c_2\,,\ldots\,, \, \lambda_{r-1}=c_{r}-c_{r-1}\,.\]

Hence, we will sometimes write $M_{r, d}^\full(\lambda)=M_{r, d}^\full(c)$. Our goal is to study the dependence on $\lambda$ of these moduli spaces and their intersection numbers.

\subsection{Space of stability conditions}\label{subsec: spacestabilities}

We define the \textit{space of (full) stability conditions}
\[S_{r,d}=\frac{\{0<c_1<\ldots<c_r<1\}}{(c_i)_{i=1}^r\sim (c_i+s)_{i=1}^r}=\{(\lambda_1, \ldots, \lambda_{r-1})\vert\, \lambda_i>0\,,\, \lambda_1+\ldots+\lambda_{r-1}<1\}\]
Given a partition $r=r_1+r_2$, $d=d_1+d_2$ and $\{1, \ldots, r\}=J_1\sqcup J_2$ with $|J_1|=r_1, |J_2|=r_2$ there is a wall
\[W_{J_1, J_2, d_1,d_2}\subseteq S_{r,d}\]
defined as follows:
\[\frac{d_1+\sum_{i\in J_1} c_i}{r_1}=\frac{d_2+\sum_{i\in J_2} c_i}{r_2}\]
Equivalently
\[r_2 d_1-r_1 d_2=\sum_{i\in J_1} r_2 c_i- \sum_{i\in J_2} r_1 c_i=\sum_{(i,j)\in J_1\times J_2} (c_i-c_j)\,.\]
This forces $\sum_{(i,j)\in J_1\times J_2} (c_i-c_j)$ to be an integer (congruent to $r_2d$ modulo $r$). Note that only finitely many walls intersect the simplex $S_{r,d}$. The following proposition is immediate from the description of the walls.

\begin{proposition}[No double walls]\label{prop: nodoublewalls}
If
\[[r]=J_1\sqcup J_2=J_1'\sqcup J_2'\textup{ and }d=d_1+d_2=d_1'+d_2'\]
are two different partitions, then the intersection of the walls $W_{J_1, J_2,d_1, d_2}\cap W_{J_1', J_2',d_1', d_2'}$ has codimension at least 2. 
\end{proposition}

For each $\alpha$ such that $c$ is on a wall, there is a distinguished decomposition $\alpha=\alpha_1+\alpha_2$ in $C(I)$ with $\mu(\alpha_1)=\mu(\alpha_2)=\mu(\alpha)$ associated to this wall, where $\alpha_i$, $i=1,2$, is the topological type of full parabolic bundles with ranks $r_i$, degree $d_i$ and weights $c_{J_i}=(c_j)_{j\in J_i}$. If $c\in S_{r,d}$ is not in any wall, there are no strictly semistable parabolic bundles of type $\alpha$, i.e. $\alpha$ is regular.

\begin{remark}\label{rmk: wcweird}
Suppose that $c^-, c^+\in S_{r,d}$ are both regular and they are in the same chamber of $S_{r,d}$, i.e. there is a continuous path in $S_{r,d}$ connecting $c^-$ and $c^+$ that does not cross any wall. Then, the corresponding moduli spaces are isomorphic. However, they define two different classes $[M_{\alpha_-}]$, $[M_{\alpha_+}]$ on the Lie algebra $\widecheck \V^\paar$ (or $\widecheck \V_{\tr}^\paar$) since they have different topological types. 
\end{remark}

To address the issue in the previous remark, we define the notion of equality after forgetting weights.

\begin{definition}\label{def: forgetweights}
Let $A, B\in \Vpar=H_\ast(\CM^\paar)$ (respectively $\widecheck \V^\paar=H_\ast((\CM^\paar)^{\rig})$). We say that $A=B$ \emph{after forgetting weights} if the images of $A, B$ along the map $H_\ast(\CM^\paar)\to H_\ast(\CM^\qpar)$ (respectively $H_\ast((\CM^\paar)^\rig)\to H_\ast((\CM^\qpar)^\rig)$) are the same.

Let $A, B\in \Vpar_{\tr}$ (respectively $\widecheck \V^\paar_{\tr}$). We say that $A=B$ \emph{after forgetting weights} if they induce the same functional on $\BD^\qpar_{\alpha_\bullet}$ (respectively $\BD^\qpar_{\alpha_\bullet, \inv}$).
\end{definition}

In particular, in the setting of Remark \ref{rmk: wcweird}, we have an equality $[M_{\alpha_-}]=[M_{\alpha_+}]$ after forgetting weights. 

\subsection{Simple wall-crossing}
\label{subsec: simplewc}
Let $\alpha_-$, $\alpha_+$ be the topological types of full parabolic bundles with rank $r$, degree $d$ and weights $c^-$ and $c^+$, respectively. We say that the wall-crossing between $c^-$ and $c^+$ is simple if there is a wall $W_{J_1, J_2, d_1,d_2}$ with the following properties:
\begin{enumerate}
 \item $c^-, c^+$ are regular;
\item We have
\[\frac{d_1+\sum_{j\in J_1}c_j^+}{r_1}<\frac{d_2+\sum_{j\in J_2}c_j^+}{r_2}\quad\textup{and}\quad\frac{d_1+\sum_{j\in J_1}c_j^-}{r_1}>\frac{d_2+\sum_{j\in J_2}c_j^-}{r_2}\,,\] 
so that $c^-$ and $c^+$ are on two different sides of the wall $W_{J_1, J_2, d_1,d_2}$.
\item The wall $W_{J_1, J_2, d_1,d_2}$ is the only one separating $c^-$ and $c^+$. I.e., there is a continuous path in $S_{r,d}$ that crosses the wall $W_{J_1, J_2, d_1,d_2}$ exactly once at $c^0$, and no other wall. 
\end{enumerate}

Let $\alpha_0$ be the topological type corresponding to $r,d, c^0$ and let $\alpha_0=\alpha_1+\alpha_2$ be the decomposition associated to the wall $W_{J_1, J_2, d_1,d_2}$. Note that the assumption that $c^0$ is only on one wall implies that $\alpha_1, \alpha_2$ are regular.

Let \[M_{\pm}=M_{\alpha_\pm}=M_{r,d}^\full(c_\pm)\,,\quad  M_{i}=M_{\alpha_i}=M_{r_i,d_i}^\full(c_{J_i})\,,\, i=1,2\,.\]
We will also denote by $\BW_+, \BW_-, \BW_1, \BW_2$ the universal parabolic bundles on these 4 moduli spaces.

\begin{theorem}[Simple wall-crossing I]\label{thm: simplewcI}
Consider the setup above. Then, for any $D\in \BD^{\qpar}_{\alpha_\bullet, \inv}$, we have
\begin{align*}\int_{M_+}&D-\int_{M_-}D\\
&=\res_{z=0}(-1)^{\chi(\alpha_1, \alpha_2)}z^{\chi_{\sym}(\alpha_1, \alpha_2)}\int_{M_1\times M_2}C_{z^{-1}}\cdot ( e^{z\bR_{-1}}\otimes \id)\Sigma^\ast D
\end{align*}
\end{theorem}

Note that, since $D\in  \BD^{\qpar}_{\alpha_\bullet, \inv}$, the integrals on the left hand side of the equality in Theorem \ref{thm: simplewcI} do not depend on the choices of universal parabolic bundles. While not obvious, it is also true that the right hand side does not depend on choices of universal parabolic bundles on $M_1, M_2$. Indeed, observe that the right hand side is equal to
\[\Res_{z=0}\int_{M_1\times M_2}Y^\vee(D)=\int_{\big[[M_1], [M_2]\big]}D\]
where $\big[[M_1], [M_2]\big]$ is the Lie bracket defined in Sections \ref{subsec: coVA}, \ref{subsec: liealgebra}; hence, independence of the choices of universal bundles is a consequence of the fact that this Lie bracket is well defined on the quotient $\widecheck \V^\paar_{\tr}$. In particular, from this discussion, Theorem \ref{thm: simplewcI} can be reformulated in the vertex algebraic language as follows:

\begin{theorem}[Simple wall-crossing II]\label{thm: simplewcII}
Consider the setup above. We have the identity
\[[M_+]-[M_-]=\big[[M_1], [M_2]\big]\]
after forgetting weights.
\end{theorem}
\begin{proof}
This is just a restatement of Theorem \ref{thm: simplewcII} and the definition of the bracket in $\Vpar_{\tr}$. \qedhere
\end{proof}

We will give the proof of Theorem \ref{thm: simplewcI} in Section \ref{sec: proofsimplewc}. The proof of simple wall-crossing utilizes the construction of a master space as in \cite{thaddeusGIT, mochizuki, Jo21, BH, ST}, which we recall below.

\subsubsection{Master space}\label{subsec: masterspace}

Moduli spaces of parabolic bundles can be constructed as GIT quotients. In this construction, 
\[M_{r,d}^{\textup{full}}(c)=X\sslash_{L_c} G\]
where $X$ is a flag bundle over (an open subset of) the Quot scheme parametrizing quotients $\CO_C^{\oplus N}\to V$, with $N$ sufficiently large; $G$ is $\textup{PSL}(N)$; $L_c$ is a rational ample line bundle with a $G$-linearization that depends on the weight vector $c$. Without loss of generality we may assume that the linear path
\[(1-t)c^-+tc^+\]
crosses the simple wall only at $t=1/2$. Denote $L_\pm=L_{c^\pm}$. The general construction of a master space \cite[Section 3]{thaddeusGIT} is as follows: the $\BP^1$ bundle $\BP_X(L_-\oplus L_+)$ over $X$ comes with a $G$ action induced from the $G$-actions of $L_\pm$ and with a $G$-linearization $\CO(1)$. Then the master space is
\[Z=\BP_X(L_-\oplus L_+)\sslash_{\CO(1)} G\,.\]

The master space $Z$ also comes equipped with a $\BC^\times$ action induced from the action of $\BC^\times$ on $L_-\oplus L_+$ that scales $L_+$, i.e. the action with weights $(0,1)$. Denote by $u\in K_{\BC^\times}$ the class of the tautological line bundle on $B \BC^\times$ and $z=c_1(u)\in H^2_{\BC^\times}$.\footnote{Note that the meanings of $u$ and $z$ are swapped compared to \cite{ST}, where $z$ is the $K$-theoretical variable and $u$ is the cohomological variable.}

 We record below the properties of the master space and refer to \cite[Section 5.2]{ST} for details:

\begin{enumerate}
\item $Z$ is a smooth projective variety with dimension $\dim M_{\pm}+1$. 
\item Both $M_-$ and $M_+$ embed into $Z$ via the inclusions $L_-\oplus 0\subseteq L_-\oplus L_+$ and  $0\oplus L_+\subseteq L_-\oplus L_+$.
\item The fixed point set $Z^{\BC^\times}$ decomposes into three connected components
\[Z^{\BC^\times}=M_-\sqcup M_+\sqcup Z_0\]
where $Z_0\simeq M_1\times M_2$. 
\item $Z$ comes with a (non-unique) $\BC^\times$-equivariant universal quasi-parabolic bundle $\BW=(\BV, \BF)$. The quasi-parabolic bundle $\BW$ is defined as in the Construction before Lemma 5.3 in \cite{ST}: by construction, $X$ comes with a quasi-parabolic bundle; the quasi-parabolic bundle on $Z$ is defined by pulling back along $\BP_X(L_-\oplus L_+)\to X$ and descending to the quotient.
\item The restrictions $\BW_{\pm}=\BW_{|M_\pm}$ are universal parabolic bundles on $M_\pm$. The restriction to $Z_0\simeq M_1\times M_2$ is identified with $\BW_1\otimes u+\BW_2$.
More precisely, this means that
\[\BV_{|Z_0}=\BV_1\otimes u+\BV_2\]
in $K_{\BC^\times}(M_1\times M_2\times C)$, where we omit the obvious pullbacks along the projections $M_1\times M_2\times C\to M_i\times C$, and 
\[\BF(t)_{|Z_0}=\BF_1(t)\otimes u+\BF_2(t)\]
in $K_{\BC^\times}(M_1\times M_2)$.
\item The normal bundles of $M_\pm$ inside $Z$ have weights $\mp 1$, i.e. 
\[N_{M_{-}/Z}=u\otimes \CL_-\textup{ and }N_{M_{+}/Z}=u^{-1}\otimes \CL_+\]
for non-equivariant line bundles $\CL_{\pm}\in \Pic(M_{\pm})$. 
\item The normal bundle of $Z_0$ is
\[N_{Z_0/Z}=\Ext_{12}[1]\otimes u^{-1}+ \Ext_{21}[1]\otimes u\,,\]
where $\Ext_{12}[1]$ is the vector bundle over $M_1\times M_2$ with fiber $\Ext^1_{\CP}(W_1, W_2)$ over $(W_1, W_2)\in M_1\times M_2$; similarly, $\Ext_{21}[1]$ is the vector bundle with fiber $\Ext^1_{\CP}(W_2, W_1)$. We observe that stability forces $\Hom_{\CP}(W_1, W_2)=0$ and $\Hom_\CP(W_2, W_1)=0$,\footnote{A homomorphism between stable objects is either an isomorphism or the zero morphism; if $W_1$ and $W_2$ were isomorphic then they would have the same set of weights, contradicting $\alpha_1+\alpha_2$ being the topological type of a full parabolic bundle.} so our use of $\Ext_{12}$ and $\Ext_{21}$ is consistent with Section \ref{subsec: joyceVA}. In particular,
\[\rk \,\Ext_{12}[1]=-\chi_\CP(\alpha_1, \alpha_2)\,,\quad \rk\, \Ext_{21}[1]=-\chi_\CP(\alpha_2, \alpha_1)\]
\end{enumerate}

\subsubsection{Proof of Theorem \ref{thm: simplewcI}} \label{sec: proofsimplewc}

Let $D\in \BD^\qpar_{\alpha_\bullet}$. By the equivariant localization formula on $Z$, we have

\begin{align}\label{eq: localization}
\int_{Z}\xi_{\BW}(D)&=\int_{M_-}\frac{\xi_{\BW}(D)_{|M_-}}{e_{\BC^\times}(N_{M_-/Z})}+\int_{M_+}\frac{\xi_{\BW}(D)_{|M_+}}{e_{\BC^\times}(N_{M_+/Z})}+\int_{Z_0}\frac{\xi_{\BW}(D)_{|Z_0}}{e_{\BC^\times}(N_{Z_0/Z})}\\
&=\int_{M_-}\frac{\xi_{\BW_-}(D)}{z+c_1(\CL_-)}+\int_{M_+}\frac{\xi_{\BW_+}(D)}{-z+c_1(\CL_+)}+\int_{Z_0}\frac{\xi_{\BW_1\otimes u+\BW_2}(D)}{e_{\BC^\times}(N_{Z_0/Z})}\nonumber
\end{align}
in $(H^\ast_{\BC^\times})_{\textup{loc}}=\BQ[z^{\pm 1}]$. The left hand side is a non-localized class, i.e. it is a polynomial in $z$, rather than a Laurent polynomial. In particular, $\Res_{z=0}$ annihilates the left hand side. Moreover, 
\[\res_{z=0}\int_{M_{\pm}}\frac{\xi_{\BW_\pm}(D)}{\mp z+c_1(\CL_\pm)}=\mp\int_{M_\pm}\xi_{\BW_\pm}(D)\,.\]

Hence, taking the residue in \eqref{eq: localization} we find that
\[\int_{M_+}\xi_{\BW_+}(D)-\int_{M_+}\xi_{\BW_+}(D)=\Res_{z=0} \int_{Z_0}\frac{\xi_{\BW_1\otimes u+\BW_2}(D)}{e_{\BC^\times}\big(\Ext_{12}[1]\otimes u^{-1}+ \Ext_{21}[1]\otimes u\big)}\,,\]
We now simplify the right hand side. First, note that
\[\xi_{\BW_1\otimes u+\BW_2}=(\xi_{\BW_1\otimes u}\otimes \xi_{\BW_2})\circ \Sigma^\ast=(\xi_{\BW_1}\otimes \xi_{\BW_2})\circ ( e^{z\bR_{-1}}\otimes \id)\circ \Sigma^\ast\]
as linear maps $\BD^{\paar}\to H^\ast_{\BC^\times}(M_1\times M_2)$. In the second equality we used Lemma \ref{lem: R-1}.

To relate the Euler class of the normal bundle to $C_{z^{-1}}$ we use the following Lemma, which is a straightforward calculation with Chern roots; since it is so fundamental in understanding Joyce's definition of the Lie bracket on $\widecheck \V^\paar$, we record it here:

\begin{lemma}
Let $A, B$ be vector bundles on $Z$. Then
\[e_{\BC^\times}(A\otimes u^{-1}+ B\otimes u)=(-1)^{\rk(A)}z^{\rk(A)+\rk(B)}c_{z^{-1}}(A^\vee+B)\,\]
holds in $H^\ast_{\BC^\times}(Z)$. 
\end{lemma}

Applying the Lemma to $A=\Ext_{12}[1]$, $B=\Ext_{21}[1]$, we obtain
\begin{align*}e_{\BC^\times}\big(\Ext_{12}[1]\otimes u^{-1}&+ \Ext_{21}[1]\otimes u\big)^{-1}\\
&=(-1)^{\chi(\alpha_1, \alpha_2)}z^{\chi_{\sym}(\alpha_1, \alpha_2)}c_{z^{-1}}(\Ext^\vee_{1,2}+\Ext_{2,1})\\
&=(-1)^{\chi(\alpha_1, \alpha_2)}z^{\chi_{\sym}(\alpha_1, \alpha_2)}(\xi_{\BW_1}\otimes \xi_{\BW_2})(C_{z^{-1}})\,.
\end{align*}
This concludes the proof of Theorem \ref{thm: simplewcI}.\qed

\begin{remark}
The same proof can be adapted to a $K$-theoretical wall-crossing formula in the style of \cite{liu}. The two fundamental properties of the residue that we have used are precisely the defining properties of what is called a residue map in \cite[Appendix A]{liu}.
\end{remark}

\subsection{General wall-crossing formula for full parabolic bundles}\,
\label{subsec: fullwcformula}

The general wall-crossing formula relates moduli spaces of full parabolic bundles with different weights that are not necessarily contained in adjacent chambers. To state it, we have to introduce some notation. Let $r\in \BZ_{\geq 1}, d\in \BZ$ and $c, c'\in S_{r,d}$. Given a subset 
\[J=\{j_1, \ldots, j_{r'}\}\subseteq [r]=\{1, \ldots, r\}\,,\] 
with $j_1< \ldots<j_{r'}$, we define the restriction of weights
\[c_J=\{c_{j_1}, \ldots, c_{j_{r'}}\}\,.\]

We will consider partitions $\vec{J}=(J_1, \ldots, J_m)\vdash [r]$ of $[r]$, i.e. a collection of disjoint subsets of $[r]$ such that $[r]=J_1\sqcup \ldots \sqcup J_m$, and partitions $\vec{d}=(d_1, \ldots, d_m)\vdash d$, i.e. integers such that $d_1+\ldots+d_m=d$; note that $\vec{J}\vdash [r]$ also induces a partition of the rank $\vec{r}=(r_1, \ldots, r_m)\vdash r$, where $r_i=|J_i|$.

\subsubsection{Wall-crossing coefficients}\label{subsubsec: wccoefficients}
Writing the general wall-crossing formula for two stabilities separated by more than a simple wall requires the wall-crossing coefficients \cite{JO06IV, Jo21}
\[U(\vec{J}, \vec{d}; c, c')\in \BQ\,\]
which depend on partitions $\vec{J}\vdash [r]$ and $\vec d\vdash d$ of the same length $m$ and weights $c, c'\in S_{r,d}$. Although stability for parabolic bundles does not exactly fit into the setup where such coefficients are defined (cf. Remark \ref{rmk: wconly1stability}), they can still be defined with minor modifications, which we now describe.

First, given a subset $J\subseteq [r]$, $r'=|J|$ and $d'\in \BZ$ we denote
\[\mu_{c}(J, d')=\frac{d'+\sum_{j\in J} c_j}{r'}\,.\]
One starts by defining constants $S(\vec J, \vec d; c, c')\in \{0,-1, 1\}$ analogously to \cite[Definition 3.10]{Jo21}; conditions $(a), (b)$ in loc. cit. should be replaced by
\begin{align*}\tag{a}
\mu_{c}(J_i, d_i)&\leq \mu_{c}(J_{i+1}, d_{i+1})\quad\textup{and}\\ \mu_{c'}(J_1\sqcup \ldots\sqcup J_i, d_1+\ldots+d_i)&>\mu_{c'}(J_{i+1}\sqcup \ldots\sqcup J_l, d_{i+1}+\ldots+d_l)
\end{align*}
\begin{align*}\tag{b}\mu_{c}(J_i, d_i)&> \mu_{c}(J_{i+1}, d_{i+1})\quad\textup{and}\\ \mu_{c'}(J_1\sqcup \ldots\sqcup J_i, d_1+\ldots+d_i)&\leq  \mu_{c'}(J_{i+1}\sqcup \ldots\sqcup J_l, d_{i+1}+\ldots+d_l)\,.
\end{align*}

If for every $i=1, \ldots, m-1$ we have either (a) or (b), $S(\vec J, \vec d; c, c')$ is defined to be $(-1)^{\#\textup{occurences of (a)}}$; otherwise, $S(\vec J, \vec d; c, c')=0$. From here, one defines coefficients\footnote{In our notation, $U^\ast$ is what Joyce calls $U$ and $U$ is reserved for what Joyce calls $\tilde U$.} $U^\ast(\vec{J}, \vec{d}; c, c')\in \BQ$ via \cite[(3.3)]{Jo21}; as in the definition of the $S$ constants, whenever there is a sum of topological types in Joyce's formulas one should replace it by a disjoint union of subsets $J_i$. The coefficients $U(\vec{J}, \vec{d}; c, c')\in \BQ$ are then defined via \cite[Theorem 3.12]{Jo21}.

\begin{theorem}\label{thm: fullwc}
Let $c, c'\in S_{r,d}$ be regular. Then we have
\begin{align}\label{eq: fullwcformula}
[M_{r,d}^{\textup{full}}(c')]=\sum_{\substack{\vec J\vdash [r]\\
\vec d\vdash d}}U(\vec J, \vec d;&\, c, c')\cdot \nonumber\\[-8mm] 
&\big[\big[\ldots\big[[M_{r_1,d_1}^{\textup{full}}(c_{J_1})], [M_{r_2,d_2}^{\textup{full}}(c_{J_2})]\big],\ldots,\big], [M_{r_m,d_m}^{\textup{full}}(c_{J_m})]\big]
\end{align}
after forgetting the weights.
\end{theorem}

The coefficients $U(\vec J, \vec d; c, c')$ are defined so that the right hand side of \eqref{eq: fullwcformula} expands as
\[\sum_{\substack{\vec J\vdash [r]\\
\vec d\vdash d}}U^\ast(\vec J, \vec d; c, c')[M_{r_1,d_1}^{\textup{full}}(c_{J_1})]\ast [M_{r_2,d_2}^{\textup{full}}(c_{J_2})]\ast \ldots \ast [M_{r_m,d_m}^{\textup{full}}(c_{J_m})]\big]\]
in the universal enveloping algebra of $\widecheck \V^\paar$ or $\widecheck \V^\paar_{\tr}$. 

\begin{proof}[Proof of Theorem \ref{thm: fullwc}]
Suppose $c=c^-$, $c'=c^+$ are separated by a simple wall so that we are in the setup of Section \ref{subsec: simplewc}. We have from the definition of the coefficients that
\begin{align*}
U^\ast([r], d; c^-, c^+)&=S([r], d; c^-, c^+)=1\\
U^\ast((J_1, J_2), (d_1, d_2); c^-, c^+)&=S((J_1, J_2), (d_1, d_2); c^-, c^+)=1\\
U^\ast((J_2, J_1), (d_2, d_1); c^-, c^+)&=S((J_2, J_1), (d_2, d_1); c^-, c^+)=-1\,
\end{align*}
and all the remaining coefficients are 0. So \eqref{eq: fullwcformula} is equivalent to
\[[M_+]=[M_-]+[M_1]\ast [M_2]-[M_2]\ast [M_1]=[M_-]+\big[[M_1], [M_2]\big]\,,\]
which is precisely  Theorem \ref{thm: simplewcII}.

In general, by Proposition \ref{prop: nodoublewalls} we can find a sequence $c^0, c^1, \ldots, c^k\in S_{r,d}$ with $c^0=c$ and $c^k=c'$ such that each consecutive weight vectors $c^i$ and $c^{i+1}$ are separated by a simple wall. By the formal property of the coefficients in \cite[Theorem 4.8]{JO06IV} (see also for instance the proof of \cite[Proposition 6.1]{GJT}), the wall-crossing formula for $c, c''$ and the wall-crossing formula for $c'', c'$ imply the wall-crossing formula for $c, c'$, so we conclude the general case from the case of simple wall-crossing.\qedhere
%\begin{align*}\mu_{c_-}(I_1, d_1)&>\mu_{c_-}(I_2, d_2)\\
%\mu_{c_+}(I_1, d_1)&<\mu_{c_+}(I_2, d_2)\,.
%\end{align*}
\end{proof}

%\mcomment{This means that extensions $0\to W_2\to W\to W_1\to 0$ are $c_-$ stable but $c_+$ unstable and vice-versa. Hence:
%\[[M_-]-\BP^{-\chi(\alpha_1, \alpha_2)}[M_1\times M_2]=[M_+]-\BP^{-\chi(\alpha_2, \alpha_1)}\]
%in particular
%\[e(M_+)-e(M_-)=\chi(\alpha_1, \alpha_2)-\chi(\alpha_2,\alpha_1)e(M_1)e(M_2)\]
%so
%\[M_+-M_-=[M_1, M_2]\]
%hopefully...}\end{example}

\subsection{Extension to non-regular weights}\label{subsec: joyceclassesfull}

Suppose that $c'\in S_{r,d}$ is not necessarily regular. We can define Joyce style classes by imposing the wall-crossing formula. 

\begin{definition}\label{def: joyceclasses}
Let $c'\in S_{r,d}$ and let $\alpha\in C(I)$ be the corresponding topological type. Choose an arbitrary $c\in S_{r,d}$ regular and define the class
\[[M_{r,d}^{\textup{full}}(c')]\in H_\ast((\CM^\paar_\alpha)^{\rig})\subseteq \widecheck \V^\paar\subseteq \widecheck \V^\paar_{\tr}\]
as the unique class that satisfies equation \eqref{eq: fullwcformula} after forgetting weights. 
\end{definition}

These classes are well defined and satisfy the general wall-crossing formula.

\begin{proposition}\label{prop: welldefinedclass}
Definition \ref{def: joyceclasses} does not depend on the choice of regular $c\in S_{r,d}$. Moreover, \eqref{eq: fullwcformula} holds for any $c,c'\in S_{r,d}$ not necessarily regular.
\end{proposition}
\begin{proof}
Both claims follow from Theorem \ref{thm: fullwc} and the fact that the wall-crossing formula is well-behaved under composition, i.e. the formal properties of the coefficients in \cite[Theorem 4.8]{JO06IV}. 
\end{proof}

\section{Partial parabolic bundles}\label{sec: partial}

One feature of parabolic bundles is that one can obtain maps between different moduli spaces by forgetting part of the data of the flag. In this section we explore this structure from the vertex algebra point of view and use it to construct Joyce style classes for partial parabolic bundles.

\subsection{Space of stability conditions for partial parabolic bundles}
\label{subsec: spacestaibilitiespartial}
When one considers partial parabolic bundles, it is natural to consider a compactification of the space of stability conditions $S_{r,d}$ from the previous section by allowing some of the weights to coincide:
\[\overline S_{r,d}=\frac{\{0\leq c_1 \leq \ldots \leq c_r\leq 1\}}{\sim}\simeq \{(\lambda_1, \ldots, \lambda_{r-1})\vert\, \lambda_i\geq 0\,,\, \lambda_1+\ldots+\lambda_{r-1}\leq 1\}\,.\]

Given $f_\bullet=(f_1, \ldots, f_{l-1})$ we let
\[S_{r,d,f_\bullet}=\frac{\{0< c_1 < \ldots <c_l< 1\}}{\sim}\,.\]
There is a natural inclusion of $S_{r,d,f_\bullet}\hookrightarrow \overline S_{r,d}$ defined by
\begin{equation}\label{eq: overlinec}S_{r,d,f_\bullet}\ni(c_1, \ldots, c_l)=c\mapsto \overline c=(c_1, \ldots, c_1, c_2, \ldots, c_2, \ldots, c_l, \ldots, c_l)\in \overline S_{r, d}
\end{equation}
where $c_j$ appears $f_{j}-f_{j-1}$ times. Conversely, given $\overline c\in \overline S_{r,d}$ we can recover the dimension function $f(t)$ (and hence $f_\bullet$ and $c$) as
\[f(t)=\#\{j\in \{1, \ldots, r\}\vert\, \overline c_j\leq t\}\,.\]

We observe that the definition of the walls $W_{J_1, J_2, d_1, d_2}$ extends to the compactification $\overline S_{r,d}$, and hence define a wall-chamber structure on $S_{r,d, f_\bullet}$. A decomposition $J_1\sqcup J_2=[r]$ induces a decomposition $f(t)=f_1(t)+f_2(t)$ where
\[f_i(t)=\#\{j\in J_i\vert\, \overline c_j\leq t\}\,,\,i=1,2\]
and 
\[\mu(r_1, d_1, f_1)=\mu(r_2, d_2, f_2) \,\Leftrightarrow \,\overline  c\in W_{J_1, J_2, d_1, d_2}\,.\]
In particular, the moduli spaces $M_{r,d,f_\bullet}(c)$ are unchanged when $c$ varies within a chamber on $S_{r,d, f_\bullet}$. A weight vector $c\in S_{r,d, f_\bullet}$ is regular if it does not lie in any wall.

The collection of subspaces $\{S_{r,d,f_\bullet}\}_{f_\bullet}$ defines a stratification of \[\overline S_{r,d}\cap\{\lambda_1+\ldots +\lambda_{r-1}<1\}\,;\] note that $S_{r,d,(1, \ldots, r-1)}=S_{r,d}$ is the open locus. Let $\overline S_{r,d,f_\bullet}$ be the closure of $S_{r,d,f_\bullet}$ inside $\overline S_{r,d}$. Note that $S_{r,d,f'_\bullet}\subseteq \overline S_{r,d,f_\bullet}$ if and only if $f'_\bullet\subseteq f_\bullet$.

\subsubsection{Flag bundles}\label{subsec: flagbundles} Suppose that $f_\bullet'\subseteq f_\bullet$, so that $S_{r,d,f'_\bullet}\subseteq \overline S_{r,d,f_\bullet}$. Concretely, let
\[f_\bullet=(f_1, \ldots, f_l)\textup{ and }f'_\bullet=(f_1', \ldots, f_k')\]
where $f'_i=f_{j_i}$ for some $1\leq j_1<\ldots<j_k\leq l$. If $(V, F_\bullet)$ is a quasi-parabolic bundle of type $(r,d,f_\bullet)$, then we let $(V, F'_\bullet)$ be the quasi-parabolic bundle of type $(r,d,f_\bullet')$ with associated flag $F'_i=F_{j_i}$. 

\begin{proposition}\label{prop: flagbundle}
Let $c'\in S_{r,d,f_\bullet'}\subseteq \overline S_{r,d,f_\bullet}$ be a regular weight vector. Let $c\in S_{r,d,f_\bullet}\subseteq \overline S_{r,d,f_\bullet}$ be sufficiently close to $c'$ so that $c$, $c'$ are in the same chamber. Then $(V, F_\bullet)$ is $c$-stable if and only if $(V, F_\bullet')$ is $c'$-stable. Hence, there is a map
\[\pi\colon M_{r,d,f_\bullet}(c)\to M_{r,d,f_\bullet'}(c')\]
with fiber over $(V, F'_\bullet)$ given by
\[\prod_{i=0}^k \mathsf{Flag}(F'_{i+1}/F'_{i}; f_{j_i+1}-f_{j_i}, f_{j_i+2}-f_{j_i}, \ldots, f_{j_{i+1}-1}-f_{j_i})\]
\end{proposition}
\begin{proof}
The proof is well-known, see \cite[Lemma 8.3]{ST} for a particular case for example. For completeness, let us briefly explain it. Since $c'$ is regular, a quasi-parabolic bundle being $c$-stable does not depend on $c$ as long as $c$ is in an open neighborhood of $c'$ that does not intersect any wall. So we may assume that $c$ is chosen so that $c_{j_i}=c_i'$ and $c_k=c_l'$. Then $W'=(V, F_\bullet', c')$ is a parabolic subbundle of $W=(V, F_\bullet, c)$. Since $c, c'$ are sufficiently close, $\mu(W')< \mu(W)\leq \mu(W')+\epsilon$ for a small enough $\epsilon$. Therefore, a destabilizing subobject of $W'$ will also destabilize $W$. On the other hand, a destabilizing subobject of $W$ produces a destabilizing subobject of $W'$ after forgetting part of the flag. 

 To identify the fibers of the morphism considered, note that the flag $F_\bullet$ contains for each $i=0, \ldots, k$ the flag
\begin{equation}
\label{eq: sequencesubspaces}F_i'=F_{j_i}\subseteq F_{j_i+1}\subseteq F_{j_i+2}\subseteq \ldots\subseteq F_{j_{i+1}-1}\subseteq F_{j_{i+1}}=F_{i+1}'\,.
\end{equation}
Given fixed $F_i', F_{i+1}'$, choosing the intermediate subspaces is the same as choosing a point in
\[\mathsf{Flag}(F'_{i+1}/F'_{i}; f_{j_i+1}-f_{j_i}, f_{j_i+2}-f_{j_i}, \ldots, f_{j_{i+1}-1}-f_{j_i})\]
obtained from \eqref{eq: sequencesubspaces} by quotienting by $F_i'$.\qedhere
\end{proof}

A particular case of the proposition above is the following: if $\gcd(r,d)=1$ and $c=(c_1, \ldots, c_l)$ is such that $c_l-c_1$ is sufficiently small, then we have a map
\[\pi\colon M_{r,d, f_\bullet}(c)\to M_{r,d}\]
to the moduli space of stable bundles $M_{r,d}$ which forgets the flag. Its fiber over $V\in M_{r,d}$ is 
\[\mathsf{Flag}(V; f_\bullet)\,.\]

\subsection{Maps of vertex algebras}\label{subsec: mapsVA}

To make use of the flag bundle structures described before, we will introduce a morphism of vertex algebras given by weight restriction. 

Let $I'\subseteq I\subseteq [0,1]$ be two sets of weights. Throughout this section, we will emphasize the role of the weights and denote the moduli stack of parabolic bundles with weights in $I, I'$ by $\CM^I, \CM^{I'}$ instead of $\CM^\paar$. Similarly, we denote by $\V^I, \V^{I'}$ the vertex algebras $H_\ast(\CM^I),H_\ast(\CM^{I'})$. 

 Given a parabolic bundle $W=(V, F)$ with weights in $I$, we have an induced parabolic bundle $W_{I'}=(V, F_{|I'})$ with weights in $I'$. Restriction of weights defines a morphism of stacks 
\[\pi=\pi_{I, I'}\colon \CM^I\to \CM^{I'}\,.\]
For example, if $I'=\{0,1\}$ this morphism is the projection of the stack of parabolic bundles with weights in $I$ onto the stack of bundles that forgets all the flag data. In general, the fibers of $\pi$ are products of flag varieties. We consider the vector bundle $\Xi=\Xi_{I, I'}$ on $\CM^I\times \CM^I$ defined by
\[\Xi=\bigoplus_{t\in I\setminus I'}\mathcal{H}om(\CF_1(t), \partial \CF_2(t))\,,\]
where $\CF_1(t), \CF_2(t)$ are the pullbacks of $\CF(t)$ by the two projections onto $\CM^I$. Note that $\Xi$ has a different rank at each connected component of $\CM^I$. It is immediate from \eqref{eq: extcomplex} that we have the following equality in the $K$-theory of $\CM^I\times \CM^I$:
\begin{equation}\label{eq: xiext}\Ext^I_{12}=(\pi\times \pi)^\ast \Ext^{I'}_{12}-\Xi_{I, I'}\,.
\end{equation}
\begin{proposition}[{\cite[Theorem 2.12]{GJT}}]
Define the linear morphism
\[\Omega=\Omega_{I, I'}\colon \V^I=H_\ast(\CM^I)\to H_\ast(\CM^{I'})=\V^{I'}\]
by
\[\Omega(u)=\pi_\ast\big(u\cap c_{\textup{top}}(\Delta^\ast \Xi)\big)\big)\,,\]
where $\Delta^\ast \Xi$ is the restriction of $\Xi$ along the diagonal map $\Delta$. Then $\Omega$ is a morphism of vertex algebras. 
\end{proposition}

We will also denote by $\Omega\colon \widecheck \V^I\to \widecheck \V^{I'}$ the induced morphism of Lie algebras.
\begin{lemma}\label{lem: compositionomega}
Let $I''\subseteq I'\subseteq I$. Then we have
\[\Omega_{I, I''}=\Omega_{I', I''}\circ \Omega_{I, I'}\,.\]
\end{lemma}
\begin{proof}
We observe that we have an equality
\[\Xi_{I, I''}=\Xi_{I, I'}+\big(\pi_{I, I'}\times \pi_{I, I'}\big)^\ast \Xi_{I', I''}\,\]
in the $K$-theory of $\CM^I\times \CM^I$. The statement of the Lemma then follows from an easy application of the push-pull formula. 
\end{proof}

\begin{remark}
A parabolic bundle with weights in $I'$ is also automatically a parabolic bundle with weights in $I\supseteq I'$, so we have an open and closed embedding of stacks $\CM^{I'}\hookrightarrow \CM^I$. Thus, we can regard $\V^{I'}\subseteq \V^I$ as a vertex subalgebra, and $\Omega_{I, I'}$ as a vertex algebra endomorphism $\V^I\to \V^I$.
\end{remark}

\subsection{Joyce style classes for partial parabolic bundles}
\label{subsec: joyceclassespartial}
The morphisms of vertex algebras $\Omega$ above can be used to define classes in $\widecheck \V^I$ for partial parabolic bundles, as we did in Section \ref{subsec: joyceclassesfull} for full parabolic bundles. Before we make the definition, we explain the relation between the flag bundles on Proposition \ref{prop: flagbundle} and the morphisms of vertex algebras $\Omega$.

\begin{proposition}\label{prop: flagbundleclass}
Suppose we are in the conditions of Proposition \ref{prop: flagbundle} and suppose further that $c$ is so that $c_{j_i}=c_i'$. Let $I'=\{0, c'_{1}, \ldots, c_k', 1\}$. Then we have the equality
\begin{equation}\label{eq: flagbundleclass}[M_{r,d, f_\bullet'}(c')]=\frac{\prod_{j=0}^k(f'_{j+1}-f'_{j})!}{\prod_{i=0}^l(f_{i+1}-f_i)!}\Omega_{I, I'}\big([M_{r,d, f_\bullet}(c)]\big)
\end{equation}
in $\widecheck \V^{I'}$.
\end{proposition}

\begin{proof}
Let $\alpha\in C(I)$, $\alpha'\in C(I')$ be the topological types associated to the two moduli spaces. Note that $\Omega_{I, I'}(\alpha)=\alpha'$ by the assumption. By Proposition \ref{prop: flagbundle}, we have a pullback square
\begin{center}
\begin{tikzcd}
M_{r,d,f_\bullet}(c)\arrow[r, "\pi"]\arrow[d] &M_{r,d,f'_\bullet}(c')\arrow[d]\\ 
(\CM_\alpha^I)^\rig \arrow[r, "\pi_{I, I'}"] &(\CM_{\alpha'}^{I'})^\rig
\end{tikzcd}
\end{center}
By \eqref{eq: xiext}, the descent of $\Delta^\ast \Xi$ to $M_{r,d,f_\bullet}(c)$ is the relative tangent bundle of the morphism $\pi_{I, I'}$. Hence, the pushforward of $c_{\textup{top}}(\Delta^\ast \Xi)$ along the smooth morphism $\pi_{I, I'}$ is equal to the Euler characteristic of the fibers, which is precisely the constant 
\[\frac{\prod_{i=0}^l(f_{i+1}-f_i)!}{\prod_{i=0}^k(f_{j_{i+1}}-f_{j_i})!}\,\]
according to the description of the fibers in Proposition \ref{prop: flagbundle}.\qedhere
\end{proof}

\subsubsection{Extension to non-regular weights}
\label{subsubsec: classesnonregular}
Finally, we would like to define
\[[M_{r,d,f_\bullet}(c)]\in \widecheck \V^I\]
for any $c\in S_{r, d, f_\bullet}$, not necessarily regular. 

Note that the definition of the wall-crossing coefficients $U(\vec J, \vec d; c, c')$ makes perfect for $c,c$ in the compactification $\overline S_{r,d}$. Hence, we can use the wall-crossing formula to extend the definition of the classes $[M_{r,d}^{\textup{full}}(c)]$ to the boundary of the space of stability conditions.

\begin{definition}\label{def: joyceclassesboundary}
Let $c\in S_{r,d, f_\bullet}$ and $\overline c\in \overline S_{r,d}$ be as in \eqref{eq: overlinec}. Choose some $\widetilde c\in S_{r,d}$ such that $\widetilde c_{f_i}=\overline c_{f_i}=c_i$. We define the class $[M_{r,d}^{\textup{full}}(\overline c)]$ by imposing the wall-crossing formula \eqref{eq: fullwcformula}:
\[[M_{r,d}^{\textup{full}}(\overline c)]
=\sum_{\substack{\vec J\vdash [r]\\
\vec d\vdash d}}U(\vec J, \vec d; \widetilde c, \overline c)\big[\big[\ldots\big[[M_{r_1,d_1}^{\textup{full}}(\widetilde c_{J_1})], [M_{r_2,d_2}^{\textup{full}}(\widetilde c_{J_2})]\big],\ldots,\big], [M_{r_m,d_m}^{\textup{full}}(\widetilde c_{J_m})]\big]
\]
Then, we let $I'=\{0, c_1, \ldots, c_l, 1\}$ and define
\[[M_{r,d, f_\bullet}(c)]=\frac{1}{\prod_{i=0}^l(f_{i+1}-f_i)!}\Omega_{I, I'}([M_{r,d}^{\textup{full}}(\overline c)])\in \widecheck \V^{I'}\subseteq \widecheck \V^{I}\,.\]
\end{definition}

As in the proof of Proposition \ref{prop: welldefinedclass}, the class $[M_{r,d}^{\textup{full}}(\overline c)]$ does not depend on the choice of $\widetilde c$ after forgetting weights. Since $\Omega_{I, I'}$ forgets all the weights in $I\setminus I'$ and we require $\widetilde c_{f_i}=c_i$, the class $[M_{r,d, f_\bullet}(c)]$ is independent of the choice of $\widetilde c$ and hence it is well defined. If $c\in S_{r,d, f_\bullet}$ is regular, then we may choose $\widetilde c\in S_{r,d}$ in the same chamber as $\overline c$, in which case the wall-crossing constants $U(\vec J, \vec d; \widetilde c, \overline c)$ are 0 except for the trivial partitions, for which it is 1. Thus, Proposition \ref{prop: flagbundleclass} implies that Definition~\ref{def: joyceclassesboundary} is consistent with the previous one in Section \ref{subsec: liealgebra} when $c$ is regular.

%\begin{example}\mcomment{Finish}
%We consider the case $r=3$ and $d=0$. There is only one wall intersecting $S_{3,0}$, which is
%\[W_{\{2\}, \{1,3\}, 0,0}=\{\lambda_1=\lambda_2\}\,.\]
%However, there are two more walls which intersect $\overline{S}_{3,0}$ at $(\lambda_1, \lambda_2)=(0,0)$, namely
%\[W_{\{1\}, \{2,3\}, 0,0}\textup{ and }W_{\{3\}, \{1,2\}, 0,0}\,.\]
%\end{example}

%Let $\alpha=(r,d,f)\in C(\CP_I)$ and let $f_\bullet=(f_1, \ldots, f_l)=\textup{im}(f)$ and $c$ be the corresponding weights. Given a subset $J=\{j_1, \ldots, j_k\}\subseteq [l]$, let $I'=\{c_{j_1},\ldots, c_{j_k}\}$. Then the map of stacks $\sigma_{I, I'}$ sends a parabolic bundle $(V, F)$, with corresponding flag
%\[0\subseteq F_1\subseteq F_2\subseteq \ldots \subseteq F_l\subseteq V_p\]
%with weights $c$ to the parabolic bundle with underlying flag
%\[0\subseteq F_{j_1}\subseteq F_{j_2}\subseteq \ldots \subseteq F_{j_k}\subseteq V_p\]
%and weights $c'=(c_{j_1},\ldots, c_{j_k})$. Let $f'=(f_{j_1}, \ldots, f_{j_k})$ and $\alpha'=(r,d,f')$. For convenience, we also set $j_0=0$ and $j_{k+1}=l+1$. 

%\begin{proposition}\label{prop: flagbundle}
%Suppose that $\alpha'$ is regular. Then there is some sufficiently small $\varepsilon>0$ such that, if $c_{i+1}-c_i<\varepsilon$ for every $i\in [l]\setminus J$, then $W$ of type $\alpha$ is stable if and only if $W_{|I'}$ is stable. In particular, the map $\sigma^\rig$ restricts to
%\[\sigma\colon M_{r,d, f_\bullet}(c)\to M_{r,d, f_\bullet'}(c')\,.\]
%\end{proposition}

\subsection{Comparison with Joyce's definition for bundles}\label{subsec: comparisonclasses}

In particular, Definition \ref{def: joyceclassesboundary} with $f_\bullet=\emptyset$ determines a class $[M_{r,d}]$ in the Lie algebra $\widecheck \V^{I'}$ with $I'=\{0,1\}$, hence providing a way to define descendent integrals on the moduli of vector bundles when $\gcd(r,d)\neq 1$, and thus we are in the presence of semistable objects. The definition is similar in spirit to the one provided by Joyce's theory \cite{Jo21}, but not trivially equivalent. We now show that the two agree.

A key input in Joyce's construction of the classes $[M_{r,d}]$ is a choice of a ``framing functor'' $\Phi$ (or a collection of such functors) from the category of coherent sheaves on $C$ to vector spaces, with some properties (see \cite[Assumption 5.1 (g)]{Jo21} for precise definitions). Loosely speaking, $\Phi$ is required to be exact when restricted to some exact subcategory containing the semistable objects. Given such $\Phi$, one constructs an auxiliary moduli space, sometimes referred to as moduli of Joyce--Song pairs, that parametrizes a vector bundle $V$ together with a section $\BC\to \Phi(V)$, and an appropriate stability condition that does not admit strictly semistable objects; see \cite[Example 5.6]{Jo21}. One of the fundamental aspects of Joyce's theory is that actually the definition of such classes does not depend on the choice of $\Phi$. 

A typical choice of $\Phi$, considered for example in \cite[Section 7.6]{Jo21} or \cite{bu}, is to take $\Phi(V)=H^0(V\otimes \CO_C(N))$ for $N$ sufficiently large. To obtain the relation with parabolic bundles, we instead consider the framing functor $\Phi(V)=H^0(V_p^\vee)$. Note that giving a section for this framing functor is the same as choosing a 1-dimensional subspace of $V_p^\vee$ or, equivalently, a $(r-1)$-dimensional subspace of $V_p$. Indeed, it is not difficult to see that the Joyce--Song moduli space associated to this framing functor is precisely
\[M_{r,d, (r-1)}(c_1, c_1+\epsilon)\]
where $\epsilon>0$ is sufficiently small. We will just write $M_{r,d, (r-1)}(\epsilon)$ for the moduli space above. The claim that our definition agrees with Joyce's amounts to the following wall-crossing formula.

\begin{theorem}\label{thm: agreement}
Let $\epsilon>0$ be sufficiently small and $I'=\{0,1\}$, $I''=\{0,c_1, c_1+\epsilon, 1\}$. We have the wall-crossing formula
\[\Omega_{I'', I'}([M_{r,d, (r-1)}(\epsilon)])=\sum_{\substack{\vec r\vdash r\,,\, \vec d\vdash d\\
d_i/r_i=d/r}}\frac{(-1)^{m+1}r_1}{m!}\big[\big[\ldots \big[[M_{r_1, d_1}], [M_{r_2, d_2}]\big], \ldots, \big], [M_{r_n, d_n}]\big]\,\] 
in $\widecheck \V^{I'}$. 
In particular, our class $[M_{r,d}]$ matches the one defined in~\cite{Jo21}. 
\end{theorem}

Before we give the proof, let us say a few words about how it will go; this proof also illustrates the general strategy to prove wall-crossing formulas for partial parabolic bundles. First, let 
\[\overline 0=(0, 0, \ldots, 0, 0)\,,\, \overline \epsilon=(0, 0, \ldots, 0, \epsilon) \in \overline S_{r,d}\,.\]
Then we have a wall-crossing formula between the classes $[M_{r,d}^\full(\overline 0)]$ and  $[M_{r,d}^\full(\overline \epsilon)]$ in $\widecheck \V^I$; this follows from the definition of the classes and the fact that the wall-crossing coefficients are well-behaved under composition (c.f. Proposition \ref{prop: welldefinedclass}). We will be able to determine the wall-crossing coefficients explicitly. Then, Theorem~\ref{thm: agreement} will be obtained by applying $\Omega_{I, I'}$ to this wall-crossing formula, after some combinatorial work with the coefficients. 

\begin{proof}[Proof of Theorem \ref{thm: agreement}]Given $\vec J\vdash [r]$ and $\vec d\vdash d$ partitions of length $m$, we begin by calculating the wall-crossing coefficients 
\[S(\vec J, \vec d)\coloneqq S(\vec J, \vec d; \overline 0, \overline \epsilon)\textup{ and }U^\ast(\vec J, \vec d)\coloneqq U^\ast(\vec J, \vec d; \overline 0, \overline \epsilon)\,.\]
We denote as usual $r_i=|J_i|$; we have
\[\mu_{\overline 0}(J_i, d_i)=\frac{d_i}{r_i}\textup{ and }\mu_{\overline \epsilon}(J_i, d_i)=\begin{cases}\frac{d_i+\epsilon}{r_i}& \textup{ if }r-1\in J_i\\
\frac{d_i}{r_i}& \textup{ otherwise.}\end{cases}\]
One can see (either directly or applying for example \cite[Proposition 3.14]{Jo21}) that if $\epsilon$ is small enough then $S(\vec J, \vec d)=0$ unless $d_i/r_i=d/r$ for $i=1, \ldots, m$. Hence, $S(\vec J, \vec d)$ can be non-zero only when we have condition (a) from Section \ref{subsubsec: wccoefficients} for each $i=1, \ldots, m-1$, so we have
\[S(\vec J, \vec d)=\begin{cases}
(-1)^{m-1}&\textup{ if }d_i/r_i=d/r\textup{ and }r-1\in J_1\\
0 &\textup{ otherwise.}
\end{cases}\]
From here, it follows (recall the definition of $U^\ast$, \cite[(3.3)]{Jo21}) that $U^\ast(\vec J, \vec d)$ in non-zero only when $d_i/r_i=d/r$, in which case
\[U^\ast(\vec J, \vec d)=\sum_{\ell\geq 1}(-1)^{\ell-1}\sum_{\substack{k_1+\ldots+k_\ell=m\\
k_1\geq j\,,\,k_i\geq 1}}\prod_{t=1}^\ell \frac{1}{k_t!}\,,\]
where $j$ is the unique index such that $r-1\in J_j$. The wall-crossing formula is\footnote{We remind the reader that the classes $[M_{r, d}^\full(\overline 0)]$, $[M_{r_i, d_i}^\full(\overline 0_{J_i})]$ are only well defined up to choice of weights, and the equality should be interpreted after forgetting weights. Because of this, we pedantically write $M_{r_i, d_i}^\full(\overline 0_{J_i})$ instead of just $M_{r_i, d_i}^\full(\overline 0)$ to emphasize that, even if $r_i=r_{i'}$ and $d_i=d_{i'}$, the corresponding classes are the same only after forgetting weights, and in particular they do not necessarily commute. After we apply $\Omega_{I, I'}$ this subtlety will be entirely removed.}
\[[M_{r,d}^\full(\overline \epsilon)]=\sum_{\substack{\vec J, \vec d\\d_i/r_i=d/r}}U^\ast(\vec J, \vec d) [M_{r_1, d_1}^\full(\overline 0_{J_1})]\ast\ldots \ast [M_{r_m, d_m}^\full(\overline 0_{J_m})]\]
in the enveloping algebra of $\widecheck \V^I$. We now apply the morphism of Lie algebras $\Omega_{I, I'}\colon \widecheck \V^I\to \widecheck \V^{I'}$, which induces a homomorphism of the corresponding enveloping algebras. By Lemma \ref{lem: compositionomega} and Proposition \ref{prop: flagbundleclass}, the left hand side becomes
\[\Omega_{I, I'}\big([M_{r,d}^\full(\overline \epsilon)]\big)=\Omega_{I'', I'}\big(\Omega_{I, I''}\big([M_{r,d}^\full(\overline \epsilon)]\big)\big)=(r-1)!\Omega_{I'', I'}\big([M_{r,d, (r-1)}(\epsilon)]\big)\,.\]
By definition of the classes,
\[\Omega_{I, I'}\big([M_{r_i, d_i}^\full(\overline 0_{J_i})]\big)=r_i![M_{r_i, d_i}]\,.\]
Therefore, applying $\Omega_{I, I'}$ to the wall-crossing formula expresses \[(r-1)!\Omega_{I'', I'}\big([M_{r,d, (r-1)}(\epsilon)]\big)\] as
\[\sum_{\substack{\vec r, \vec d\\
d_i/r_i=d/r}}\Big(\prod_{i=1}^m r_i!\Big)c(\vec r)[M_{r_1, d_1}]\ast \ldots \ast [M_{r_m, d_m}]\]
where the combinatorial coefficients $c(\vec r)$ are the sum of $U^\ast(\vec J, \vec d)$ over all partitions $\vec J\vdash [r]$ with $|J_i|=r_i$. Note that, if we fix $j$, the number of such partitions with $r-1\in J_j$ is given by
\[\binom{r-1}{r_1, \ldots, r_j-1, \ldots, r_m}=\frac{r_j}{r}\binom{r}{r_1, \ldots, r_m}\,.\]
Hence, 
\begin{align*}c(\vec r)&=\frac{1}{r} \binom{r}{r_1, \ldots, r_m}\sum_{j=1}^m r_j\sum_{\ell\geq 1}(-1)^{\ell-1}\sum_{\substack{k_1+\ldots+k_\ell=m\\
k_1\geq j\,,\,k_i\geq 1}}\prod_{t=1}^\ell \frac{1}{k_t!}\\
&=\frac{(-1)^m(r-1)!}{m!r_1!\ldots r_m!} \sum_{j=1}^m (-1)^jr_j\binom{m-1}{j-1}\,,
\end{align*}
where in the second line we used an identity that we shall prove in Lemma \ref{lem: identity}. Thus, we conclude that
\begin{align*}\Omega_{I'', I'}&\big([M_{r,d, (r-1)}(\epsilon)]\big)\\
&=\sum_{\substack{\vec r, \vec d\\
d_i/r_i=d/r}}\frac{(-1)^m}{m!}\left(\sum_{j=1}^m (-1)^j r_j \binom{m-1}{j-1}\right)[M_{r_1, d_1}]\ast \ldots \ast [M_{r_m, d_m}]\\
&=\sum_{\substack{\vec r\,,\, \vec d\\
d_i/r_i=d/r}}\frac{(-1)^{m+1}r_1}{m!}\big[\big[\ldots \big[[M_{r_1, d_1}], [M_{r_2, d_2}]\big], \ldots, \big], [M_{r_n, d_n}]\big]
 \end{align*}
 where the last line uses the combinatorial identity \cite[(9.54)]{Jo21}. 
 
 To conclude that the classes $[M_{r,d}]$ agree with the ones defined by Joyce, note that the wall-crossing formula just proved is exactly Theorem 5.7 (iii) in loc. cit. for the choice of framing functor discussed before. The class $[M_{r,d, (r-1)}(\epsilon)]$ agrees with Joyce's, since this moduli space does not have strictly semistable objects and hence it is the pushforward of the virtual fundamental class. Since the right hand side has main term $r[M_{r,d}]$ and the remaining terms only involve classes in lower ranks, this formula determines all the classes $[M_{r,d}]$ recursively, establishing the match between the two definitions.\qedhere 
\end{proof}

\begin{lemma}\label{lem: identity}
Let $1\leq j\leq m$ be integers. We have the following identity:
\[\sum_{\ell\geq 1}(-1)^{\ell-1} \sum_{\substack{k_1+\ldots+k_\ell=m\\
k_1\geq j\,,\,k_i\geq 1}}\binom{m}{k_1, k_2, \ldots, k_\ell}=(-1)^{m+j}\binom{m-1}{j-1}\]
\end{lemma}
\begin{proof}
Denote 
\[F(m, j)=\sum_{\ell\geq 1}(-1)^{\ell-1+m+j} \sum_{\substack{k_1+\ldots+k_\ell=m\\
k_1\geq j\,,\,k_i\geq 1}}\binom{m}{k_1, k_2, \ldots, k_\ell}\,.\]
It is straightforward to see that
\begin{equation}\label{eq: recursion}F(m, j)+F(m, j+1)=\binom{m}{j}F(m-j, 1)\,.\end{equation}
We now argue by induction on $m$ that $F(m, j)=\binom{m-1}{j-1}$. If $j=m$, this is immediate. Otherwise, use the recursion \eqref{eq: recursion}, together with the induction hypothesis that $F(m-j, 1)=1$, to show that the claim for $j+1$ implies the claim for $j$.\qedhere
\end{proof}

\section{Hecke modifications}\label{sec: hecke}

Hecke modifications are certain automorphisms of the category of parabolic bundles that can be used to identify moduli spaces of parabolic bundles with different degrees. In this section, we fit the theory of Hecke modifications together with the vertex algebra structure from Section \ref{sec: vertexalgebras}. 

\begin{remark}
What we call Hecke modification here is called tautological Hecke correspondence in \cite{ST} and parity shift in \cite{DP}. The term Hecke correspondence in \cite{DP} is reserved for something different. 
\end{remark}

\subsection{Definition of Hecke modifications}
\label{subsec: heckedefinition}
The definition of Hecke modification is made cleaner by describing the category of parabolic bundles in an alternative way.\footnote{The description below is the general way to define parabolic sheaves on $(X, D)$ where $X$ is a non-singular projective variety and $D$ is a divisor, see \cite{MY}.} Given a parabolic bundle $W=(V, F)$, define for each $t\in I$ the vector bundle $W(t)$ as the kernel
\[0\to W(t)\to V\to \iota_\ast\big(V_p/F(t))\to 0\,.\]
In particular $W(0)=V\otimes \,\!\CO_C(-p)$ and $W(1)=V$. We extend $W(t)$ to every $t\in I+\BZ$ by requiring that $W(t+1)=W(t)\otimes \CO_C(p)$. Such a function is increasing, meaning that for any $t<t'$ we have a distinguished injection of locally free sheaves $W(t)\hookrightarrow W(t')$, lower semicontinuous, and periodic in the sense that $W(t+1)=W(t)\otimes \CO_C(p)$.

This construction defines an equivalence between the category of parabolic bundles with weights in $I$ and the category whose objects are increasing, lower semicontinuous and periodic functions
\[W\colon I+\BZ\to \{\textup{vector bundles over }C\}\,.\]

\begin{remark}\label{rem: slope}
Under this interpretation of parabolic bundles, the slope admits the natural looking formula
\[\mu(W)=\frac{\int_{0}^1 \deg(W(t)) \,\textup{d}t}{\rk(V)}\,.\]
\end{remark}

For the construction of Hecke modifications, we will assume for convenience that $I=[0,1]$ or $I=\frac{1}{N}\BZ\cap [0,1]$, so that $I+\BZ\subseteq \BR$ is a subgroup.\footnote{For general $I$, the Hecke modification would not be an automorphism of the exact category $\CP_I$, but rather an isomorphism between categories $\CP_I$ and $\CP_{I'}$ for a different set of weights $I'$, obtained by ``translating $I$ modulo 1''.} Given $\tau\in \frac{1}{N}\BZ$, the Hecke modification is the automorphism of the  category $\CP_I$ that sends a parabolic bundle $W=(V, F)$ to a new parabolic bundle $\widetilde W=(\widetilde V, \widetilde F)$ such that 
\[\widetilde W(t)\coloneqq W(t-1+\tau)\,.\]

We now make it more explicit how the Hecke modification transforms $V$ and $F$ when $\tau\in [0,1]$. The underlying vector bundle of $\widetilde W$ is $\widetilde V=\widetilde W(1)=W(\tau)$. It follows from the definition of $\widetilde V$ that there is a map $\widetilde V\to V$ and that the restriction $\widetilde V_p\to V_p$ has image $F(\tau)$. We denote by $\phi\colon \widetilde V_p\to F(\tau)$ the surjection onto the image, and let $K=\ker(\phi)$. 

If $0\leq t\leq 1-\tau$ then $\widetilde F(t)$ is the kernel of the following composition:
\[\widetilde F(t)=\ker\big(\widetilde V_p\to V_p\to V_p/F(t+\tau)\big)\]
Note that since $F(\tau)\subseteq F(t+\tau)$ this composition is surjective. Clearly $\widetilde F(t)$ is increasing and $\widetilde F(1-\tau)=K$. 

If $1-\tau\leq t\leq 1$, we have $F(t+\tau-1)\subseteq F(\tau)$ and
\[\widetilde F(t)=\phi^{-1}(F(t+\tau-1))\,.\]
Note that $\widetilde F(t)$ is increasing and that $\widetilde F(1-\tau)=K$. 

In particular, the Hecke modification sends an object of topological type $\alpha=(r,d,f)$ to an object of topological type $h_\tau(\alpha)$, where $h_\tau\colon K(I)\to K(I)$ is given by
\[h_\tau(r,d,f)=(r, d-r+f(\tau), \widetilde f)\]
with
\[\widetilde f(t)=\begin{cases}
f(t+\tau)-f(\tau)&\textup{ if }0\leq t\leq 1-\tau\\
f(t+\tau-1)-f(\tau)+r&\textup{ if }1-\tau \leq t\leq 1\,.
\end{cases}\]
\bigskip

\includegraphics[scale=0.8, trim=1.5cm 0 -2cm 0cm]{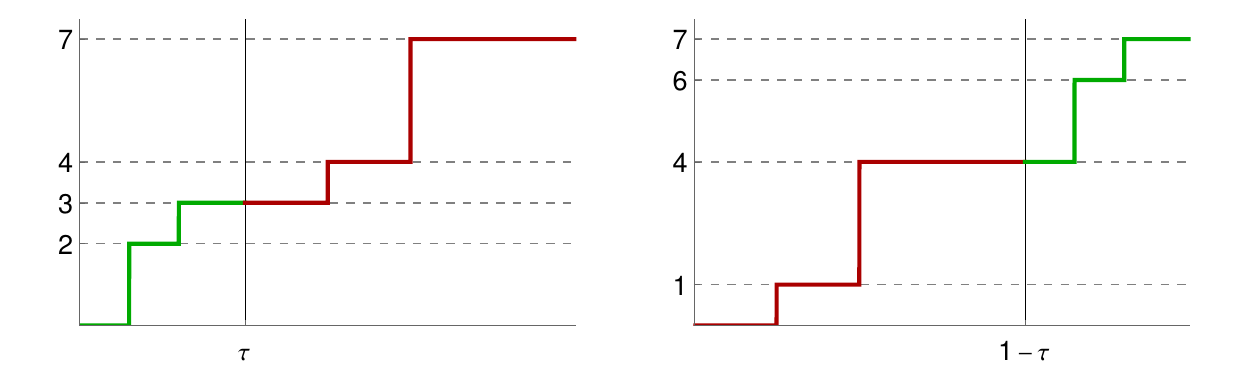}
\captionof{figure}{An example of the graphs of functions $f$ (on the left) and $\widetilde f$ (on the right). In this example, $r=7$ and $\tau$ is between $c_2$ and $c_3$. We have $h_\tau(7, d, f)=(7, d-4, \widetilde f)\,.$}
\setlength{\parindent}{1em}
\bigskip

The Hecke modification defines an automorphism of the stack $\CM^\paar$
\[h_\tau\colon \CM^\paar\to \CM^\paar\]
and hence induces an automorphism of the vertex algebra
\[(h_{\tau})_\ast\colon H_\ast(\CM^\paar)=\Vpar\to \Vpar=H_\ast(\CM^{\paar})\,.\]
We also denote $(h_{\tau})_\ast\colon \widecheck \V^\paar\to \widecheck \V^\paar$ the induced homomorphism of Lie algebras. We call $(h_\tau)_\ast$ a Hecke operator. 

We may define also the Hecke operator at the level of the descendent algebra, and hence on $\Vpartr$. We let $h^\dagger_\tau\colon \BD^\paar\to \BD^\paar$ be the algebra homomorphism defined on generators by
\begin{align*}
h^\dagger_\tau\ch_k(e(t))=\begin{cases}
\ch_k\big(e(t+\tau)-e(\tau) \big)& \textup{ if }0\leq t\leq 1-\tau\\
\ch_k\big(e(t+\tau-1)+e(1)-e(\tau)\big) & \textup{ if }1-\tau\leq t\leq 1
\end{cases}
\end{align*}
and 
\begin{align*}
h^\dagger_\tau\ch_k(\gamma)=\begin{cases}
\ch_k(\gamma) & \textup{ if }\gamma\in H^{>0}(C)\\
\ch_k(\gamma) + \ch_{k-1}(e(\tau)-e(1))& \textup{ if }\gamma=\1
\end{cases}
\end{align*}

It is clear that $h_\tau^\dagger$ descends to
\[h_\tau^\dagger\colon \BD^\paar_{h_\tau\alpha}\to \BD^\paar_{\alpha}\,.\]
By taking duals, we obtain an operator
\[(h_\tau)_\dagger\colon \Vpartr\to \Vpartr\,.\]

The next lemma establishes the relation between this formally defined operator and the Hecke operator $(h_\tau)_\ast$.

\begin{lemma}\label{lem: hecketr}
The following diagram commutes:
\begin{center}
\begin{tikzcd}
\Vpar\arrow[r, "(h_\tau)_\ast"]\arrow[d, hookrightarrow]& \Vpar\arrow[d, hookrightarrow]\\
\Vpartr\arrow[r, "(h_\tau)_\dagger"]& \Vpartr
\end{tikzcd}
\end{center}
\end{lemma}

\begin{proof}
The lemma is equivalent to the commutativity of the dual square 
\begin{center}
\begin{tikzcd}
H^\ast(\CM^{\paar}_\alpha)& H^\ast(\CM^{\paar}_{h_\tau \alpha}) \arrow[l, swap,"h_\tau^\ast"]\\
\BD^\paar_{\alpha}\arrow[u,"\xi_{\CW}"]&\BD^\paar_{h_\tau\alpha}\arrow[u,"\xi_{\CW}"] \arrow[l, swap,"h_\tau^\dagger"]
\end{tikzcd}
\end{center}
where $\xi$ is the geometric realization homomorphism with respect to the universal parabolic bundle $\CW=(\CV, \CF)$. By the previous description of $h_\tau$, we obtain
\[(h_\tau\times \id)^\ast \CV=\CV-\iota_\ast(\CV_{|\CM^\paar\times \{p\}})+\iota_\ast(\CF(\tau))\]
in the $K$-theory of $\CM^\paar\times C$, where $\iota\colon \CM^\paar=\CM^\paar\times \{p\}\hookrightarrow \CM^\paar\times C$ is the inclusion of the fiber over $p\in C$. We also find that
\begin{align*}
h^\ast_\tau \CF(t)=\begin{cases}
\CF(t+\tau)-\CF(\tau) & \textup{ if }0\leq t\leq 1-\tau\\
\CF(t+\tau-1)+\CF(1)-\CF(\tau) & \textup{ if }1-\tau\leq t\leq 1\,,
\end{cases}
\end{align*}
and the conclusion follows.\qedhere
\end{proof}

\begin{remark}\label{rem: compositionhecke}
It is immediate to see from the definition that the composition of two Hecke operators in a Hecke operator; precisely,
\[h_{\tau}\circ h_{\tau'}=h_{\tau+\tau'-1}\,.\]
In particular, note that the $N$-fold composition of $h_{1-1/N}$ is $h_{-1}$, which is the auto-equivalence $\CO_C(-p)\otimes -$.
\end{remark}

\begin{remark}
Under the equivalence between the categories of parabolic bundles and of bundles on the root stack $\mathcal C=\sqrt[N]{(C,p)}$, the Hecke operator $h_{1-1/N}$ corresponds to twisting by the orbifold line bundle $\CO_\CC(-\frac{1}{N}p)$, which explains the periodicity in the previous remark.
\end{remark}

\subsection{Hecke modifications and stability}
\label{subsec: heckestability}
The Hecke modifications preserve stability, as shown in \cite[Proposition 7.1]{ST}. For the reader's convenience, we give our own proof in the slightly more general setting that we consider.

\begin{lemma}\label{lem: heckestability}
The Hecke modification $h_\tau$ sends (semi)stable objects in $\CP_I$ to (semi)
stable objects in $\CP_I$.

Hence, for any regular $\alpha\in C(I)$ we have an identity
\begin{equation}
\label{eq: heckeclass}(h_\tau)_\ast [M_\alpha]=[M_{h_\tau \alpha}]\,.
\end{equation}
\end{lemma}
\begin{proof}
This follows from $h_\tau$ being an automorphism of the exact category $\CP_I$ together with the fact that $h_\tau$ changes the slope by a constant, which we show here. Using the expression for the slope in Remark \ref{rem: slope},
\begin{align*}\mu(\widetilde W)&=\frac{1}{r}\int_0^1 \dim \widetilde W(t)\,\mathrm{d}t\\
&=\frac{1}{r}\left(\int_{0}^{1-\tau}(\dim W(t+\tau)-r)\,\mathrm{d}t+\int_{1-\tau}^{1}\dim W(t-1+\tau)\,\mathrm{d}t\right)\\
&=\frac{1}{r}\int_{0}^{1}\dim W(t)\,\mathrm{d}t-(1-\tau)=\mu(W)-1+\tau\,.\qedhere
\end{align*}
\end{proof}

We remark that \eqref{eq: heckeclass} also holds for non-regular $\alpha$. This is a consequence of the equality for regular $\alpha$ and the wall-crossing formula; note that the wall-crossing coefficients are compatible with the Hecke operator, since the Hecke modification changes the slope monotonically, by the proof of the previous lemma.

%The same holds true even if $\alpha$ is not regular:

%\begin{proposition}
%For any $\alpha\in C(I)$, not necessarily regular, we have
%\[(h_\tau)_\ast [M_\alpha]=[M_{h_\tau \alpha}]\,.\]
%\end{proposition}
%\begin{proof}
%As explained above, this holds true for $\alpha$ regular. Suppose now that $\alpha$ is an arbitrary topological type of a full parabolic bundle and let $r,d,c$ be the rank, degree and weights of $\alpha$, respectively. \mcomment{write}
%\end{proof}

\subsection{Hecke symmetries between spaces of stability conditions}
\label{subsec: heckespacesstability}
Let $r>0$ and $d, d'\in \BZ$. We can identify the spaces of stability conditions
\[S_{r,d}\simeq S_{r, d'}\,,\]
and their compactifications
\[\overline S_{r,d}\simeq \overline S_{r, d'}\,,\]
via the map
\[\lambda\mapsto \mu=(\lambda_{i+1}, \lambda_{i+2}, \ldots, \lambda_{r-1}, 1-\lambda_1-\ldots-\lambda_{r-1}, \lambda_1, \ldots, \lambda_{i-1})\,\]
where $0\leq i\leq r-1$ is $d'-d \mod r$. This identification preserves the wall-chamber structures. 

Consider a given $\alpha\in C(I)$ with associated vector of weight differences $\lambda\in S_{r,d}$ and $0\leq i\leq r-1$. Let $\tau$ be such that $c_{i-1}<\tau<c_i$. Then $h_\tau(\alpha)$ has rank $r$, degree $d'\coloneqq d-r+i$ and vector of weight differences $\mu \in S_{r,d'}$.

\begin{corollary}\label{cor: heckeiso}
Let $r>0$ and $d, d'\in \BZ$. Let $\lambda\in S_{r,d}$ and let $\mu\in S_{r, d'}$ be the corresponding stability condition. Then there is $\tau\in \BR$ such that $h_\tau$ induces an isomorphism
\[M_{r,d}^{\full}(\lambda)\simeq M_{r,d'}^{\full}(\mu)\,.\]
\end{corollary}

\section{The affine Weyl symmetry}\label{sec: weyl}

The next ingredient in the determination of integrals on moduli of parabolic bundles is the affine Weyl symmetry. This idea was crucially used by Teleman--Woodward \cite{TW} to obtain the formula for the index of certain $K$-theory classes on moduli spaces of bundles, generalizing the Verlinde formula. It appears again in the work of Szenes--Trapeznikova \cite{ST} to prove the Verlinde formula for parabolic bundles, in a more elementary form which is also better suited for intersection numbers. We follow their approach here.

\begin{remark}Teleman and Woodward use a non-abelian localization theorem  to relate $K$-theoretical invariants on the the moduli spaces of semistable bundles and on the moduli stack of all bundles, without any stability. The advantage of removing stability is that the cohomology of the stack of all parabolic bundles admits itself an action of the affine Weyl group. Since intersection numbers cannot be defined on such a stack, this approach does not work without first moving to $K$-theory. 
\end{remark}

\begin{comment}
As explained in Section \ref{subsec: flagbundles}, for certain choices of stability condition $c$ and $\gcd(r,d)=1$, the moduli of full parabolic bundles $M_{r,d}^\full(c)$ if a full flag bundle over $M_{r,d}$. From this description, the cohomology of $H^\ast(M_{r,d}^\full(c))$ comes with an action of the symmetric group $\Sigma_r$ (regarded as the Weyl group of $\mathfrak{gl}_r$) which permutes the first Chern classes of the tautological line bundles; moreover, integrals are anti-symmetric with respect to this action, see Corollary \ref{cor: weylantisymmetry}. 

This action alone is not enough to constraint integrals enough to determine them. The insight from \cite{ST} is that we can use the Hecke isomorphisms to combine two of these $\Sigma_r$ actions and obtain an action of the affine Weyl group, which is much larger and completely determines the integrals. More precisely: by the Hecke isomorphism of Corollary \ref{cor: heckeiso}, given another $d'$ such that $\gcd(r, d')=1$, there is a different stability condition $c'$ such that $M_{r,d}^\full(c')$ is a full flag bundle over $M_{r,d'}$, and thus a second $\Sigma_r$ action. If it were the case that $c, c'$ were in the same chamber of $S_{r,d}$, we would get two different actions of $\Sigma_r$ on the cohomology of $M_{r,d}^\full(c')$. Although this is never the case, wall-crossing tells us the difference between integrals on the two moduli spaces, so it is still possible to combine the two actions to prove Theorem~\ref{thm: reconstruction}.
\end{comment}

\subsection{Weyl symmetry for flag bundles} We start with the elementary observation that full flag bundles carry an action of $\Sigma_r$ on their cohomology. Let $M$ be a projective smooth variety (e.g. $M=M_{r,d}$) and let $\CE$ be a rank $r$ bundle on $M$ (e.g. $\CE=\CV_{|M\times \{p\}}$). By \cite[Théorème 1]{Gro}, the cohomology of the full flag variety $\mathrm{Fl}_M(\CE)\coloneqq \mathrm{Fl}_M(\CE; 1, 2, \ldots, r-1)$ is given by
\begin{equation}\label{eq: cohflagbundle}
H^\ast(\mathrm{Fl}_M(\CE))\simeq \frac{H^\ast(M)[t_1, \ldots, t_r]}{\big(c_i(\CE)=e_i(t_1, \ldots, t_r)\big)_{i=1}^r}
\end{equation}
where $e_i(t_1, \ldots, t_r)$ are the elementary symmetric polynomials in $t_i$. Under this isomorphism, the variables $t_i$ are the first Chern classes of the tautological line bundles obtained as successive quotients $\CF_i/\CF_{i-1}$ in the tautological flag.

\begin{proposition}\label{prop: weylflagbundle}
Let $M, \CE$ be as before. Then $\Sigma_r$ acts on the cohomology $H^\ast(\mathrm{Fl}_M(\CE))$ by permuting the first Chern classes $t_1, \ldots, t_r$ and fixing the cohomology classes pulled back from $M$.

Moreover, we have
\[\int_{\mathrm{Fl}_M(\CE)}\sigma\cdot D=(-1)^{\textup{sgn}(\sigma)}\int_{\mathrm{Fl}_M(\CE)}D\]
for any $D\in H^\ast(\mathrm{Fl}_M(\CE))$ and $\sigma\in \Sigma_r$. 
\end{proposition}
\begin{proof}
The existence of the symmetry is clear from \eqref{eq: cohflagbundle}. The anti-symmetry follows from the following integral formula: if $F(t_1, \ldots, t_r)$ is a polynomial with coefficients being cohomology classes of $M$, then
\[\int_{\mathrm{Fl}_M(\CE)}F(t_1, \ldots, t_r)=\Res_{z_1}\ldots \Res_{z_r} \int_M F(z_1, \ldots, z_r)\prod_{1\leq i<j\leq r}(z_i-z_j)\prod_{i=1}^r \hat s_{z_i}(\CE)\]
where $\hat s_z(\CE)=z^r s_{1/z}(\CE)$ and $s_\ast(\CE)$ is the Segre polynomial of $\CE$. A very similar formula (but taking only residues in $z_1, \ldots, z_{r-1}$) can be found in \cite[Page 3]{DP}. Our formula can be derived from theirs by using the same identity as in the proof of \cite[Theorem 1.1]{DP}.
\end{proof}

In particular, if $\alpha_\bullet=(r,d, f_\bullet)$ is a full type, i.e. $f_\bullet=(1, 2, \ldots, r-1)$, we have a $\Sigma_r$ action on $H^\ast(\CM_{\alpha_\bullet}^\qpar)$. This motivates the definition of a $\Sigma_r$ action on the quasi-parabolic descendent algebra $\BD^\qpar_{\alpha_\bullet}$, which will make the realization map $\Sigma_r$ equivariant. Recall that $\BD^\qpar_{\alpha_\bullet}$ is generated by symbols $\ch_k(\gamma)$ for $\gamma\in H^\ast(C)$ and $\ch_k(e_j)$ for $j=1, \ldots, r$. Define $d_j=e_{j}-e_{j-1}$ for $j=1, \ldots, r$. Then $\BD^\qpar_{\alpha_\bullet}$ is the free supercommutative algebra generated by $\ch_k(\gamma), \ch_k(d_j)$ modulo the relations
\[\ch_k(\pt)=\ch_k(d_1+\ldots+d_r)\,,\ch_1(\1)=d\,,\ch_0(\gamma^{<2})=0\,, \ch_0(d_j)=1\,.\]
We define the action of $\sigma\in \Sigma_r$ on $\BD^\qpar_{\alpha_\bullet}$ by $\ch_k(\gamma)\mapsto \ch_k(\gamma)$ for $\gamma\in H^\ast(C)$ and $\ch_k(d_i)\mapsto \ch_k(d_{\sigma(i)})$. The action is well defined since it preserves the relations above. Moreover, it is easy to see that the action commutes with the operator $\bR_{-1}$, and hence it preserves the subalgebra $\BD_{\alpha_\bullet, \inv}^{\qpar}$.

\begin{corollary}\label{cor: weylantisymmetry}
Suppose that $\gcd(r,d)=1$ and let $c\in S_{r,d}$ be sufficiently close to $(0,\ldots, 0)\in \overline S_{r,d}$, so that Proposition \ref{prop: flagbundle} applies. Then integration over $M_{r,d}^\full(c)$ is $\Sigma_r$ anti-invariant, i.e.
\begin{equation}\label{eq: weylantisymmetry}
\int_{M_{r,d}^\full(c)}\sigma\cdot D=(-1)^{\textup{sgn}(\sigma)}\int_{M_{r,d}^\full(c)}D\,
\end{equation}
for every $D\in \BD_{\alpha_\bullet, \inv}^{\qpar}$ and $\sigma\in \Sigma_r$. 
\end{corollary}
\begin{proof}
Let $\pi\colon M_{r,d}^\full(c)\to M_{r,d}$ be the full parabolic bundle given by Proposition \ref{prop: flagbundle}. Fix a universal vector bundle $\CV$ on $M_{r,d}$ and consider the induced universal parabolic bundle $\CW=(\pi^\ast\CV, \CF_\bullet)$ on $M_{r,d}^\full(c)$. The result is immediate from Proposition \ref{prop: weylflagbundle} after we observe that 
\[\xi_{\CV}(\ch_k(d_i))=\ch_k(\CF_i/\CF_{i-1})=\frac{1}{k!}c_1(\CF_i/\CF_{i-1})^k\,.\qedhere\]
\end{proof}

\subsection{Affine Weyl symmetry}\label{subsec: affineweyl} As sketched in the introduction (Section \ref{subsec: ingredientsintro}), conjugating the Weyl symmetry with Hecke isomorphisms produces an extremely strong constraint on the integrals over moduli spaces of full parabolic bundles.

Let $r>2$ and let $\alpha_\bullet, \alpha_\bullet'$ be the full topological types $(r,1,f_\bullet)$ and $(r,-1, f_\bullet)$, respectively. As explained in Section \ref{subsec: heckestability}, given $c'\in S_{r,-1}$ there is $c''\in S_{r,1}$ that makes the moduli spaces $M_{r,-1}^\full(c')$ and $M_{r,1}^\full(c'')$ isomorphic via a Hecke modification. In such a case, by Proposition \ref{lem: hecketr}, we have a commutative diagram 
\begin{center}
\begin{tikzcd}
\BD^{\qpar}_{\alpha_\bullet'} \arrow[r]\arrow[d, "H"',  "\sim" {sloped,inner sep=.4mm}] & H^\ast\big(M_{r,-1}^\full(c'')\big)\arrow[d,  "\sim" {sloped,inner sep=.4mm}]\\
\BD^{\qpar}_{\alpha_\bullet}\arrow[r]& H^\ast\big(M_{r,1}^\full(c')\big)
\end{tikzcd}
\end{center}
where $H$ is the Hecke homomorphism of the quasi-parabolic descendent algebra. Explicitly, $H$ is defined on generators by
\begin{align*}
H(\ch_k(d_j))&=\ch_k(d_{\tau(j)})\\
H(\ch_{k+1}(\1))&=\ch_{k+1}(\1)-\ch_k(d_{r-1})-\ch_k(d_r)\\
H(\ch_{k}(\gamma))&=\ch_{k}(\gamma)\quad\textup{ for } \gamma\in H^{>0}(C)\,.
\end{align*}
where $\tau$ is the permutation $\tau(i)=i+r-2 \mod r$. If $c'$ is chosen so that $M_{r,-1}^\full(c')$ is a full flag bundle over $M_{r,-1}$, then by Corollary \ref{cor: weylantisymmetry}, integration over $M_{r,1}^\full(c'')$ is anti-symmetric with respect to the action of $\widetilde \Sigma_r \coloneqq H\circ \Sigma_r\circ H^{-1}$. Note that $\widetilde \Sigma_r$ is isomorphic again to the symmetric group, but we use the tilde to distinguish the different action on $\BD^\qpar_{\alpha_\bullet}$. The following elementary proposition is at the heart of the proofs of Theorems \ref{thm: reconstruction}, \ref{thm: B} and \ref{thm: C}.

\begin{proposition}\label{prop: affineweylvanishing}
Let $\alpha_\bullet$ be a full topological type with rank $r>2$. Let $f\colon \BD^{\qpar}_{\alpha_\bullet}\to \BQ$ be a functional with the property that 
\begin{equation}\label{eq: sillyfunctionalassumption}
f(\ch_k(d_i)\cdot D)=\frac{1}{k!}f(\ch_1(d_i)^k\cdot D)\,.
\end{equation}
If $f$ is anti-invariant with respect to the actions of both $\Sigma_r$ and $\widetilde \Sigma_r$, then $f=0$. 
\end{proposition}

Note that, given a class $u\in H_\ast(\CM^\qpar_{\alpha_\bullet})$, the functional
\[\BD_{\alpha_\bullet}^{\qpar}\longrightarrow H^\ast(\CM_{\alpha_\bullet}^\qpar)\xlongrightarrow{\int_u} \BQ\]
satisfies \eqref{eq: sillyfunctionalassumption} since $\ch_k(d_i)$ and $\ch_1(d_i)^k/k!$ both map to
\[\frac{1}{k!}c_1(\CF_{i}/\CF_{i-1})^k\in H^\ast\big(\CM_{\alpha_\bullet}^\qpar\big)\,,\]
via the realization morphism.

\begin{proof}
A straightforward calculation shows that $\tau \sigma^{-1}\tau^{-1}H\sigma H^{-1}$ sends
\[\ch_{k+1}(\1)\mapsto\ch_{k+1}(\1)+\ch_k(d_r)+\ch_k(d_{r-1})-\ch_k(d_{\tau \sigma(1)})-\ch_k(d_{\tau\sigma(2)})\,\]
and preserves every other descendent, i.e. $\ch_k(g)\mapsto \ch_k(g)$ for $g=d_j$ or $g\in H^{>0}(C)$. Let $S_{i,j}$ be the automorphism of $\BD_{\alpha_\bullet}^{\qpar}$ which sends
\[\ch_{k+1}(\1)\mapsto\ch_{k+1}(\1)+\ch_k(d_i)-\ch_k(d_j)\,\]
and preserves every other descendent. If $j\neq r-1$, we can pick $\sigma$ so that $\tau\sigma(1)=j$, $\tau \sigma(2)=r-1$ and then the previous calculation shows that $\tau \sigma^{-1}\tau^{-1}H\sigma H^{-1}=S_{r,j}$. Hence, $S_{r,j}$ is in the subgroup of automorphisms generated by $\Sigma_r$ and $\widetilde \Sigma_r$, and the functional $f$ is invariant with respect to $S_{r,j}$. Similarly, we conclude the same for $S_{r-1, j}$ for any $j\neq r$. Using that
\[S_{i,j}^{-1}=S_{j,i}\textup{ and }S_{i,j}\circ S_{j,k}=S_{i,k}\,,\]
and the assumption that $r>2$, we can conclude that all the $S_{i,j}$ are in the subgroup generated by $\Sigma_r, \widetilde \Sigma_r$. 

Let $t_i=\ch_1(d_i)\in \BD^\qpar_{\alpha_\bullet}$. For the proof, we introduce the equivalence relation $\sim$ on $\BD_{\alpha_\bullet}^{\qpar}$ which is defined as the smallest $\BQ$-linear equivalence relation such that 
\begin{equation}\label{eq: relationssim}\ch_k(d_i)\cdot D\sim \frac{1}{k!}t_i^k \cdot D\,,\quad D\sim -\sigma \cdot D\,,\quad D\sim S_{i,j}(D)
\end{equation}
for any $D\in \BD^\qpar_{\alpha_\bullet}$, $\sigma\in \Sigma_r$, $1\leq i, j\leq r$. By $\BQ$-linear we mean that we impose that $D_1\sim D'_1$ and $D_2\sim D_2'$ imply that $\lambda_1D_1+\lambda_2D_2\sim \lambda_1D_1'+\lambda_2D_2'$. Since by assumption $f(D)=f(D')$ for $D\sim D'$, the proposition is equivalent to the statement that $D\sim 0$ for any $D$.

It is enough to prove that $D\sim 0$ if $D$ is in the subalgebra of $\BD^\qpar_{\alpha_\bullet}$ generated by $t_1, \ldots, t_r$ and $\ch_k(\gamma)$ for $k\geq 0, \gamma\in H^\ast(C)$. To show it, we introduce a filtration $0=F_{-1}\subseteq F_0\subseteq F_1\subseteq \ldots$ of this subalgebra and prove the claim for $D\in F_k$ by induction on $k$. The $k$-th step of the filtration $F_k$ is defined to be the span of monomials of the form
\[D=\prod_{i=1}^\ell \ch_{k_i+1}(\1)\cdot D'\]
where $D'$ is a monomial in $t_1, \ldots, t_r$ and $\ch_k(\gamma)$ for $\gamma\in H^{>0}(C)$, and $\sum_{i=1}^\ell k_i\leq k$. Note that $\ch_0(\1)=0$, so $F_k=0$ for $k<0$. In particular, the statement is tautological for $m<0$. We now induct on $k$, so assume that $D\sim 0$ for every $D\in F_{k-1}$. We begin by showing that $t_i \cdot D\sim t_j\cdot D$ for every $i,j$ and $D\in F_k$. We may write $D$ in the form 
\[D=\ch_2(\1)^{m}\cdot \prod_{i=1}^{\ell'} \ch_{k_i+1}(\1)\cdot D'\in F_k\]
where $m\geq 0$, all the $k_i$ are $>1$ and $\ell'=\ell-m$. Then,
\begin{align*}\ch_2(\1)\cdot D&\sim S_{i,j}(\ch_2(\1)\cdot D)\\
&\sim (\ch_2(\1)+t_i-t_j)^{m+1}\cdot \prod_{i=1}^{\ell'} \left(\ch_{k_i+1}(\1)+\frac{1}{k!}t_i^{k_i}-\frac{1}{k!}t_j^{k_i}\right)\cdot D' \\
&=\ch_2(\1)\cdot D+(m+1)(t_i-t_j)\cdot D+\ldots
\end{align*}
where all the remaining terms in $\ldots$ are in $F_{k-1}$. By induction, they are equivalent to $0$ so we conclude that $t_i\cdot D\sim t_j\cdot D$. We now write $D$ as
\[D=\prod_{j=1}^r t_j^{n_j}\cdot D''\]
for $D''\in F_k$ a monomial in $\ch_k(\gamma)$ for $\gamma\in H^\ast(C)$. Then $D\sim D'''\coloneqq t_1^{n_1+\ldots+n_r}\cdot D''$. If $\sigma$ is the two cycle $(2\,3)$ then 
\[D\sim D'''=\sigma(D''')\sim -D'''=-D\]
so we finally conclude that $D\sim 0$, which proves the induction step and the proposition.\qedhere
\end{proof}

\begin{remark}\label{rmk: affineweyl}
It can be seen from the proof of the previous proposition that, for $r>2$, the subgroup of $\Aut(\BD^\qpar_{\alpha_\bullet})$ generated by $\Sigma_r$ and $\widetilde \Sigma_r$ is isomorphic to the affine Weyl group
\[\Lambda\rtimes \Sigma_r\]
where $\Lambda=\{(a_1, \ldots, a_r)\colon a_1+\ldots+a_r=0\}\subseteq \BZ^r$ and $\Sigma_r$ acts on $\Lambda$ by permuting the coordinates. The element 
\[(0, \ldots, 1, 0,\ldots, 0, -1, \ldots, 0)\in \Lambda\]
corresponds to the automorphism $S_{i,j}$. 
\end{remark}

The next Corollary is essentially the same statement, but for the algebra of weight zero descendents.

\begin{corollary}\label{cor: affineweylvanishing}
Let $\alpha_\bullet$ be a full topological type with rank $r>2$. Let $u\in H_\ast(\CM^{\qpar, \rig}_{\alpha_\bullet})$. If the linear functional
\[\BD_{\alpha_\bullet, \inv}^{\qpar}\longrightarrow H^\ast(\CM^{\qpar, \rig}_{\alpha_\bullet})\xlongrightarrow{\int_u} \BQ\]
is anti-invariant with respect to the actions of both $\Sigma_r$ and $\widetilde \Sigma_r$, then $u=0$. 
\end{corollary}
\begin{proof}
Since $\bR_{-1}\ch_0(\pt)=r\neq 0$, the proof of \cite[Lemma 5.11]{blm} shows that there is a unique lift $\tilde u\in H_\ast(\CM_{\alpha_\bullet}^\qpar)$ of $u$ with the property that $\ch_1(\pt)\cap \tilde u=0$. Then the induced functional $\int_{\tilde u}$ on $\BD^\qpar_{\alpha_\bullet}$ is the unique extension of $\int_u$ on $\BD^\qpar_{\alpha_\bullet, \inv}$ (see also \cite[Proposition 1.13, Lemma 1.16]{klmp}) that vanishes on the ideal generated by $\ch_1(\pt)$. Clearly this extension is still anti-invariant, and the conclusion follows from Proposition \ref{prop: affineweylvanishing} applied to $\int_{\tilde u}$.
\end{proof}

\subsection{Proof of the reconstruction theorem} \label{subsec: proofreconstruction}
We now give a proof of Theorem~\ref{thm: reconstruction}. Suppose that we have two collections of classes
\[u_\alpha, v_\alpha\in H_\ast(\CM^{\paar, \rig}_\alpha)\subseteq \widecheck \V^\paar\]
that satisfy (0)-(3), and we want to show that $u_\alpha=v_\alpha$ for every $\alpha$. By the flag bundle formula, the classes $u_\alpha$ with $\alpha$ being a full parabolic type determine all the other ones. Denote $u_\alpha$ by $u^\full_{r,d}(c)$ if $\alpha=(r,d,f_\bullet, c)$ is full. We induct on the rank and assume that $u_{\alpha'}=v_{\alpha'}$ for $\rk(\alpha')<r$. Our base case is $r=1$, which holds by hypothesis, but going from $r=1$ to $r=2$ requires a modification of the argument (due to Corollary \ref{cor: affineweylvanishing} not holding for $r=2$). So we assume for now that $r>2$ and sketch below the necessary modifications for the $r=2$ case.

By the Hecke compatibility, it is enough to prove that $u_{r,1}^\full(c)=v_{r,1}^\full(c)$ for $d=1$. It follows from the wall-crossing formula and the induction hypothesis that, for any $c,c'\in \overline S_{r,1}$,
\[u_{r,1}^\full(c)-u_{r,1}^\full(c')=v_{r,1}^\full(c)-v_{r,1}^\full(c')\]
after forgetting weights. In other words, the difference $w=u_{r,1}^\full(c)-v_{r,1}^\full(c)$ does not depend on the stability condition $c\in \overline S_{r,1}$ after forgetting weights, i.e. $w$ does not depend on $c$ when regarded as an element of $H_\ast(\CM^{\qpar, \rig}_{\alpha_\bullet})$. By the Weyl symmetry and the Hecke compatibility assumptions, the functional $\int_w$ is anti-invariant under $\Sigma_r$ and $\widetilde \Sigma_r$, since $u_{r,1}^\full(c)$ and $v_{r,1}^\full(c)$ are anti-invariant for particular choices of $c$. Then Corollary \ref{cor: affineweylvanishing} implies that $w=0$.

\subsubsection{The rank 2 case}\label{subsubsec: rank2punctures}
We now sketch the necessary modifications to deal with the rank 2 case. Following \cite[Section 9]{ST}, we need to consider parabolic bundles with two punctures. The introduction of a second puncture allows extra flexibility to define Hecke modifications, since it becomes possible to modify the parabolic bundle at the two points, and this will recover the affine Weyl symmetry in rank 2.

Fix $p_1=p$ and a second distinct point $p_2\in C$. There are moduli spaces 
\[M_{2,d}^{2\textup{-pts}}(c)\]
that parametrize a vector bundle $V$ of rank 2 and degree $d$, together with two 1-dimensional subspaces
\[0\subsetneq F^1 \subsetneq V_{p}\textup{ and }0\subsetneq F^2\subsetneq V_{q}\,.\]
These moduli spaces depend on a choice of weight vector $c$ on some space of stability conditions $S_{2,d}^{2\textup{-pt}}$, which can be identified with the square
\[\{(\lambda_1, \lambda_2)\colon 0<\lambda_1, \lambda_2<0\}\,.
\]
As explained in \cite[Section 9.1]{ST}, there are two chambers in this space of stability conditions and a wall-crossing formula can be proven by using the exact same techniques as in Section \ref{sec: proofsimplewc}. In particular, the difference
\[\int_{M_{2,d}^{2\textup{-pts}}(c)}D-\int_{M_{2,d}^{2\textup{-pts}}(c')}D\]
can be written as an integral over a product of two Jacobians. The entire vertex algebra technology can easily be adapted to this setting with two punctures, and the wall-crossing formula is written as usual as an equality
\[[M_{2,d}^{2\textup{-pts}}(c)]- [M_{2,d}^{2\textup{-pts}}(c')]=\big[[\textup{Jac}_1],[\textup{Jac}_2] \big]\]
in the associated Lie algebra.

We can define the quasi-parabolic descendent algebra $\BD^{2\textup{-pts}}_{2,d}$ adapted to the 2 punctured situation. It is the supercommutative algebra generated by
\[\ch_k(\gamma)\,\textup{ for }\gamma\in H^\ast(C)\quad \textup{ and \quad }\ch_k(d_{i}^{\hspace{0.1em}j})\,,\textup{ for }i,j=1,2\,.\]
There is a realization map $\BD^{2\textup{-pts}}_{2,d}\to H^\ast(M_{2,d}^{2\textup{-pts}}(c))$ which is defined as usual for $\ch_k(\gamma)$ and sends 
\[\ch_k(d^{\hspace{0.1em}j}_1)\mapsto \ch_k(\BF^{\hspace{0.1em}j})\,,\,\ch_k(d^{\hspace{0.1em}j}_2)\mapsto \ch_k(\BV_{|p_j}/\CF^{\hspace{0.1em}j})\]
where by $\BV_{|p_j}$ we mean the restriction of $\BV$ to the fiber over $p_j\in C$. 

There is a natural Weyl group $G\coloneqq \Sigma_2\times \Sigma_2$ action on $\BD^{2\textup{-pts}}_{2,d}$ where $\sigma=(\sigma_1, \sigma_2)$ acts as
\[\sigma\cdot \ch_k(d_i^{\hspace{0.1em}j})=\ch_k(d_{\sigma_j(i)}^{\hspace{0.1em}j})\,.\]
If $d$ is odd and $c$ is in one of the two chambers of the space of stability condition, the moduli space $M_{2,d}^{2\textup{-pts}}(c)$ is a $\BP^1\times \BP^1$ bundle over $M_{2,d}$. More precisely,
\[M_{2,d}^{2\textup{-pts}}(c)=\BP(\CV_{|p_1})\times_{M_{2,d}} \BP(\CV_{| p_2})\,.\]
For such $c$, the action of $G$ descends to the cohomology and makes integrals anti-invariant in the sense that 
\[\int_{M_{2,d}^{2\textup{-pts}}(c)}\sigma\cdot D=\textup{sgn}(\sigma_1)\textup{sgn}(\sigma_2)\int_{M_{2,d}^{2\textup{-pts}}(c)} D\,.\]
This follows from applying Proposition \ref{prop: weylflagbundle} to the two descriptions of $M_{2,d}^{2\textup{-pts}}(c)$ as a $\BP^1$ bundle. 

To get our second copy of the Weyl group acting, we consider Hecke modifications. Since we now have two points, we can either modify a parabolic bundle at $p_1$ or at $p_2$. Composing a modification at $p_1$ with a modification at $p_2$ and choosing stability conditions $c', c'', c'''$ appropriately, we get a $\BP^1\times \BP^1$-bundle over $M_{2,-1}$
\[M_{2,1}^{2\textup{-pts}}(c')\simeq M_{2,0}^{2\textup{-pts}}(c'')\simeq M_{2,-1}^{2\textup{-pts}}(c''')\to M_{2,-1}\,.\]
Let $H\colon \BD^{2\textup{-pts}}_{2,-1}\to \BD^{2\textup{-pts}}_{2,1}$ be the corresponding Hecke operator on the descendent algebra; it sends 
\[\ch_{k+1}(\1)\mapsto \ch_{k+1}(\1)-\ch_k(d_2^1)-\ch_k(d_2^2)\,\textup{ and }\ch_k(d_i^{\hspace{0.1em}j})\mapsto \ch_k(d_{2-i}^{\hspace{0.1em}j})\,.\]
As in Section \ref{subsec: affineweyl}, the action of $\widetilde G\coloneqq H\circ G\circ H^{-1}$ descends to the cohomology of $M_{2,1}^{2\textup{-pts}}(c')$ for $c'$ chosen as before. 

Finally, it is easy to adapt the argument given in the proof of Proposition \ref{prop: affineweylvanishing} to show that any functional on $\BD^{2\textup{-pts}}_{2,1}$ (imposing the analog of \eqref{eq: sillyfunctionalassumption}) which is anti-invariant with respect to the actions of both $G$ and $\widetilde G$ must vanish. With this, the argument in the proof of the reconstruction theorem also shows that integrals over $M_{2,1}^{2\textup{-pts}}(c)$ are determined from integrals over the Jacobian. Since $M_{2,1}^{2\textup{-pts}}(c)$ is a $\BP^1$-bundle over the moduli of rank 2 degree 1 parabolic bundles punctured only at $p$, we recover also those integrals.

\section{Newstead vanishing and the Chern filtration}\label{sec: newstead}

Based on low genus calculations by himself and Ramanan, Newstead conjectured in \cite{newsteadrk2} that, in the moduli space $M_{2,1}$ or rank 2 vector bundles, we have the vanishing 
\[\beta^g\in H^{4g}(M_{2,1})\]
where $\beta$ is (up to normalization) the realization of the tautological class $\ch_2(\pt)$. This was later shown by Thaddeus \cite[Section 5]{thaddeuscft}. The result was generalized by Earl--Kirwan in \cite{EKvanishing} and independently by Jeffrey--Weitsman \cite{JW}. They showed that any polynomial in the descendents of the point class $\ch_k(\pt)$ vanishes in cohomological degree $>2r(r-1)(g-1)$. 

It was observed in \cite{LMP} that these sort of vanishing results can be understood in terms of the Chern filtration on the cohomology of $M_{r,d}$. Using the methods of \cite{EKvanishing}, it is shown in \cite[Theorem 0.2, Corollary 1.7]{LMP} that the top Chern degree of $M_{r,d}$ is $r(r+1)(g-1)+2$. This implies the higher rank Newstead conjecture, but it is actually a stronger statement. For example, it implies as well the vanishing of
\[\beta^a \gamma^b \in H^{4a+6b}(M_{2,1})\]
for $a+b\geq g$, also observed by Thaddeus, where $\gamma$ is quadratic in the descendents $\ch_2(\textup{class in }H^1(C))$.

In this section, we formulate and prove an analogous result for parabolic bundles. We point out that our proof is independent and fundamentally different from the previous ones, so in particular we reprove the aforementioned results in the case of vector bundles.

\subsection{The Chern filtration for parabolic bundles}
\label{subsec: chern}

We start by defining the Chern grading on the parabolic descendent algebra, and induced Chern filtration.  

\begin{definition}\label{def: chernfiltration}
The Chern grading on the descendent algebra $\BD_{\alpha_\bullet}^\qpar$ is the multiplicative grading for which the generators
\[\ch_k(g),\,g\in H^\ast(C)\cup \{e_1, \ldots, e_l\}\]
are homogeneous of degree $k$. We write $\deg^C(D)=k$ if $D$ is homogeneous of Chern degree~$k$. Denote by $C_\bullet \BD_{\alpha_\bullet}^\qpar$ the induced filtration, i.e.
\[C_k \BD_{\alpha_\bullet}^\qpar=\textup{span}_\BQ\big\{\ch_{k_1}(g_1)\ldots \ch_{k_m}(g_m)\colon k_1+\ldots+k_m\leq k\big\}\,.\] 
We write $\deg^C(D)\leq k$ if $D\in C_k \BD_{\alpha_\bullet}^\qpar$. Denote by \[C_\bullet \BD_{\alpha_\bullet, \inv}^\qpar=\BD_{\alpha_\bullet, \inv}^\qpar\cap C_\bullet\BD_{\alpha_\bullet}^\qpar\] 
the induced subspace filtration. The same definitions can be made in the parabolic descendent algebra $\BD^\paar_\alpha$ or $\BD^\paar_{\alpha, \inv}$.
\end{definition}

The Chern filtration on the descendent algebra induces a Chern filtration on the cohomology of a moduli space $M_\alpha$, with $\alpha$ regular, defined by
\[C_\bullet H^\ast(M_\alpha)=\textup{im}\big(\xi\colon C_\bullet \BD_{\alpha_\bullet, \inv}^\qpar \to H^\ast(M_\alpha )\big)\,.\]
Alternatively, the Chern filtration can be defined in a more concrete way that avoids the use of weight 0 descendents, and instead uses a point normalization:
\[C_\bullet H^\ast(M_\alpha)=\textup{im}\big(\xi_\pt\colon C_\bullet \BD_{\alpha_\bullet}^\qpar \to H^\ast(M_\alpha )\big)\,\]
where $\xi_\pt$ is the unique extension of $\xi$ to $\BD_{\alpha_\bullet}^\qpar $ that annihilates $\ch_1(\pt)$; see \cite[Propositions 1.13, 6.5]{klmp}. This is closer to the way that tautological classes are treated in the literature, for example in \cite{JK}. 

The main result in this section is the following:

\begin{theorem}\label{thm: newstead}
Let $\alpha=(r, d, f_\bullet, c)$ and $[M_\alpha]$ be the corresponding class in the Lie algebra $\widecheck{\V}^\paar\subseteq \widecheck{\V}^\paar_{\tr}$, and let $D\in \BD^\qpar_{\alpha_\bullet, \inv}$. If $\deg^C(D)\leq \dim M_\alpha+r(g-1)$ then 
\[\int_{[M_\alpha]}D =0\,.\]
If $\deg^C(D)=\dim M+r(g-1)+1$, then
\[(-1)^{(r-1)d}\int_{[M_\alpha]} D\]
does not depend on $d$ and on $c$.%\footnote{We clarify the statement that the integral does not depend on $d$: for $d'\neq d$, we have different descendent algebras $\BD^\qpar_{\alpha_\bullet}$ and $\BD^\qpar_{\alpha_\bullet'}$. Recall that we have the relation $\ch_1(\1)=d$  (resp. $=d'$) in $\BD^\qpar_{\alpha_\bullet}$ (resp. in $\BD^\qpar_{\alpha_\bullet'}$). We identify the two algebras by identifying $\ch_k(g)$ with $\ch_k(g)$ for any generator other than $\ch_1(\1)$. Then the claim is that $\int_{[M_\alpha]} \xi_\pt(D)=\int_{[M_\alpha']} \xi_\pt(D')$ where $D'\in \BD^\qpar_{\alpha_\bullet'}$ is the corresponding descendent.}
\end{theorem}

We will give a proof in Section \ref{subsec: proofnewstead}. Note that when $\alpha$ is regular the first statement is equivalent to 
\[C_{\dim M_\alpha+r(g-1)}H^{\textup{top}}(M_\alpha)=0\,.\]

Recall that $\ch_k(e_j)$ and $\ch_k(\gamma)$ have cohomological degrees 
\[\deg(\ch_k(e_j))=2k \textup{ and }\deg(\ch_k(\gamma))=2k+2-\deg(\gamma)\,.\]
Thus, the Chern degree is always between half of the cohomological degree and the cohomological degree. For a polynomial in classes $\ch_k(\pt), \ch_k(e_j)$, the Chern degree is half of the cohomological degree. In particular, \eqref{eq: vanishingchernflag} is equivalent to the statement that 
\[C_{k}H^{2k}(M_\alpha)=0 \quad \textup{ for }k\geq \dim M-r(g-1)\,.\]
More generally, we have the following:

\begin{corollary}\label{cor: newstead}
Let $\alpha\in C(I)$ be regular. If
\[m-k\geq \dim M_\alpha-r(g-1)\]
then $C_k H^m(M_\alpha)=0$. In particular, any polynomial in the classes 
\[\xi_\pt(\ch_i(\pt)), \xi_\pt(\ch_i(e_j))\in H^{2i}(M_\alpha)\] 
of cohomological degree $\geq 2(\dim M_\alpha-r(g-1))$ vanishes.
\end{corollary}
\begin{proof}
Let $\delta=\dim M_\alpha$. Suppose $D\in C_k H^m(M_\alpha)$ is non-zero. By Poincaré duality, there is $E\in H^{2\delta-m}(M_\alpha)$ such 
\[\int_{M_\alpha} D\cdot E\neq 0\,.\]
Recall that the Chern degree is bounded by the cohomological degree, i.e. \[H^{2\delta-m}(M_\alpha)=C_{2\delta-m}H^{2\delta-m}(M_\alpha)\,.\]
But then
\[D\cdot E\in C_{2\delta+k-m }H^{2\delta}(M_\alpha)\,,\]
so by Theorem \ref{thm: newstead} it follows that
\[2\delta+k-m>\delta +r(g-1)\]
and we are done. For the second part of the statement, a polynomial in such classes of cohomological degree $2k$ is in $C_k H^{2k}(M_\alpha)$, so the statement follows immediately.\qedhere
\end{proof}

\begin{remark}\label{rmk: optimal}
The bounds in Theorem \ref{thm: newstead} and Corollary \ref{cor: newstead} are optimal, i.e. the \textit{top Chern degree} of $M_\alpha$ is $\dim M_\alpha-r(g-1)+1$. In the case of moduli of bundles $M_{r,d}$ and $\gcd(r,d)=1$, it is shown in \cite{EKvanishing} that there are classes $\eta, \theta\in \BD^\paar$ such that: 
\begin{enumerate}
\item $\eta$ is a polynomial in descendents of the form $\ch_k(\pt)$, and has
\[\frac{1}{2}\deg(\eta)=\deg^C(\eta)=r(r-1)(g-1)=\dim M_{r,d}-r(g-1)-1\,,\]
\item $\theta$ is given by\footnote{Note that Earl--Kirwan work with moduli of bundles with fixed determinant, so they consider only the class $\ch_2(\1)^{(r-1)(g-1)}$. The realization of $\prod_{i=1}^{2g}\ch_1(\gamma_i)$ is a multiple of the Poincaré dual to the loci of $M_{r,d}$ with a fixed determinant.}
 \[\theta=\ch_2(\1)^{(r-1)(g-1)}\prod_{i=1}^{2g}\ch_1(\gamma_i)\]
where $\{\gamma_i\}$ is a basis of $H^1(C)$. In particular, \[\deg(\theta)=\deg^C(\theta)=2r(g-1)+2\]
\item We have
\[\int_{M_{r,d}}\xi_{\pt}(\eta\cdot \theta)\neq 0\,.\]
\end{enumerate}
If $c\in S_{r,d,f_\bullet}$ is chosen so that there is a flag bundle morphism $\pi\colon M_{r,d, f_\bullet}(c)\to M_{r,d}$ and Proposition \ref{prop: flagbundleclass} holds, then 
\[\int_{M_{r,d, f_\bullet}(c)}D\neq 0\textup{ for }D=\xi_\pt(\eta\cdot \theta)\cdot c_{\textup{top}}(\Delta^\ast \Xi)\,.\]
Recall the definition of $\Xi$ in Section \ref{subsec: mapsVA} and that $\Delta^\ast \Xi$ is the relative tangent bundle of $\pi$. Clearly, $c_{\textup{top}}(\Delta^\ast \Xi)$ is a polynomial in descendents $\xi_\pt(\ch_k(e_i))$. Hence, $D$ has Chern degree equal to $\dim M_\alpha+r(g-1)+1$. By the second part of Theorem~\ref{thm: newstead}, $\int_{[M_{r,d, f_\bullet}(c)]}D\neq 0$ actually holds for any $d$ and $c\in S_{r,d,f_\bullet}$. The class $\xi_\pt(\eta)\cdot c_{\textup{top}}(\Delta^\ast \Xi)$ then provides an example of a non-zero polynomial in $\xi_\pt(\ch_k(\pt)), \xi_\pt(\ch_k(e_i))$ with cohomological degree equal to $2(\dim M_\alpha-r(g-1)-1)$. 
\end{remark}

\begin{remark}
The part of Theorem \ref{thm: newstead} regarding independence of $d$ is analogous to $\chi$-independence phenomena conjectured for moduli of 1-dimensional sheaves on del Pezzos surfaces, see \cite[Remark 1.11]{LMP}. Note that $(-1)^{(r-1)d}=(-1)^{r-1}$ unless $r,d$ are both even, so in the case of $M_{r,d}$ the sign only plays a role when we allow $\gcd(r,d)\neq 1$, i.e. moduli with strictly semistable objects. This is why our result is not incompatible with \cite[Proposition 1.8]{LMP}, where there is no sign. Note that the $(-1)^{(r-1)d}$ sign appears in some of the formulas of \cite{bu}, for example in Theorem~5.15 ($r=2$) or in Theorem 6.6 ($g=1$).
\end{remark}

\begin{remark}
In the same paper \cite{newsteadrk2}, Newstead also makes a second conjecture -- in rank 2, but later generalized to higher $r$ -- regarding the vanishing of Chern classes $c_k(TM_{r,d})$ of the tangent bundle for $k>r(r-1)(g-1)$. This conjecture was proved in \cite{TW}. Although our methods seem almost fit to reprove it, we could not quite figure out how to do it due to the subtle but important point that $c_k(TM_{r,d})$ does not behave well with respect to the Weyl action. 
\end{remark}

\subsection{Proof of Theorem \ref{thm: newstead}}\label{subsec: proofnewstead}
We start by reducing to the case of full parabolic bundles. Let $c\in S_{r,d, f_\bullet}$. Then,
\[\int_{[M_{r,d, f_\bullet}(c)]}D=\textup{const}\cdot \int_{[M_{r,d}^\full(\overline c)]}c_N(\Delta^\ast \Xi) \pi^\ast D\]
where $N=\dim M_{r,d}^\full(\overline c)-\dim M_{r,d, f_\bullet}(c)$ and $\pi$ is the flag bundle morphism. By the definition of $\Xi$, it is clear that $c_N(\Delta^\ast \Xi)$ is a polynomial in classes $\xi_\pt(\ch_k(e_j))$, hence $c_N(\Delta^\ast \Xi) $ has Chern degree equal to $N$. It is also clear that $\pi^\ast$ preserves the Chern filtration, and that shows that the statement for $M_{r,d}^\full(\overline c)$ implies the statement for $M_{r,d, f_\bullet}(c)$.

We now focus on full parabolic bundles and induct on the rank $r=\rk(\alpha)$. When $r=1$, $M_\alpha$ is a Jacobian and our claim is straightforward. Indeed, the descendents $\xi_\pt(\ch_k(\gamma))$ all vanish except when $k=1$ and $\gamma\in H^1(C)$ or $k=2$, $\gamma\in H^0(C)$; in both cases, the Chern degree is equal to the cohomological degree, so $C_{2g-1}H^{2g}(M_\alpha)=0$. There is nothing to prove regarding wall-crossing invariance and the independence of $d$ is easy since all the moduli spaces $M_{1,d}$ are isomorphic.

Suppose now that the claim holds for $r'<r$ and assume that $r>2$ (once again, the case $r=2$ requires the use of 2 punctures, cf. Section \ref{subsubsec: rank2punctures}). We use the induction hypothesis to prove the wall-crossing independence first, and then make use of the affine Weyl symmetry. Suppose that $c^\pm$ are stability conditions in $\overline S_{r,d}$. To ease the notation, we assume they are both regular and in adjacent chambers so that we are in the setting of simple wall-crossing, as described in Section \ref{subsec: simplewc}.

The simple wall-crossing formula (cf. Theorem \ref{thm: simplewcI}) expresses the difference between integrals over $M_{r,d}^\full(c_\pm)$ as
\[\int_{M_{r,d}^\full(c_+)}D-\int_{M_{r,d}^\full(c_-)} D=\sum_{i}\int_{M_1}D_i^1\cdot \int_{M_2}D_i^2\]
where $\sum_{i}D_i^1\otimes D_i^2=\Res_{z=0} Y^\vee(D)$. The following Lemma controls the Chern degrees appearing in the wall-crossing term, and is arguably the heart of the proof.

\begin{lemma}\label{lem: cdegreewc}
The dual state-to-field correspondence
\[Y^\vee\colon \BD_{\alpha+\beta}^\paar\to (\BD_{\alpha}^\paar\otimes \BD_\beta^\paar)\llbracket z^{-1}, z\rrbracket\]
has degree $\chi^\sym(\alpha, \beta)$ with respect to the Chern degree, where on the right hand side we consider the natural Chern degree on the tensor and regard $z$ as having Chern degree 1. In particular, the dual Lie bracket
\[\Res_{z=0} Y^\vee\colon \BD_{\alpha+\beta}^\paar\to \BD_{\alpha}^\paar\otimes \BD_\beta^\paar\]
has Chern degree $\chi^\sym(\alpha, \beta)+1$.
\end{lemma}
\begin{proof}
It is clear that the operator $\Sigma^\ast$ has degree 0 and that $\bR_{-1}$ has degree $-1$, therefore it is enough to prove that $C_{z^{-1}}$ has degree 0 or, equivalently, $c_k(\Theta)$ has Chern degree $k$. This is equivalent to $c_k(\Theta)$ being a polynomial in the (realization of) classes 
\begin{equation}\label{eq: classesDtensorD}
\ch_j(\pt)\otimes 1,\, 1\boxtimes \ch_j(\pt),\, \ch_j(e(t))\otimes 1,\, 1\otimes \ch_j(e(t))\,,
\end{equation}
which we now prove; interestingly, this is true only for the symmetrization $\Theta=\Ext_{12}^\vee+\Ext_{21}$, but not for each of the Ext complexes. 

 Recall our formula for the Ext complex \eqref{eq: extcomplex}. By Grothendieck--Verdier duality we have 
\[(Rp_\ast \CR\CH om(\CV_1, \CV_2)\big)^\vee=-Rp_\ast \CR\CH om(\CV_2, \CV_1\otimes \omega_C)\]
in $K$-theory. Therefore,
\begin{align*}\Theta=&Rp_\ast \CR\CH om\big(\CV_2, \CV_1\otimes (\CO_C-\omega_C)\big)\\
&-\sum_{t\in I}\big(\CH om(\CF_1(t), \partial\CF_2(t))^\vee+\CH om(\CF_2(t), \partial\CF_1(t))\big)\,.
\end{align*}
It is clear that the Chern character of the terms in the second line are polynomials in \eqref{eq: classesDtensorD}. The Chern character of the first line on the right hand side can be calculated via Grothendieck--Riemann--Roch as
\[p_\ast\big(\ch(\CV_2^\vee)\ch(\CV_1)\ch(\CO_C-\omega_C)\td(C)\big)=(2g-2)\sum_{a,b\geq 0}(-1)^b \ch_a(\pt)\otimes \ch_b(\pt)\]
since $\ch(\CO_C-\omega_C)=(2g-2)\pt$. It follows that $\ch_k(\Theta)$ is a polynomial in the classes  \eqref{eq: classesDtensorD}, and therefore the same is true for $c_k(\Theta)$ by Newton's identities. \qedhere
\end{proof}

\begin{remark}
Note that $\chi^\sym(\alpha, \beta)+1$ is equal to
\[1-\chi(\alpha, \alpha)+1-\chi(\beta, \beta)-(1-\chi(\alpha+\beta, \alpha+\beta))=\dim M_\alpha+\dim M_\beta-\dim M_{\alpha+\beta}\,.\]
\end{remark}

By the lemma and the previous remark, we have
\begin{align*}\deg^C(D_i^1)+\deg^C(D_i^2)&\leq \deg^C(D)-(\dim M_\pm-\dim M_1-\dim M_2)\\
&\leq
\dim M_1+r_1(g-1)+\dim M_2+r_2(g-1)+1\,,
\end{align*}
so we must have $\deg^C(D_i^j) \leq \dim M_j+r_j(g-1)$ for $j=1$ or for $j=2$. By the induction hypothesis,  it follows that $\int_{M_{r,d}^\full(c_+)}D=\int_{M_{r,d}^\full(c_-)} D$. Since any two regular stability conditions can be connected by a sequence of stability conditions which are related by simple wall-crossing, it follows that $\int_{M_{r,d}^\full(c)}D$ does not depend on $c$ for $c\in S_{r,d}$ regular. It is easy to extend this argument and take into account the iterated brackets that appear in the general wall-crossing formula \eqref{eq: fullwcformula}, showing that the same is true for non-regular weights. 

Both the Hecke operators and the Weyl action preserve the Chern filtration. Therefore, by the wall-crossing independence and the same argument used in Section~\ref{subsec: proofreconstruction}, the functional  $\int_{[M_{r,1}^\full(c)]}$ is anti-invariant with respect to the actions of $\Sigma_r$ and $\widetilde \Sigma_r$ when restricted to $C_{\dim M_\alpha+r(g-1)+1}\BD^{\qpar}_{\alpha_\bullet, \inv}$. The following refinement of Proposition \ref{prop: affineweylvanishing} now gives the vanishing in Chern degree $\dim M_\alpha+r(g-1)$:

\begin{corollary}\label{cor: vanishingweylsymmetrychern}
Let $\alpha_\bullet$ be a full topological type with rank $r>2$. Let $f\colon \BD^{\qpar}_{\alpha_\bullet}\to \BQ$ be a functional satisfying \eqref{eq: sillyfunctionalassumption}. If the restriction of $f$ to $C_{k+1}  \BD^{\qpar}_{\alpha_\bullet}$ is anti-invariant with respect to the actions of both $\Sigma_r$ and $\widetilde \Sigma_r$, then the restriction of $f$ to $C_k \BD^{\qpar}_{\alpha_\bullet}$ vanishes. The same is true if we replace $\BD^{\qpar}_{\alpha_\bullet}$ by $\BD^{\qpar}_{\alpha_\bullet, \inv}$.
\end{corollary}
\begin{proof}
The proof is exactly the same as Proposition \ref{prop: affineweylvanishing} (and Corollary \ref{cor: affineweylvanishing}). Let us just point out that when we show $t_i\cdot D\sim t_j \cdot D$ we use the anti-invariance applied to $\ch_2(\1)\cdot D$, which has Chern degree 1 larger than $t_1\cdot D, t_2\cdot D$, and this is why we need anti-invariance on $C_{k+1}$ to ensure vanishing on $C_k$.\qedhere
\end{proof}

Finally, let us address the invariance with respect to $d$. By the same argument as before, it is enough to consider full regular parabolic bundles. Let $\alpha=(r, d, f_\bullet, c)$ and $\alpha'=(r,1, f_\bullet, c')$ be two full parabolic types related by the Hecke isomorphism, i.e. such that $\alpha=h_\tau(\alpha')$. Recall the definition\footnote{Before, we defined the parabolic version of $h_\tau^\dagger$ from $\BD^\paar_{\alpha}\to \BD^\paar_{\alpha'}$, but we are now using a quasi-parabolic version $\BD^\qpar_{\alpha_\bullet}\to \BD^\qpar_{\alpha'_\bullet}$.}  of $h^\dagger\colon \BD^\qpar_{\alpha_\bullet}\to \BD^\qpar_{\alpha'_\bullet}$. The fundamental observation is that if $\deg^C(D)=\dim M_\alpha+r(g-1)+1$ then
\[h_\tau^\dagger(D)=\tau\cdot D+\ldots\]
where the extra terms have $\deg^C(\ldots)\leq \dim M_\alpha+r(g-1)$ and $\tau\in \Sigma_r$ is the permutation $\tau(i)=i+d-1 \mod r$. Thus, we have
\begin{align*}
\int_{M_\alpha} D=\int_{M_{\alpha'}}h_\tau^\dagger(D)=\int_{M_{\alpha'}}\tau\cdot D=\textup{sgn}(\tau)\int_{M_{\alpha'}} D 
\end{align*}
where we have used Proposition \ref{lem: heckestability} (and Lemma \ref{lem: hecketr}) in the first equality and the vanishing already proved in the second. In the third equality we have used the Weyl anti-symmetry -- Weyl anti-symmetry holds when $c'$ is small, but we have already shown wall-crossing independence when $\deg^C(D)=\dim M_\alpha+r(g-1)+1$, so it actually holds for any $c'$. 

The cycle decomposition of $\tau$ consists of $m$ cycles of length $r/m$ each, where $m=\gcd(r, d-1)$. Therefore, 
\[\textup{sgn}(\tau)=(-1)^{m(r/m-1)}=(-1)^{(r-1)(d-1)}\,.\]
Thus, the equality of integrals above can be rewritten as
\[(-1)^{(r-1)d}\int_{M_\alpha} D=(-1)^{(r-1)\cdot 1}\int_{M_\alpha'}D\,,\]
which shows the independence on $d$.\hfill $\Box$

\section{Virasoro constraints}\label{sec: virasoro}

The Virasoro constraints for moduli spaces of parabolic bundles were conjectured in \cite[Conjecture 1.4.9]{thesis}. In this section we will state them and give a proof. 

To write down the Virasoro operators it is convenient to introduce the Hodge shifted descendents $\ch_k^H(g)$, which are defined as
\[\ch_k^H(\gamma)=\ch_{k+1-p}(\gamma)\in \BD^\paar\]
where $\gamma\in H^{p,q}(C)$ and 
\[\ch_k^H(e(t))=\ch_k(e(t))\]
for $t\in I$. For notational convenience, we assume in this section that $I$ is finite. Given $t\in I\setminus \{1\}$, we write
\[\partial e(t)=e(t')-e(t)\in H_I(C)\]
where 
\[t'=\min \{t'>t\colon t'\in I\}\,.\]
When $t=1$ we simply set $\partial e(1)=0$. Note that the realization $\xi(\ch_k(\partial e(t))$ is, by definition, $\ch_k(\partial \CF(t))$. 

\begin{definition}
We define operators $\bL_n\colon \BD^\paar\to \BD^\paar$ for $n\geq -1$ as a sum $\bL_n=\bR_n+\bT_n$, where
\begin{enumerate}
\item $\bR_n$ is a derivation defined on generators $\chh_k(g)$, for $k\geq 0$ and $g\in H_I(C)$, by
\[\bR_n(\ch_k^H(g))=\left(\prod_{j=0}^n (k+j)\right)\ch_{n+k}^H(g)\,.\]
\item The operator $\bT_n$ is multiplication by
\begin{align*}\bT_n=\sum_{a+b=n}a!b!\left((1-g)\ch_a(\pt)\ch_b(\pt)+\sum_{t\in I}\ch_a(e(t))\ch_b(\partial e(t))\right)\,,
\end{align*}
where the sum runs over $a,b$ non-negative integers. In particular, $\bT_{-1}=0$ and $\bL_{-1}=\bR_{-1}$. 
\end{enumerate}
\end{definition}

The operators $\{\bL_n\}_{n\geq -1}$ satisfy the Virasoro bracket relation
\[[\bL_n, \bL_m]=(m-n)\bL_{n+m}\,.\]
To compare with the definition in \cite[Section 2.3]{blm} note that
\[(1-g)\ch_a(\pt)\ch_b(\pt)=(-1)^{p^L-1}\ch_a^H\ch_b^H(\td(C))\,\]
where the right hand side is defined in loc. cit.

The operators, which are a priori defined in $\BD^\paar$, induce well-defined operators $\bL_n\colon \BD^\paar_\alpha\to \BD^\paar_\alpha$, which we denote in the same way. Taking duals and summing over all $\alpha$  produces operators on the vertex algebra $\V^\paar_{\tr}$:
\begin{equation}\label{eq: operatorsva}L_n\colon \V^\paar_{\tr}=\bigoplus_{\alpha\in K(I)} (\BD^\paar_\alpha)^\vee\xrightarrow{\bL_n^\vee}\bigoplus_{\alpha\in K(I)} (\BD^\paar_\alpha)^\vee=\V^\paar_{\tr}\,.
\end{equation}

To write down the Virasoro constraints we need one further operator, namely $\bL_\inv\colon \BD^\paar\to \BD^\paar_\inv$, defined by
\[\bL_\inv=\sum_{n\geq -1}\frac{(-1)^n}{(n+1)!}\bL_{n}\circ \bR_{-1}^{n+1}\,.\]
The operator mapping to $\BD^\paar_\inv$ is equivalent to the statement that $\bR_{-1}\circ \bL_\inv=0$, which follows from the Virasoro commutator relations, see \cite[Section 2.5]{blm}.

\begin{remark}\label{rmk: quasiparoperators}
We point out that one can also define analogues of these operators in the quasi-parabolic descendent algebras $\BD^{\qpar}_{\alpha_\bullet}$; recall that this is a supercommutative algebras generated by symbols $\ch_k(\gamma)$ for $\gamma\in H^\ast(C)$ and $\ch_k(e_j)$ for $0\leq j\leq l+1$, modulo certain simple relations; alternatively, it is generated by $\ch_k^H(\gamma), \ch_k^H(e_j)$. One defines  $\bL_n^\qpar\colon \BD^\qpar_{\alpha_\bullet}\to \BD^\qpar_{\alpha_\bullet}$ by $\bL_n^\qpar=\bR_n^\qpar+\bT_n^\qpar$ where $\bR_n^\qpar$ is defined exactly in the same way as in the parabolic case, and $\bT_n$ is 
\begin{align*}\bT_n^\qpar=\sum_{a+b=n}a!b!\left((1-g)\ch_a(\pt)\ch_b(\pt)+\sum_{j=0}^l\ch_a(e_j)\ch_b(e_{j+1}-e_j)\right)\,.
\end{align*}
In particular, the Virasoro constraints do not depend on a choice of weights within a fixed chamber.
\end{remark}

The main result of this section are the Virasoro constraints for moduli spaces of parabolic bundles:

\begin{theorem}\label{thm: virasoro}
Let $\alpha=(r,d,f_\bullet, c)\in C(I)$. For every $D\in \BD^\paar$ we have
\[\int_{[M_\alpha]}\bL_\inv(D)=0\,.\]
\end{theorem}

These generalize the Virasoro constraints for moduli spaces of stable bundles shown in \cite{blm}. For the proof, we will need compatibility statements between the Virasoro constraints and the structures that we have discussed earlier in the paper: wall-crossing, flag bundles, and Hecke operators. We develop these compatibilities first, and then give the proof in Section \ref{subsec: proofvirasoro}.

\subsection{Primary states and wall-crossing compatibility}\label{subsec: wccompatibility} Compatibility with wall-crossing is already the key ingredient in \cite{blm}. To state it, let us introduce the Lie subalgebra of primary states. 

\begin{definition}\label{def: primarystate}
The space of primary states of weight $i$ on $\V^\paar_{\tr}$ is
\[P_i=\big\{v\in \V^\paar_{\tr}\colon L_0(v)=iv\textup{ and }L_n(v)=0\textup{ for }n>0\big\}\,.\]
The space of primary states on the Lie algebra $\widecheck \V^\paar_{\tr}$ is
\[\widecheck P_0=P_1/T(P_0)\subseteq \widecheck \V^\paar_{\tr}\,.\]
\end{definition}

By \cite[Corollary 5.7]{lmquivers}, a class $u\in \widecheck \V^\paar_{\tr}$ (of non-trivial topological type) is a primary state if and only if $L_\inv(u)=0$, where
\[L_\inv=\sum_{n\geq -1}\frac{(-1)^n}{(n+1)!}L_{-1}^{n+1}\circ L_n\colon \widecheck \V^\paar_{\tr}\to  \V^\paar_{\tr}\] 
is the dual of $\bL_\inv$.\footnote{When $\V^\paar_{\tr}$ admits a conformal element $\omega$, $L_\inv$ is a canonical lift of the operator $[-, \omega]$, see \cite[Proposition 3.13, Lemma 3.15]{blm}.} 
In particular, the moduli space $M_\alpha$ satisfies the Virasoro constraints (i.e. Theorem \ref{thm: virasoro} holds) if and only if $[M_\alpha]\in \widecheck P_0$ is a primary state. 

We recall the reader that $\V^\paar_{\tr}$ is a lattice vertex algebra (cf. Proposition \ref{prop: latticeva}). When the symmetrized Euler pairing $\chi_{\CP}^\sym$ is non-degenerate, this vertex algebra carries a natural conformal element, which is defined as explained in \cite[Section 3.3, Section 4.3]{blm}. In the case of vector bundles, it is shown in \cite[Theorem 4.12]{blm} that the Virasoro operators coming from this conformal element coincide with \eqref{eq: operatorsva}; the general parabolic case can be shown by a straightforward adaptation of the calculations. 

More generally, if the symmetrized Euler pairing is degenerate, the vertex algebra $\V^\paar_{\tr}$ may not include a conformal element; however, it can always be embedded into a larger (lattice) vertex algebra admitting a conformal element which induce Virasoro operators that restrict to \eqref{eq: operatorsva}. See \cite[Section 5.1.1, Proposition 5.23]{lmquivers} for a more detailed discussion.

A fundamental property of primary states, observed by Borcherds himself in his foundational paper on vertex algebras, is that they are closed under the Lie bracket.

\begin{lemma}[\cite{Borcherds}]\label{lem: liesubalgebra}
$\widecheck P_0$ is a Lie subalgebra of $\widecheck \V^\paar_{\tr}$.
\end{lemma}

Since wall-crossing formulas are stated using the Lie bracket on $\widecheck \V^\paar_{\tr}$, this statement implies a compatibility between the Virasoro constraints and wall-crossing: if every moduli space appearing on the right hand side of \eqref{eq: fullwcformula} satisfies the Virasoro constraints, then so does the moduli space on the left.

\subsection{Flag bundle compatibility}\label{subsec: flagcompatibility}
We now turn to the compatibility with the flag bundle maps between different moduli spaces of parabolic bundles that were discussed in Section \ref{sec: partial}. We recall that we have interpreted these flag bundle structures in terms of a certain map of vertex algebras $\Omega_{I, I'}\colon \V^I\to \V^{I'}$ (and the induced homomorphism of Lie algebras $\widecheck\V^I\to \widecheck\V^{I'}$), see Proposition \ref{prop: flagbundleclass} and Definition \ref{def: joyceclassesboundary}. However, the Virasoro operators \eqref{eq: operatorsva} were defined on the larger vertex algebra  $\V^I_{\tr}$. Before we prove a compatibility statement between the Virasoro operators and the vertex algebra homomorphisms $\Omega_{I, I'}$, we show that the Virasoro operators preserve~$\V^I$, which we think is a statement of independent interest.

\begin{proposition}\label{prop: operatorsdescend}
The Virasoro operators $\{L_n\}_{n\geq -1}$ preserve the vertex subalgebra $\V^\paar\subseteq \V^\paar_{\tr}$. 
\end{proposition}
\begin{proof}
This statement is equivalent to the claim that the operators $\bL_n$ preserve the ideal 
\begin{equation}\label{eq: kernelrealization}
\ker\big(\BD^\paar_\alpha\twoheadrightarrow H^\ast(\CM_\alpha^\paar)\big)\,.
\end{equation}
Since $\bT_n$ is a multiplication operator, this is equivalent to $\bR_n$ preserving the said ideal. Recall that $\CM_\alpha^\paar$ is a flag bundle over the moduli stack of vector bundles $\CM_{r,d}$. Moreover, it is well-known \cite{AB, hs} that $H^\ast(\CM_{r,d})$ is isomorphic to the free supercommutative algebra generated by $p_\ast(c_k(\CV)q^\ast \gamma)$ for $k=1, \ldots, r$ and $\gamma\in H^\ast(C)$. Hence, the kernel \eqref{eq: kernelrealization} is generated by relations of two types:
\begin{enumerate}
\item $c_k(\CF(t)-\CF(s))$ for $t>s$ and $k>f(t)-f(s)$;
\item $p_\ast(c_k(\CV)q^\ast \gamma)=0$ for every $k>r$, $\gamma\in H^\ast(C)$. 
\end{enumerate}
The relations of the first type are analogous to the relations among generators in the cohomology ring of the Grassmannian, which were shown to be preserved by similar Virasoro operators \cite[Section 8.7]{lmquivers}. See also the proof of \cite[Proposition 3.8]{klmp}. The same proof can be used to show that each of the ideals generated by the relations (1) for fixed $t,s$ is preserved by the Virasoro action.

On the other hand, the relations of the second type can be seen as generalized Mumford relations, as in \cite[Section 2.2]{klmp}. To do so, observe that 
\[\CV=(Rp_{12})_\ast R\CH om(p_{23}^\ast\CK, p_{13}^\ast \CV)\]
where $p_{ij}$ is the projection of $\CM_\alpha\times C\times C$ onto the $i$-th and $j$-th components and $\CK$ is the complex $\Delta_\ast(\omega_C^\vee)[-1]$ on $C\times C$; this follows from
\begin{align*}R\CH om(\CV, \CO_{\CM\times C})&=(Rp_{12})_\ast\Delta_\ast R\CH om(\Delta^\ast p_{13}^\ast \CV, \CO_{\CM\times C})\\
&=(Rp_{12})_\ast R\CH om( p_{13}^\ast \CV, p_{23}^\ast \Delta_\ast \CO_{C})
\end{align*}
and an application of Grothendieck--Verdier duality. By the formal argument in the proof of \cite[Theorem 3.10]{klmp}, $\bR_n$ preserves the ideal generated by relations of the second type provided that $\bR_n$ preserves the ideal of relations
\[\ker\big(\xi_{\CK}\colon \BD\twoheadrightarrow H^\ast(C)\big)\,,\]
where $\BD$ is the algebra generated by symbols $\ch_k(\gamma)$ and $\xi_{\CK}$ is the realization homomorphism sending $\ch_k(\gamma)$ to $p_\ast(\ch_k(\CK)q^\ast \gamma)$. But it is easy to see that we have $\xi_{\CK}\circ \bR_n=0$ for every $n\geq 1$; this follows from the fact that $\bR_n(\ch_0^H(\gamma))=0$, by definition, and $\xi_{\CK}\big(\ch_k^H(\gamma)\big)=0$ for every $k>1$ by Hodge degree reasons.\qedhere
\end{proof}

\begin{theorem}\label{thm: voamap}
Let $I'\subseteq I$. The Virasoro operators commute with $\Omega_{I, I'}$, i.e. the square
\begin{center}
\begin{tikzcd}
\V^I \arrow[r, "\Omega_{I, I'}"]\arrow[d, swap,"L_n^I"]& \V^{I'} \arrow[d, "L_n^{I'}"]\\
\V^I \arrow[r, "\Omega_{I, I'}"]& \V^{I'}
\end{tikzcd}
\end{center}
commutes for every $n\geq -1$. In particular, $\Omega_{I, I'}$ sends primary states to primary states. 
\end{theorem}
If the vertex algebras $\V^I, \V^{I'}$ had conformal elements inducing their respective Virasoro operators (which is never the case), this commutativity property would mean that $\Omega_{I, I'}$ was a weak vertex operator algebra homomorphism. Being a weak vertex operator algebra homomorphism is implied by being a strong vertex operator algebra homomorphism, which is a vertex algebra homomorphism respecting the conformal elements. Recall that, crucially, $\Omega_{I, I'}$ cannot be extended to a homomorphism of vertex algebras $\V^I_{\paar}\to\V^{I'}_{\paar}$ (which often contains a conformal element), so there is no natural way to produce an actual homomorphism of vertex operator algebras.

\begin{proof}[Proof of Theorem \ref{thm: voamap}]
The commutativity is equivalent to the commutativity of the dual square
\begin{center}
\begin{tikzcd}
H^\ast(\CM^{I'}) \arrow[rr, "c_{\textup{top}}\cdot \pi^\ast(-)"]\arrow[d, swap,"\bL_n^I"]& & H^\ast(\CM^{I}) \arrow[d, "\bL_n^{I'}"]\\
H^\ast(\CM^{I'}) \arrow[rr, "c_{\textup{top}}\cdot \pi^\ast(-)"]& &H^\ast(\CM^{I}) 
\end{tikzcd}
\end{center}
where $c_{\textup{top}}\coloneqq c_{\textup{top}}(\Delta^\ast \Xi_{I, I'})$. Note that we have\footnote{We are committing a small abuse of notation here by omitting the realization homomorphism $\xi$ and regarding $\ch_k(g)$ as an element of $H^\ast(\CM^\paar)$.}
\[\pi^\ast \ch_k(e(t))=\ch_k(e(t))\,\textup{ for }t\in I'\,,\textup{ and }\pi^\ast \ch_k(\gamma)=\ch_k(\gamma)\,\textup{ for }\gamma\in H^\ast(C)\,,\]
so $\pi^\ast$ intertwines $\bR_n^I$, $\bR_n^{I'}$. Hence, 
\begin{align*}
\bL_n^I\big(c_{\textup{top}}\cdot \pi^\ast(D)\big)-c_{\textup{top}}\cdot \pi^\ast(\bL_n^{I'}(D))\big)=\big(\bR_n^I(c_{\textup{top}})+c_{\textup{top}}\cdot \bT_n^I-c_{\textup{top}}\cdot \pi^\ast\bT_n^{I'}\big)\cdot \pi^\ast(D)\,.
\end{align*}
Therefore, commutativity of the diagram is equivalent to the identity 
\begin{equation}\label{eq: identityvoamap}
\bR_n^I(c_{\textup{top}})=-c_{\textup{top}}\cdot \bT_n^I+c_{\textup{top}}\cdot \pi^\ast\bT_n^{I'}=-c_{\textup{top}}\cdot\left(\sum_{a+b=n}\sum_{t\in I\setminus I'}a!b!\ch_a(e(t))\ch_b(\partial e(t))\right)
\end{equation}
in $H^\ast(\CM^I)$. By Newton's identities, the total Chern class of $\Xi$ is
\[c(\Xi)=(\xi\otimes \xi)\exp\left(\sum_{\substack{a,b\geq 0\\ (a,b)\neq (0,0)}}\sum_{t\in I\setminus I'}(-1)^{b}(a+b-1)!\ch_a(e(t))\otimes \ch_b(\partial e(t))\right)\,.\]
Denote by $C_j\in \BD^I\otimes \BD^I$ the degree $2j$ part of the exponential on the right hand side, so that $c_j(\Xi)=(\xi\otimes \xi)C_j$. It is shown in \cite[Theorem A]{klmp} that we have a formal identity
\begin{align*}(\bR_n^I\otimes \id)(C_N)+&\sum_{k=-1}^n \binom{n+1}{k+1}(\id\otimes \,\bR_k^I)(C_{N+n-k})=\\
&-\sum_{0\leq a+b\leq n}\frac{a!(n-a)!}{(n-a-b)!}\sum_{t\in I\setminus I'}\big(\ch_a(e(t))\otimes \ch_b(\partial e(t))\big) C_{N+n-a-b}\,
\end{align*}
in $\BD^I\otimes \BD^I$ for every $N$. Now take $N=\rk\, \Xi$ in the identity above and apply the realization map $\xi\otimes \xi$. Since $\Xi$ is a vector bundle, the realization of $C_{N+n-a-b}$ in $H^\ast(\CM^I\times \CM^I)$ is 0 when $a+b<n$. Similarly, the realization of $C_{N+n-k}$ is 0 if $k<n$; since by Proposition \ref{prop: operatorsdescend} the operators $\bR_k^I$ descend to $H^\ast(\CM^I)$, the realization of $(\id\otimes \bR_k^I)(C_{N+n-k})$ is also $0$ for $k<n$. Therefore, the formal identity implies that
\[(\xi\otimes \xi)(\bR_n^I\otimes \id+\id\otimes\, \bR_n^I)C_N=-(\xi\otimes \xi)\left(\sum_{a+b=n}a!b!\sum_{t\in I\setminus I'}\ch_a(e(t))\otimes \ch_b(\partial e(t))\right) C_N\,\]
in $H^\ast(\CM^I\times \CM^I)$. By further applying $\Delta^\ast$ we deduce the identity \eqref{eq: identityvoamap}, which proves the commutativity statement. The conclusion regarding primary states follows trivially.
\end{proof}

\begin{remark}
This type of statement holds in other scenarios where the Virasoro constraints have been studied and there are natural maps between vertex algebras. The projective bundle compatibility in \cite[Theorem 5.13]{blm} is an earlier and slightly simpler statement of this kind. An analogue of Theorem \ref{thm: voamap} can be proved, exactly in the same way, for the maps between vertex algebras associated to different quivers in \cite[Definition 5.10]{GJT} and the quiver Virasoro operators studied in \cite{bojkoquivers, lmquivers}. 
\end{remark}

\subsection{Compatibility with Hecke operators}
\label{subsec: heckecompatibility}
The next ingredient is compatibility of the Virasoro constraints with Hecke operators, which we prove in the next lemma.

\begin{lemma}\label{lem: heckevirasoro}
 The Virasoro operators commute with the Hecke operators $(h_\tau)_\dagger$, i.e. the square
\begin{center}
\begin{tikzcd}
\V^\paar_{\tr} \arrow[r, "(h_\tau)_\dagger"]\arrow[d, swap,"L_n"]& \V^\paar_{\tr} \arrow[d, "L_n"]\\
\V^\paar_{\tr} \arrow[r, "(h_\tau)_\dagger"]& \V^\paar_{\tr}
\end{tikzcd}
\end{center}
commutes for every $n\geq -1$. In particular, $(h_\tau)_\dagger$ sends primary states to primary states. 
\end{lemma}
\begin{proof}
Recall that the Hecke operator $(h_\tau)_\dagger$ is defined as the dual of operators $h_\tau^\dagger$ on the descendent algebra. In turn, these are induced by an automorphism $\phi_\tau\colon H_I(C)\to H_I(C)$ in the sense that $h_\tau^\dagger(\ch_k^H(g))=\ch_k^H(\phi_\tau(g))$. As explained in the proof of Proposition \ref{prop: latticeva}, $\V^\paar_{\tr}$ is isomorphic to a lattice vertex algebra whose underlying ``lattice'' is the dual of $H_I(C)$, which can be identified with $K(I)\oplus K^1(C)$. The dual of the automorphism $\phi_\tau$ is easily seen to be
\[h_\tau\oplus \textup{id}_{K^1(C)}\colon K(I)\oplus K^1(C)\to K(I)\oplus K^1(C)\,.\]
Therefore, the automorphism $(h_\tau)_\dagger$ of $\V^\paar_{\tr}$ is induced by the above automorphism of the underlying lattice. For the lattice vertex algebras without an odd part, it was shown in \cite[Proposition 5.2]{lmquivers} that Virasoro operator are compatible with maps of vertex algebras induced by embeddings (in particular, automorphisms) of lattices. In the more general case of lattices with an odd part, the construction of Virasoro operators requires a further choice of an isotropic splitting of the odd part, see \cite[3.6.14]{Ka98} or the discussion in \cite[Section 3.3]{blm}. The result mentioned above \cite[Proposition 5.2]{lmquivers} is immediately generalized to the case of super lattices if we consider embedding of lattices that preserve the $\BZ/2$-grading of the lattice and the isotropic decomposition mentioned above.

In our case, the isotropic splitting of the odd part of the lattice is
\[K^1(C)\simeq H^1(C)= H^{0,1}(C)\oplus H^{1,0}(C)\,.\]
Since our lattice automorphism acts as the identity on the odd part, it clearly preserves this splitting, and we conclude the proof.\qedhere
\end{proof}

\begin{corollary}\label{cor: heckecompatibility}
Let $\alpha\in C(I)$ be regular. Theorem \ref{thm: virasoro} holds for $M_\alpha$ if and only if it does for $M_{h_\tau\alpha}$. 
\end{corollary}
\begin{proof}
By Lemma 5.5, we have $(h_\tau)_\ast[M_\alpha]=[M_{h_\tau\alpha}]$ in the Lie algebra $\widecheck \V^\paar$. Embedding $\widecheck \V^\paar\subseteq \widecheck \V^\paar_{\tr}$ and applying Lemma \ref{lem: hecketr} we obtain the equality $(h_\tau)_\dagger[M_\alpha]=[M_{h_\tau\alpha}]$ in $\widecheck\V^\paar_{\tr}$. The conclusion then follows from Lemma \ref{lem: heckevirasoro}.
\end{proof}

\subsection{Compatibility with the Weyl action}
\label{subsec: weylcompatibility}
This very short section will consist only of the proof of the following easy proposition:

\begin{proposition}\label{prop: weylcompatibility}
Let $\alpha_\bullet$ be a full topological type. The operators $\bL_n^\qpar\colon \BD^\qpar_{\alpha_\bullet}\to \BD^\qpar_{\alpha_\bullet}$ are equivariant with respect to the $\Sigma_r$ action.
\end{proposition}
\begin{proof}
The fact that $\bR_n^\qpar$ is $\Sigma_r$ equivariant is clear. The $\Sigma_r$ equivariance of $\bT_n^\qpar$ is also clear once we write it as
\[\bT_n^\qpar=\sum_{a+b=n}a!b!\left((1-g)\ch_a(\pt)\ch_b(\pt)+\frac{1}{2}\sum_{1\leq i\neq j\leq r} \ch_a(d_i)\ch_b(d_j)\right)\,.\qedhere
\]

\end{proof}

\subsection{Proof of Theorem \ref{thm: virasoro}}
\label{subsec: proofvirasoro}
We now have all the ingredients for the proof of the Virasoro constraints for parabolic bundles. Our strategy will be exactly the same as in the proof of the reconstruction theorem (cf. Section \ref{subsec: proofreconstruction}) and the Newstead vanishing (cf. Section \ref{subsec: proofnewstead}), by induction on the rank. The base case are the Virasoro constraints for the Jacobian; these can be shown by a direct calculation, or alternatively they are shown in \cite[Proposition 6.2, Corollary 6.4]{blm}.

As in the previous proofs, the case of rank 2 requires utilizing 2 punctures, as sketched in Section \ref{subsubsec: rank2punctures}, but otherwise the argument is similar, so we focus on $r>2$. We first observe that it is enough to prove the Virasoro constraints for $M_{r,1}^\full(c)$ for any $c\in S_{r,1}$. First, the statement for regular full parabolic bundles in any other degree $d$, i.e. $[M_{r,d}^\full(c)]\in \widecheck P_0$ for every $c\in S_{r,d}$ regular, follows from the compatibility with Hecke operators, cf. Corollary \ref{cor: heckecompatibility}. By compatibility with wall-crossing, cf. Proposition \ref{lem: liesubalgebra}, and the induction hypothesis, we get $[M_{r,d}^\full(c)]\in \widecheck P_0$ for every $c\in \overline S_{r,d}$. Finally, the statement for partial parabolic bundles is deduced from Proposition \ref{prop: flagbundleclass} (or, in the case of non-regular $c$, Definition \ref{def: joyceclassesboundary}) and Theorem \ref{thm: voamap}. 

We now claim that the class $L_\inv([M_{r,1}^\full(c)])$ does not depend on $c$ after forgetting weights. Given any two stability conditions $c, c'\in S_{r,1}$, recall the wall-crossing formula in Theorem \ref{thm: fullwc}. Every term with $l>1$ on the right hand side is an iterated bracket of classes $[M_{r_i, d_i}^\full(c_{J_i})]$ with $r_i<r$. By the induction hypothesis and the fact that primary states form a Lie subalgebra (cf. Lemma \ref{lem: liesubalgebra}) it follows that any such bracket is annihilated by $L_\inv$. Hence, applying $L_\inv$ to the wall-crossing formula gives an equality
\[L_\inv([M_{r,1}^\full(c)])=L_\inv([M_{r,1}^\full(c')])\]
after forgetting weights. Equivalently, this is to say that the functional $\BD^\qpar_{\alpha_\bullet}\to \BQ$ defined by
\[D\mapsto \int_{[M_{r,1}^\full(c)]}\bL_\inv^\qpar(D)\]
does not depend on $c$. This functional is anti-invariant with respect to the action of $\Sigma_r$ for some choice of $c$ (recall Proposition \ref{prop: weylcompatibility}) and anti-invariant with respect to the action of $\widetilde \Sigma_r$ for some other choice of $c''$. Hence, it must vanish by Proposition~\ref{prop: affineweylvanishing}.

\bibliographystyle{mybstfile.bst}
\bibliography{refs.bib} 
\end{document}